\documentclass[12pt]{amsbook} 
\usepackage{amssymb,latexsym,amsmath,amscd} 

 \usepackage{chngcntr}

 \counterwithout{section}{chapter} 

\swapnumbers 

\numberwithin{equation}{section}

\input cyracc.def

\newtheorem{theorem}{Theorem}[section] 
 
\newtheorem{thm}[theorem]{Theorem} 
\newtheorem{comment}[theorem]{Comment}
\newtheorem{construction}[theorem]{Construction} 
\newtheorem{constructions}[theorem]{Constructions} 
\newtheorem{construdef}[theorem]{Construction--Definition} 
\newtheorem{conv}[theorem]{Convention} 
\newtheorem{corollary}[theorem]{Corollary} 
\newtheorem{cory}[theorem]{Corollary} 
\newtheorem{definition}[theorem]{Definition} 
\newtheorem{df}[theorem]{Definition} 
\newtheorem{example}[theorem]{Example} 
\newtheorem{examples}[theorem]{Examples} 

\newtheorem{lemma}[theorem]{Lemma}
\newtheorem{notation}[theorem]{Notation} 
\newtheorem{prop}[theorem]{Proposition} 
\newtheorem{propcon}[theorem]{Proposition-Construction} 
\newtheorem{proposition}[theorem]{Proposition}
\newtheorem{remark}[theorem]{Remark} 
\newtheorem{remarks}[theorem]{Remarks} 
\newtheorem{theordef}[theorem]{Theorem--Definition}

\newcommand\theoref{Theorem~\ref}
\newcommand\lemref{Lemma~\ref}
\newcommand\propref{Proposition~\ref}
\newcommand\corref{Corollary~\ref}
\newcommand\defref{Definition~\ref}
\newcommand\remref{Remark~\ref}

\newcommand\exref{Example~\ref}

\newcommand\comref{Comment~\ref}
\newcommand\secref{Section~\ref}
\newcommand\chapref{Chapter~\ref}

\def \p{{\po \m \it Proof. }} 
\def \mp{\po \m} 
\def\po{\parindent 0pt} 
 
\def \m{\medskip} 
\def \sm{\smallskip} 
\def \bs{\bigskip} 
 
\def \eps{\varepsilon} 
\def \lp{ $\boxplus$} 
 
\def \Wh{\operatorname {Wh}} 
\def \pt{\operatorname {pt}} 
 
\def \rel{\operatorname {rel}} 
\def \Hom{\operatorname {Hom}} 
\def \Ext{\operatorname {Ext}} 
\def \Ker{\operatorname {Ker}} 
\def \IM{\operatorname {Im}} 
 
\def \Lift{\operatorname {Lift}} 
\def \aut{\operatorname {aut}} 
\def \ev{\operatorname {ev}} 
\def \ad{\operatorname {ad}}
\def \Wh{\operatorname {Wh}}

\def \Num{\operatorname {Num}}
\def \Nor{\operatorname {Nor}}
\def \bs{\operatorname {bs}}

\def \tors{\operatorname {tors}}

\def \sirc{{\raise0.2ex \hbox{$\scriptstyle \circ$}}} 
\def \suset{{\raise0.2ex \hbox{$\subset$}}} 
\def \equity{{\raise0.2ex \hbox{=}}} 
\def \wt{\widetilde} 
\def \wh{\widehat} 
\def \ov{\overline} 
\def\la{\langle}
\def\ra{\rangle}
\def\ts{\times}
\def\pa{\partial}
\def\lr{\xrightarrow{\hspace*{17pt}}}
 
\def\ga{\alpha} 
\def\gb{\beta} 
\def\gga{\gamma} 
\def\gd{\delta} 
 
\def\gf{\varphi} 
\def\gk{\kappa} 
\def\gkk{\varkappa} 
\def\gl{\lambda} 
 
\def\gs{\sigma}

\def\gS{\Sigma} 
 
\def\CC{\mathbb C} 
\def\NN{\mathbb N} 
\def\PP{\mathbb P} 
\def\QQ{\mathbb Q} 
\def\RR{\mathbb R} 
\def\ZZ{\mathbb Z}
\def\cp2{\CC \PP^2}

\def\II{\mathfrak I}
\def\HH{\mathcal H} 
 
\def\C #1{\mathcal{#1}} 
\def\spl{{\C S}_{PL}} 
 
\def\tpl{{\C T}_{PL}} 
\def\td{{\C T}_{D}} 
\def\sd{{\C S}_{D}} 
\def\pd{{\C P}_{D}}

\def\jtop{j_{TOP}} 
\def\jf{j_{F}} 
 
\def\bul{^{\scriptscriptstyle {\bullet}}} 
 
\def\lp{$\diamond$} 
 
\long\def\forget#1\forgotten{} %
 
\begin{document}
 

\title[PL structures]{Piecewise Linear Structures on Topological 
Manifolds}

\date{\today} 
 
\maketitle 

{\topmargin -1in \textheight 10in \tableofcontents }

\chapter*{Preface}
 

In his paper Novikov \cite[p.409]{N3}  wrote: 
 
\begin{quote} Sullivan's Hauptvermutung theorem was 
announced first in early 1967. After the careful analysis made by Bill 
Browder and myself in Princeton, the first version in May 1967 (before 
publication), his theorem was corrected: a necessary restriction on 
the 2-torsion of the group $H_3(M)$ was missing. This gap was found 
and restriction was added. Full proof of this theory has never been 
written and published. Indeed, nobody knows whether it has been 
finished or not. Who knows whether it is complete or not? This 
question is not clarified properly in the literature. Many pieces of 
this theory were developed by other topologists later. In particular, 
the final Kirby--Siebenmann classification of topological 
multi-dimensional manifolds therefore is not proved yet in the 
literature. 
\end{quote}

\m 
I do not want to discuss here whether the situation is so dramatic as 
Novikov wrote. However, it is definitely true that, up to now, there is no 
detailed enough and well-ordered exposition of Kirby--Siebenmann 
classification, such that it can be recommended to advanced                
students which are willing to learn the subject. The fundamental book 
of Kirby--Siebenmann \cite{KS2} was written by pioneers and, in a sense, 
in hot pursuit. It contains all the necessary ingredient for the proof, but it is really 
``Essays'', and one have to do a certain work in order to make it easy readable for 
general audience. 
 
\m 

\m {\bf Acknowledgments.} The work was partially supported by Max-Planck of Mathematics, Bonn,  
and by a grant from the Simons Foundation (\#209424 to Yuli Rudyak).
I express my best thanks to Andrew Ranicki who read the
whole manuscript and did many useful remarks and comments. I am also grateful
to Hans-Joachim Baues for useful discussions.     

\newpage

\chapter*{Notation and Conventions}

\m We work mainly with $CW$-spaces and topological manifolds. However, when we 
quit these classes by taking products or functional spaces, we equip the last 
ones with the compactly generated topology, (following Steenrod 
\cite{Steenrod} and McCord \cite{McCord}, see e.g.\cite{Rud} for the 
exposition). All maps are supposed to be continuous. All neighborhoods 
are supposed to be open. 
 
 \m We denote the one-point space by $\pt$.
 
 \m A {\it pointed space} is pair $(X,\{x_0\}$ where $x_0$ is a point of $X$. We also use that notation $(X,x_0)$ and call $x_0$ the {\it base point} of $X$. If we do not need to indicate the base point, we can write $(X,*)$ (or even $X$ if it is clear that $X$ is a pointed space). Given two pointed spaces $(X,x_0)$ and $(Y,y_0)$, a pointed map is a map $f: X \to Y$ such that $f(x_0)=y_0$.  
  
\m Given two topological spaces $X,Y$, we denote by $[X,Y]$ the set of 
homotopy classes of maps $X \to Y$. We also use the notation $[X,Y]\bul$ for 
the set of pointed homotopy classes of pointed maps $X \to Y$ of pointed 
spaces.  
 
\m It is quite standard to denote by $[f]$ the homotopy class of a map $f$. 
However, frequently we do not distinguish a map and its homotopy class and 
use the same symbol, say $f$ for a map as well as for the homotopy class. In 
this paper this does not lead to any confusion.   
 
\m We use the term {\it inessential map} for null-homotopic maps; otherwise a map is called {\it essential}.  
 
\m We use the sign $\simeq$ for homotopy of maps or homotopy equivalence 
of spaces. We use the sign $\cong$ for bijection of sets or isomorphism of groups. 
We use the notation $:=$ for ``is defined to be''.

\m We reserve the term {\it bundle} for locally trivial bundles and 
the term {\it fibration} for Hurewicz fibrations.

  Given a space $F$, an  
{\it $F$-bundle} is a bundle whose fibers are homeomorphic to $F$, and an {\it 
$F$-fibration} is a fibration whose fibers are homotopy equivalent to $F$. 

We denote the {\it trivial} $F$-bundle $X\ts F\to X$ over $X$ by $\theta_X^F$ or merely $\theta^F$, Also, we denote the trivial $\RR^n$-bundle over $X$ by $\theta^n$ or $\theta^n$.

\m {\footnotesize We do not mention {\it microbundles} at all, because in the topological and in the PL category every $n$-dimensional microbundle over a space $X$ contains an $\RR^n$-bundle over $X$, and these bundles are unique up to equivalence, see Kister~\cite{Kis} for the topological category and Kuiper-Lashof~\cite{KL} for the PL category. For this reason, any claim on microbundles can be restated in terms of bundles. The reader should keep it in the mind when we cite (quote about) something concerning microbundles.  } 
 
\m Given a bundle or fibration $\xi =\{p: E \to B\}$, the space $B$ is called the 
{\it base} of $\xi$ and denote also by $\bs(\xi)$, i.e.  $\bs(xi)=B$. The space $E$ is called the {\it total space} of $\xi$.  
Furthermore, given a space $X$, we set   
$$ 
\xi \times X =\{p\times 1: E \times X \to B \times X\}.  
$$ 
 
\m Given two bundles $\xi=\{p: E \to B\}$ and $\eta=\{q: Y \to 
X\}$, a {\it bundle morphism} $\gf: \xi \to \eta$ is a commutative diagram 
$$ 
\CD 
E @>g>> Y\\ 
@VpVV @VVqV\\ 
B @>f>> \,X. 
\endCD 
$$\ 
We say that $f$ is the {\it base} of the morphism $\gf$ or that $\gf$ is a 
{\it morphism over $f$}. We also say that $g$ {\it is a map over} $f$.   
If $X=B$ and $f=1_B$ we say that $g$ is a  map over $B$ (and $\gf$ is a 
morphism over $B$).  
 
\m Given a map $f: Z \to B$ and a bundle (or fibration) $\xi=\{p: E \to B\}$, we use 
the notation $f^*\xi$ for the induced bundle over $Z$. Recall that $f^*(\xi)=\{r:D\to Z\}$ where 
\[
D=\{(z,e)\in Z\ts E\, \bigm| f(z)=p(e)\}
\]. 
There is a canonical bundle morphism  
\[
I=\mathfrak I_{f}=\mathfrak I_{f, \xi} :f^*\xi \to \xi 
\]
given by the map $ D\to E, (z,e)\mapsto e$
 over $f$, see \cite{Rud} (or \cite{FR} where it is denoted by $\ad(f)$). Following 
\cite{FR}, we call $\mathfrak I_{f,\xi}$ the {\it adjoint morphism of $f$}, or 
just the {\it $f$-adjoint} morphism. Furthermore, given a bundle morphism $\gf: 
\xi \to \eta$ with the base $f$, there  exists a unique bundle morphism 
$c(\gf): \xi \to f^*\eta$ over the base of $\xi$ such that the composition 
$$ 
\CD 
\xi @>c(\gf)>> f^*\eta @> \mathfrak I_{f, \eta} >> \eta 
\endCD 
$$ 
coincides with $\gf$. Following \cite{FR}, we call $c(\gf)$ the {\it 
correcting morphism}. 
 
\m Given a subspace $A$ of a space $X$ and a bundle $\xi$ over $X$, we denote by 
$\xi|_A$ the bundle $i^*\xi$ where $i: A \subset X$ is the inclusion.  

\m Given a map $p: E \to B$ and a map $f: X \to B$,  a {\it $p$-lifting of $f$} is any map $g: X\to E$ with $pg=f$. Two $p$-lifting $g_0, g_1$ of $f$ are {\it vertically homotopic} if there exists a homotopy $G: X\ts I \to E$ between $g_0$ and $g_1$ such that $pg_t=f$ for all $t\in I$. The set of vertically homotopic $p$-liftings of $f$ is denoted by $[\Lift_pf]$. 

\m We denote by $p_k, w_k$ and $L_k$ the Pontryagin,  Stiefel--Whitney, and 
Hirzebruch characteristic classes, respectively. We denote by $\gs(M)$ the 
signature of a manifold $M$. See~\cite{MS} for the definitions.

\newpage

\chapter*{Introduction} 
 
Throughout the paper we use abbreviation PL for ``piecewise linear''.  
 
\m 
{\it Hauptvermutung} (main conjecture) is an abbreviation for  
{\it die Hauptvermutung der kombinatorischen Topologie}  
(the main conjecture of combinatorial topology). It seems that the conjecture 
was first formulated in the papers of Steinitz \cite{Steinitz} and Tietze 
\cite{Tietze} in 1908. 
 
\m The conjecture claims that the topology of a simplicial complex determines 
completely its combinatorial structure. In other words, two simplicial 
complexes are simplicially isomorphic whenever they are homeomorphic. This 
conjecture was disproved by Milnor~\cite{Mi2} in 1961.  
 
\m However, for manifolds one can state a refined version of the {\it 
Hauptvermutung} by considering  simplicial 
complexes with additional restrictions.  A {\it PL manifold} is defined to be a simplicial complex 
such that the star of every point (the union of all closed 
simplexes containing the point) is simplicially isomorphic to the 
$n$-dimensional ball. Such simplicial complexes are also called {\it combinatorial triangulations}.
Equivalently, a PL manifold can also be defined a manifold equipped with a   
maximal PL atlas.

There exist  topological manifolds that are homeomorphic to a simplicial 
complex but do not admit a PL structure (non-combinatorial triangulations), see \exref{tri-nonpl}. 
Furthermore, there exist topological manifolds that are not homeomorphic to any simplicial complex, see \exref{non-tri}.
 
\m  Now, the {\it Hauptvermutung for manifolds} asks if any two 
homeomorphic PL manifolds are PL homeomorphic. Furthermore, the related 
question asks whether every topological manifold is homeomorphic to a PL 
manifold. Both these questions were solved (negatively) by Kirby and 
Siebenmann \cite{KS1, KS2}. In fact, Kirby and Siebenmann classified PL 
structures on high-dimensional topological manifolds. It turned out that a 
topological manifold can have different PL structures, as well as not to have any. 
Now we give a brief description of these results.  
 
\m Let $BTOP$ and $BPL$ be the classifying spaces for 
stable topological and PL bundles, respectively. We regard the forgetful map 
$\ga: BPL \to BTOP$ as a fibration and denote its homotopy fiber by $TOP/PL$.  
 
\m Let $f: M \to BTOP$ classify the stable tangent bundle of a topological 
manifold $M$. By main properties of classifying spaces, every PL structure on $M$ gives us a $\ga$-lifting 
of $f$ and that every two such liftings for the same PL structure are fiberwise homotopic. 
 
\m It is remarkable that the inverse is also true provided that $\dim 
M\geqslant 5$.  
In greater detail, $M$ admits a PL structure if $f$ admits a $\ga$-lifting (the 
Existence Theorem \ref{exist}), and PL structures on $M$ are in a bijective 
correspondence with fiberwise homotopy classes of $\ga$-liftings of $f$ (the 
Classification Theorem \ref{classif}). Kirby and Siebenmann proved these 
theorems and, moreover, they proved the following Main Theorem:
\begin{equation*}
TOP/PL \textbf{  is the Eilenberg--MacLane 
space }K(\ZZ/2,3). 
\end{equation*}

Thus, there is only one obstruction  
$$ 
\varkappa (M)\in 
H^4(M;\ZZ/2) 
$$  
to an $\ga$-lifting of $f$, and the set of fiberwise homotopic 
$\ga$-liftings of $f$ (if they exist) is in bijective correspondence with 
$H^3(M;\ZZ/2)$. In other words, a topological manifold $M,\, \dim 
M\geqslant 5$ admits a PL structure if and only if $\varkappa (M) =0$. 
Furthermore, every homeomorphism $h: V \to M$ of two PL manifolds
assigns a class  
$$ 
\varkappa (h) \in H^3(M;\ZZ/2), 
$$  
and $\varkappa (h)=0$ if and only if $h$ is concordant to a PL homeomorphism (or, equivalently, to the 
identity map $1_M$, see \remref{conc}(2)). Finally, every class $a\in H^3(M;\ZZ/2)$ has the form 
$a=\varkappa(h)$ for some homeomorphism $h: V\to M$ of two PL manifolds. 

These results give us the complete 
classification of PL structures on a topological manifold of dimension $\geqslant 5$. In particular, the situation with Hauptvermutung turns out to be understandable. See \secref{ksks} for more detailed exposition
 
\m We must explain the following. It can happen that two different PL 
structures on $M$ yield PL homemorphic PL manifolds (like that two $p$-liftings 
$f_1, f_2: M \to BPL$ of $f$ can be non-fiberwise homotopic). Indeed, roughly 
speaking, a PL structure on a topological manifold $M$ is a concordance class 
of PL atlases on $M$ (see Section \ref{structures} for accurate definitions). 
However, a PL automorphism of a PL manifold can turn the atlas into a 
non-concordant to the original one, see Example \ref{non-conc}. So, in fact, 
the  set of pairwise non-isomorphic PL manifolds which are homeomorphic to a 
given PL manifold is in a bijective correspondence with the set 
$H^3(M;\ZZ/2)/R$ where $R$ is the following equivalence relation: two PL 
structure are equivalent if the corresponding PL manifolds are PL homeomorphic. 
The {\it Hauptvermutung} for manifolds claims that the set $H^3(M;\ZZ/2)/R$ is 
a singleton for all $M$. But this is wrong in general. 
 
\m Namely, there exists a PL manifold $M$ which is homeomorphic but not PL 
isomorphic to $\RR\PP^{n},n\geqslant 5$, see \exref{hptv}. So, here we have a 
counterexample to the {\it Hauptvermutung}. 
 
\m To complete the picture, we mention again that there are 
topological manifolds that do not admit any PL structure, see 
\exref{non-triang}. Moreover, there are manifold that cannot be triangulated as simplicial complexes, see \exref{non-tri}.

\m Comparing the classes of smooth, PL and topological 
manifolds, we see that there is a big difference between first and second 
classes, and not so big difference between second and third ones. From the 
homotopy-theoretical point of view, one can say that the space $PL/O$ (which 
classifies smooth structures on PL manifold, see \remref{pl/o}) has many 
non-trivial homotopy groups, while the space $TOP/PL$ is an Eilenberg--MacLane space. 
Geometrically, one can mention that there are many smooth manifolds which are 
PL homeomorphic to $S^n$ but pairwise non-diffeomorphic, while any PL manifold 
$M^n, n \geqslant 5$ is PL homeomorphic to $S^n$ provided that it is homeomorphic to 
$S^n$.  
 
\m It is worthwhile to go one step deeper and explain the following. Let $M^{4k}$ 
be a closed connected almost parallelizable manifold (i.e. $M$ becomes 
parallelizable after deletion of a point). Let $\gs_k$ denote the minimal 
natural number which can be realized as the signature of the manifold 
$M^{4k}$. In fact, for every $k$ we have three numbers $\gs_k^S, 
\gs_k^{PL}$ and $\gs_k^{TOP}$ while $M^{4k}$ is a smooth, PL or 
topological manifold, respectively. Milnor and Kervaire \cite{MK} proved that   
$$   
\gs_k^S = c_k(2k-1)!  
$$  
where $c_k \in \NN$. On the other hand,  
$$  
\gs_1^{PL}=16 \text { and } \gs_k^{PL} =8 \text{ for } k>1. 
$$ 
Finally, 
$$ 
\gs_k^{TOP} =8 \text{ for all } k. 
$$ 
 
\m So, here we can see again the big difference between smooth and PL cases. 
On the other hand, $\gs^{PL}_k=\gs^{TOP}_k$ for $k>1$. Moreover, we will see 
below that the number  
$$ 
2=16/8=\gs^{PL}_1/\gs^{TOP}_1 
$$ 
is another guise of the number
\[
2=\text{ the order of the group }\pi_3(TOP/PL).
\]  
 
\m In this context, it makes sense to notice about low dimensional manifolds, because of the 
following remarkable contrast. There is no difference between PL and smooth 
manifolds in dimension $<7$: every PL manifold $V^n, n<7$ admits a unique 
smooth structure. However, there are infinitely many smooth manifolds which 
are homeomorphic to $\RR^4$ but pairwise non-diffeomorphic, see \secref{s:ex}, Summary.

\m Concerning the description of the homotopy type of 
$TOP/PL$, we have the following. Because of the Classification Theorem, if $k+n \geqslant 5$ then the group 
$\pi_n(TOP/PL)$ is in a bijective correspondence with the set of PL structures 
on $\RR^k \times S^n$. However, this set of PL
structures looks wild and uncontrollable. In order to make the situation more
manageable, we consider PL structures on the {\it compact} manifold $T^k\ts S^n$ and then
extract the necessary information on the universal covering $\RR^k \times S^n$ from here. We can't do it directly,
but there is a trick (the Reduction Theorem \ref{t:reduction}) which allows us
to estimate PL structures on $\RR^k \times S^n$ in terms of so-called {\it
homotopy PL structures} on $T^k \times S^n$ (more precisely, we should
consider the homotopy PL structures on $T^k \times D^n$ modulo the boundary),
see Section~\ref{structures} for the definitions. Now, using results of Hsiang
and Shaneson~\cite{HS} or Wall~\cite{Wall2,Wall3} about homotopy PL structures
on $T^k \times D^n$, one can prove that $\pi_i(TOP/PL)=0$ for $i\ne 3$ and that
$\pi_3(TOP/PL)$ has at most 2 elements. Finally, there exists a
high-dimensional topological manifold which does not admit any PL structure.
Hence, by the Existence Theorem, the space $TOP/PL$ is not contractible.
Thus, $TOP/PL = K(\ZZ/2,3)$.

\m For better arrangement of the previous paragraph, look at the graph located at the end of the Introduction. 
Here we formulate without proofs the  boxed claims (and provide the necessary preliminaries and references), 
while in \chapref{arch} we explain how a claim (box) can be deduced from other ones, accordingly with the arrows in the graph.

\m Let me tell you something more about the graph. As we have already seen, the 
classification theory of PL structures on topological manifolds splits into two 
parts. The first one reduces the original geometric problem to a homotopy one 
(a classification of $p$-liftings of a map $M \to BTOP$ to $BPL$), the second part 
solves this homotopy problem by proving that $TOP/PL=K(\ZZ/2,3)$.  
 
\m The Product Structure Theorem 
\ref{product} is a very important ingredient for the proof. Roughly speaking, this theorem establishes a bijection between 
PL structures on $M$ and $M \times \RR$. The Classification Theorem 
\ref{classif} and the Existence Theorem \ref{exist} are the consequences 
of the Product Structure Theorem. 
 
\m  
Now I say some words about the top box of the above graph. Let $F_n$ be the  
monoid of pointed homotopy equivalences $S^n \to S^n$, let $BF_n$ be the  
classifying space for $F_n$, and let $BF=\lim_{n \to \infty} BF_n$.   There is
an obvious forgetful map $BPL \to BF$, and we denote by $F/PL$   the homotopy
fiber of this map. For every homotopy equivalence of closed   PL manifolds 
$h: V \to M$ Sullivan \cite{Sul1, Sul2} defined the {\it   normal invariant} of
$h$ to be a certain homotopy class $j_F(h)\in    [M,F/PL] $, see Section
\ref{from}. 

Let $M, \dim M\geqslant 5$ be a closed PL manifold.
Sullivan proved that, for  every {\it  homeomorphism} $h: V \to M$,
we have $j_{F}(h)=0$ whenever $H_3(M)$ is  2-torsion free. 
Moreover, this theorem implies that if, in addition, $M$ is simply-connected 
then $h$ is  homotopic to a PL homeomorphism. Thus the {\it Hauptvermutung} 
holds for simply-connected closed manifolds $M, \dim M \geqslant 5$
 with $H_3(M)$ 2-torsion free, see \secref{secnorm}.   

\m Definitely, the above formulated Sullivan Theorem on the Normal  Invariant of a Homeomorphism is
interesting by itself. However, in the paper  on hand this theorem plays also
an additional important role. Namely, the  Sullivan Theorem for $T^k \times
S^n$ is a lemma in classifying of homotopy  structures on $T^k \times D^n$.
For this reason, we first prove the Sullivan  Theorem for $T^k \times S^n$,
then use it in the proof of the Main Theorem,  and then (in \chapref{c:norminv}) use the Maim Theorem in order to prove
the Sullivan Theorem in full generality.    

\m You can also see that the proof of the Main Theorem uses the  difficult  Freedman's example of a 4-dimensional
almost parallelizable topological  manifold of signature 8. This example
provides the equality   $\gs^{TOP}_1=8$. Actually, the original proof of the
Main Theorem appeared  before Freedman's Theorem and therefore did not use the
last one. However, as we already mentioned, the Freedman's results clarify the
relations between PL and  topological manifolds, and thus they should be
incorporated in the exposition  of the global picture.      

\m  The paper is organized as follows. The first chapter contains the  
architecture of the proof of the Main Theorem: $TOP/PL \simeq K(\ZZ/2,3)$. 
In fact, here we comment the above mentioned graph.

\m The second chapter contains a proof of the Sullivan Theorem on the normal 
invariant of a homeomorphism for $T^k \times S^n$, i.e. we attend the top box of the graph. 
We also discuss  the Browder--Novikov Theorem \ref{t:bn} about homotopy properties of  normal
bundles: we need this discussion in order to clarify the concept of normal invariant.  
 
\m The third chapter contains some applications if the Main Theorem. We complete the proof of the Sullivan Theorem on the normal 
invariant of a homeomorphism and tell more on classification of PL manifolds an, in particular, on Hauptvermutung. Several interesting examples are considered. Finally, we discuss the homotopy and topological invariants of certain characteristic classes.

\vfill
\pagebreak

\vskip -3cm
                        
\section*{The Graph}


 
 
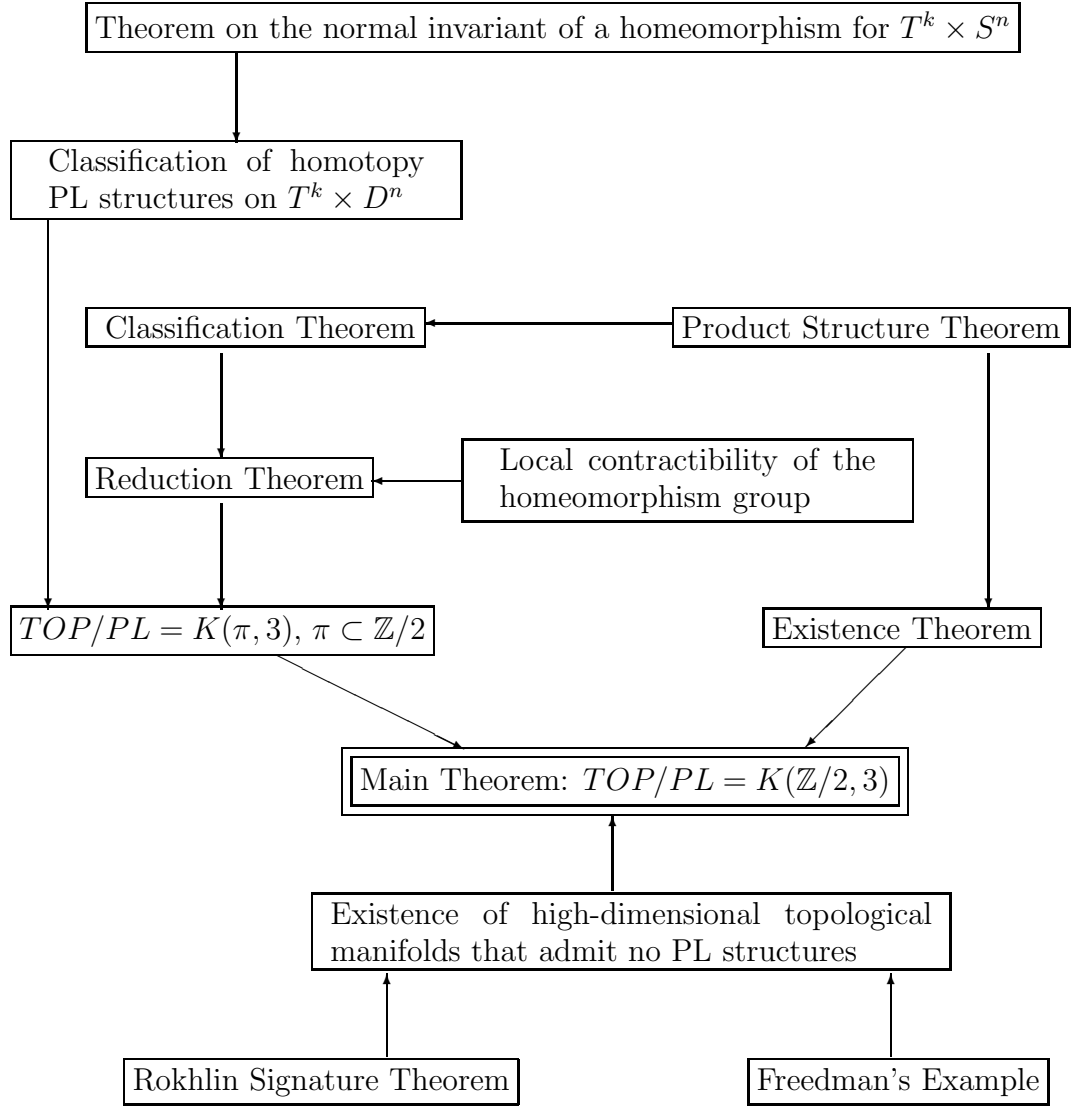
\unitlength=1mm
 
\begin{picture}(180,190)  

\put(10,170){\fbox{Theorem on the normal invariant of a homeomorphism for $T^k \times S^n$}} 
 
\put(0,150){\framebox[60mm]{\parbox{50mm}{Classification of homotopy PL 
structures on $T^k \times D^n$}}} 

\put(10,130){\fbox{ Classification Theorem}} 

\put(88,130){\fbox{Product Structure Theorem}} 
 
\put(10,110){\fbox{Reduction Theorem}} 
 
\put(60,110){\framebox[60mm]{\parbox{50mm}{Local contractibility of the homeomorphism group}}} 
 
\put(0,90){\fbox{$TOP/PL = K(\pi,3)$, $\pi \subset \ZZ/2$}} 

\put(100,90){\fbox{Existence Theorem}} 
 
\put(44,70){\fbox{\fbox{Main Theorem: $TOP/PL=K(\ZZ/2,3)$}}} 
 
\put(40,50){\framebox[85mm]{\parbox{80mm}{Existence of high-dimensional topological manifolds that admit no PL structures}}} 
 
\put(15,30){\fbox{Rokhlin Signature Theorem}} 
 
\put(98,30){\fbox{Freedman's Example}} 
 
\put(30,168){\vector(0,-1){12}} 

\put(88,132){\vector(-1,0){33}} 
 
\put(28,128){\vector(0,-1){14}} 
 
\put(28,108){\vector(0,-1){14}} 
 
\put(60,111){\vector(-1,0){12}} 
 
\put(5,146){\vector(0,-1){52}} 
 
\put(130,128){\vector(0,-1){34}} 
 
\put(119,88.9){\vector(-1,-1){13.3}} 
 
\put(35.25,87.85){\vector(2,-1){25}} 
 
\put(80,56.5){\vector(0,1){10}} 
 
\put(50,34.5){\vector(0,1){11}} 
 
\put(117,34.5){\vector(0,1){11}}

\end{picture}

\chapter{ Architecture of the Proof}\label{arch}

\section{Principal Fibrations} 
 
 Recall that an $H$-space is a space $F$ with a base point $f_0$ and a 
multiplication map $\mu: F \times F \to F$ such that $f_0$ is a homotopy unit, 
i.e. the maps $f\mapsto \mu(f,f_0)$ and  $f\mapsto \mu(f_0,f)$ are homotopic 
to the identity rel $\{f_0\}$. For details, see \cite{BV}. 
 
\begin{definition}\rm\label{d:prin} 
(a) Let $(F,f_0)$ be an $H$-space with the multiplication $\mu: F \times F 
\to F$. A {\it principal $F$-fibration} is an $F$-fibration $p: E \to B$ 
equipped with a map $m: E \times F \to E$ such that the following holds:  
\par 
(i) the diagrams  
$$  
\CD 
E\times F \times F @>m\times 1>> E \times F @. \phantom{OOO} @. E\times F @>m>> E\\ 
@V1\times \mu VV @VVmV @. @Vp_1VV @VVpV\\ 
E \times F @>m>> E @.  \phantom{OOO}  @. E @>p>> B 
\endCD 
$$ 
commute; 
\par(ii) the map 
$$ 
E \longrightarrow E, \quad e\mapsto m(e,f_0) 
$$ 
is a homotopy equivalence; 
\par (iii) for every $e_0\in E$, the map 
$$ 
F \longrightarrow p^{-1}(p(e_0)), \quad f\mapsto m(e_0, f) 
$$ 
is a homotopy equivalence. 

{\mp (b)} A {\it trivial} principal $F$-fibration is the fibration $p_2: 
X\times F \to F$ with the action $m: E\times F \to E$ of the form 
$$ 
m: X \times F \times F \to X \times F, \quad m(x, f_1, f_2)=(x, 
\mu(f_1,f_2)). 
$$ 
\end{definition} 
 
\m It is easy to see that if the fibration $\eta$ is induced from a principal 
fibration $\xi$ then $\eta$ turns into a principal fibration in a canonical 
way.  
\begin{definition} \rm 
Let $\pi_1: E_1 \to B$ and $\pi_2: E_2 \to B$ be two principal $F$-fibrations 
over the same base $B$. We say that a map $h: E_1 \to E_2$ is an {\it 
$F$-equivariant map over $B$} if $h$ is a map over $B$ and the diagram  
$$  
\CD 
E_1\times F  @>h\times 1>> E_2 \times F\\ 
@Vm_1VV @VVm_2V\\ 
E_1 @>h>> E_2 
\endCD 
$$ 
commutes up to homotopy over $B$. 
\end{definition} 
 
\m Note that, for every $b\in B$, the map  
$$ 
h_b: \pi_1^{-1}(b) \to  \pi_2^{-1}(b), \quad h_b(x)=h(x) 
$$  
is a homotopy equivalence. 
 
\m Now, let $p: E \to B$ be a principal $F$-fibration, and let $f: X \to 
B$ be an arbitrary map. Given a $p$-lifting $g: X \to E$ of $f$ and a 
map $u: X \to F$, consider the map 
$$ 
\CD 
g_u: X @>\Delta>> X \times X @>g\times u>> E \times F @>m>>E. 
\endCD 
$$ 
It is easy to see that the correspondence $(g,u) \mapsto g_u$ yields a 
well-defined map (right action) 
\begin{equation}\label{action} 
[\Lift_p f] \times [X,F] \to [\Lift_p f] .
\end{equation} 
In particular, for every $p$-lifting $g$ of $f$ the correspondence $u 
\mapsto g_u$ induces a map 
$$ 
T_g : [X,F] \to [\Lift_p f]. 
$$  
 
\begin{theorem} \label{lift} 
Let $\xi=\{p: E \to B\}$ be a principal $F$-fibration, and let $f: X 
\to B$ be a map where $X$ is assumed to be paracompact and locally 
contractible. If $F$ is a homotopy associative $H$-space with a homotopy inversion, then $[X,F]$ is a group the 
above action $\eqref{action}$ is free and transitive provided 
$[\Lift_p f] \ne \emptyset$. In particular, for every $p$-lifting $g: 
X \to E$ of $f$ the map $T_g$ is a bijection. 
\end{theorem} 
 
\p We start with the following lemma. 
\begin{lemma}\label{l1} 
The theorem holds if $X=B$, $f=1_X$ and $\xi$ is the trivial principal 
$F$-fibration. \end{lemma} 
 
\p In this case every $p$-lifting $g: X \to X \times F$ of 
$f=1_X$ determines and is completely determined by the map 
$$ 
\CD 
\ov g: X: @>g>> X \times F @>p_2>> F. 
\endCD 
$$ 
In other words, we have the bijection $[\Lift_p f] \cong [X,F]$, and 
under this bijection the action \eqref{action} turns into the 
multiplication 
$$ 
[X,F] \times [X,F] \to [X,F]. 
$$ 
Now the result follows since $[X,F]$ is a group. 
\lp 
 
\m We finish the proof of the theorem. Consider the induced fibration $f^*\xi 
=\{q: Y \to X\}$ and note that there is an $[X,F]$-equivariant bijection 
\begin{equation}\label{bi-lift} 
 [\Lift_p f] \cong [\Lift_q1_X]. 
\end{equation} 
Now, suppose that $[\Lift_p f]\ne \emptyset$ and take a $p$-lifting $g$ of 
$f$. Regarding $Y$ as the subset of $X\times E$, define the $F$-equivariant 
map  
$$  
h: X \times F \to Y, \quad h(x,a)=(x,g(x)a),\, x\in X, a\in F. 
$$  
It is easy to see that the diagram 
$$ 
\CD 
X \times F @>h>> Y\\ 
@Vp_1VV @VV qV\\ 
X @= X 
\endCD 
$$ 
commutes, i.e. $h$ is a map over $X$. Since $X$ is a locally contractible 
paracompact space, and by a theorem  of Dold \cite{Dold}, there exists a map 
$k: Y \to X \times F$ over $X$ which is homotopy inverse over $X$ to $h$. It 
is easy to see that $k$ is an equivariant map over $X$. Indeed, if $m_1: X 
\times F \times F \to X \times F$ and $m_2: Y \times F \to Y$ are the 
corresponding actions then 
$$ 
m_1(k\times 1) \simeq khm_1(k\times 1) \simeq 
km_2(h\times 1)(k\times 1) = \simeq km_2(hk\times 1) \simeq km_2, 
$$ 
where 
$\simeq$ denotes the homotopy over $X$. 
 
\m In particular, there is an $[X,F]$-equivariant bijection  
$$ 
 [\Lift_q1_X] \cong [\Lift_{p_1}1_X] 
$$ 
where $p_1: X \times F \to X$ is the projection. 
Now we compose this bijection with \eqref{bi-lift} and get $[X,F]$-equivariant 
bijections  
$$ 
[\Lift_p f]] \cong [\Lift_q1_X] \cong [\Lift_{p_1}1_X], 
$$ 
and the result follows from Lemma \ref{l1}. 
\qed 
 
\begin{example}\rm\label{loops} 
If $p: E \to B$ is an $F$-fibration then $\Omega p: \Omega E \to 
\Omega B$ is a principal $\Omega F$-fibration. Here $\Omega$ denotes 
the loop functor. 
\end{example}

\section{Preliminaries on Classifying Spaces} 
 
Here we give a brief recollection on $\RR^n$ bundles, spherical 
fibrations, and their classifying spaces. For details, see \cite[Chapter IV]{Rud}. 
 
\begin{df}\label{d:morph} \rm
We define a {\it topological $\RR^n$-bundle} over a space $B$ to be an 
$\RR^n$-bundle $p: E \to B$ equipped with a fixed section $s: B \to 
E$ (the zero section). Given two topological $\RR^n$-bundles $\xi=\{p:E\to B\}$ and 
$\eta=\{q: Y\to X\}$, we define a {\it topological $\RR^n$-morphism} $\gf: \xi \to 
\eta$ to be a commutative diagram 
\begin{equation}\label{morph} 
\CD 
E @>g>> Y\\ 
@VpVV @VVqV\\ 
B @>f>> X 
\endCD 
\end{equation} 
where $g$ preserves the sections and induces a homeomorphism on each of fibers. 
The last one means that, for every $b\in B$, the map 
$$ 
g_b: \RR^n = p^{-1}(b) \to q^{-1}(f(b))=\RR^n, \quad g_b(a)=g(a) 
\text{ 
for all } a\in p^{-1}(b) 
$$ 
is a homeomorphism. As usual, we call $f$ the {\it base} of the morphism 
$\gf$. and denote it also $\bs(\gf)$, i.e. $\bs(\gf)=f$. 

We say that  topological $\RR^n$-morphism $\gf$ is a morphism over $B$ if the map $f$ in \eqref{morph} is equal to $1_B$ 
 
 A topological $\RR^n$-morphism is a topological  {\it $\RR^n$-isomorphism} if the above mention $g$ is a homeomorphism. 
 
 We define two topological $\RR^n$-morphisms $\gf_0, \gf_1:\xi\to \eta$ to be 
 {\it bundle homotopic} if  there exists a topological $\RR^n$-morphism $\Phi: \xi \times 
I \to \eta$ such that $\Phi|_{\xi} \times \{i\}=\gf_i, i=0,1$. 

A topological $\RR^n$-morphism $\gf: \xi \to \eta$  is a {\it bundle homotopy equivalence} if there exists a topological $\RR^n$-morphism $\psi: \eta to \xi$ such that $\gf\psi$ and $\psi\gf$ are bundle homotopic to the corresponding identity maps.

Frequently, we will say merely ``homotopy'' instead of ``bundle homotopy'', etc. if this does not lead to confusions.
\end{df}

\begin{theordef}\label{univ} 
There exists a topological $\RR^n$-bundle $\gga^n_{TOP}$ with the following {\em universal propert}: For every topological $\RR^n$-bundle $\xi$ over a $CW$-space $B$,  every $CW$-subspace $A$ of $B$ and every morphism  
$$ 
\psi: \xi_A \to \gga^n_{TOP} 
$$  
of topological $\RR^n$-bundles, there exists a morphism $\gf:\xi \to 
\gga^n_{TOP}$ which is an extension of $\psi$. The base of $\gga^n_{TOP}$ is called the {\em classifying space} for topological 
$\RR^n$-bundles ans denoted by $BTOP_n$. 
\end{theordef} 
 
 \m 
We can regard topological $\RR^n$-bundles as $(TOP_n, \RR^n)$-bundles, 
i.e. $\RR^n$-bundles with the structure group $TOP_n$. 
Here $TOP_n$ is the topological group of self-homeomorphism $f: \RR^n 
\to \RR^n, f(0)=0$. The classifying space $BTOP_n$ of the group 
$TOP_n$ turns out to be a classifying space for topological 
$\RR^n$-bundles.
 
 \forget
 \m It is known that $BTOP_n$ can be chosen to be a CW space, and in this case $BTOP_n$ is determined uniquely up to homotopy equivalence.
\forgotten
  
Consider a topological $\RR^n$-bundle $\xi$ over 
a CW space $B$. By the definition of universal bundle,  there exists a  topological 
$\RR^n$-morphism $\gf: \xi \to \gga^n_{TOP}$. We call such $\gf$ a {\it classifying morphism} for 
$\xi$. The base $f: B\to BTOP_n$ of $\gf$ is called a {\it  
classifying map} for $\xi$. It is clear that $\xi$ is  
isomorphic over $B$ to $f^*\gga^n_{TOP}$.  
 
\begin{proposition}\label{clas-hom} 
If $\gf_0, \gf_1: \xi \to \gga^n_{TOP}$ be two classifying morphisms 
for $\xi$, then there are homotopic. 
In particular, a classifying map $f$ for $\xi$ is determined by $\xi$  
uniquely up to homotopy.  
\end{proposition} 
 
\p This follows from the universal property \ref{univ} applied to $\xi \times 
I$, if we put $A=X\times \{0,1\}$ where $X$ denotes the 
base of $\xi$. 
\qed

\m
\begin{remark}\rm
This is important to understand the difference between classifying maps and classifying morphisms. Non-homotopic classifying morphisms can induce homotopic classifying maps. On the other hand, not every map $X \to BTOP_n$ is a classifying map, while every morphism $\xi \to \gga^N$ is a classifying morphism. To feel the difference, consider a trivial $\RR^n$-bundle $\theta^N$ over a space $X$. Then there is a classifying map $X\to*\in BTOP_n)$, while a corresponding morphism is a trivialization of $\theta^N$, i.e. a morphism $X\times R^n \to R^n$. 
\end{remark}
 
\m 
A {\it piecewise linear} (in future PL) {\it$\RR^n$-bundle} is a topological 
$\RR^n$-bundle $\xi=\{p: E \to B\}$ such that $E$ and $B$ are 
polyhedra and $p: E \to B$ and $s: B \to E$ are PL maps. Furthermore, 
we require that, for every simplex $\Delta \subset B$, there is a PL 
homeomorphism $h: p^{-1}(\Delta) \cong \Delta \times \RR^n$ with 
$h(s(\Delta))=\Delta \times \{0\}$. (For definitions of PL maps, see 
\cite{Hudson, RS}.) 
 
A PL morphism of PL $\RR^n$-bundles is a topological 
$\RR^n$-morphism where the maps $g$ and $f$ in \eqref{morph} are 
PL maps. 
 
There exists a universal PL $\RR^n$-bundle $\gga^n_{PL}$ over a 
certain space $BPL_n$. This means that the universal property 
\ref{univ} remains valid if we replace $\gga^n_{TOP}$ by 
$\gga^n_{PL}$ and ``topological $\RR^n$ bundle'' by ``PL 
$\RR^n$-bundle'' there. So, $BPL_n$ is a classifying space for PL 
$\RR^n$-bundles. 

This is worthy to mention that $BPL_n$ can be chosen to be a locally finite simplicial complex,~\cite[Essay IV \S8]{KS2}.
 
   Note that $BPL_n$ can also be regarded as the classifying space of a 
certain group $PL_n$ (which is constructed as the geometric realization of a 
certain simplicial group), \cite{KL, LR}. 
 
\m 
A {\it sectioned $S^n$-fibration} is defined to be an $S^n$-fibration $p: E 
\to B$ equipped with a section $s: B\to E$. Morphisms of 
sectioned $S^n$-fibrations are defined similarly  to
\defref{d:morph} where each map $g_b$ and the total map $g$ is assumed to be a pointed proper homotopy 
equivalence.  We shall use the brief term ``$(S^n,*)$-fibration'' 
for sectioned $S^n$-fibrations and ``$(S^n,*)$-morphism'' for morphisms of 
sectioned $S^n$-fibrations.
 
There exists a universal sectioned $S^n$-fibration $\gga^n_{F}$. To define it, replace $\gga^n_{TOP}$ by  
$\gga^n_{F}$ and ``topological $\RR^n$ bundle'' by ``sectioned  
$S^n$-fibration'' in \ref{univ}. The base $BF_n$ of $\gga^n_{F}$ is called the 
classifying space for sectioned  $S^n$-fibrations.
The space $BF_n$ can also be regarded as the classifying space for the 
monoid $F_n$ of pointed homotopy equivalences $(S^n,*) \to (S^n,*)$. 

\forget
Furthermore, we will also use the term 
``$F_n$-morphism'' for morphism of sectioned $S^n$-fibrations.  Finally, an 
$F_n$-morphism over a space is called an $F_n$-equivalence. 
\forgotten
 
\m We need also to recall the space $BO_n$ which classifies $n$-dimensional 
vector bundles. This well-known space is described in many sources, e.g. 
\cite{MS}. The universal vector bundle over $BO_n$ is denoted by 
$\gga^n_O$. 
              
\m It is worthy to notice that the spaces $BO_n, BPL_n, BTOP_n$ and $BF_n$ are defined 
uniquely up to weak homotopy equivalence. 

\m We regard $\gga^n_{PL}$ as the (underlying) 
topological $\RR^n$-bundle and get the classifying morphism  
\begin{equation}\label{omegapl} 
\omega=\omega^{PL}_{TOP}(n):\gga^n_{PL} \to \gga^n_{TOP}. 
\end{equation} 
 We denote by $\ga=\ga^{PL}_{TOP}(n): BPL_n \to BTOP_n$ the 
base of this morphism. 
 
\m 
Given a topological $\RR^n$-bundle $\xi=\{p: E \to B \}$, 
let $\xi\bul$ denote the $S^n$-bundle 
\begin{equation}\label{xibul}
\xi\bul=\{p\bul: E\bul \to B\}
\end{equation}
 where $E\bul$ is the fiberwise one-point compactification of $E$. 
Note that the added points (``infinities'') give us a certain 
section of $\xi\bul$. 
 
\m 
{\footnotesize In other words, the $TOP_n$-action on $\RR^n$ extends uniquely 
to a $TOP_n$-action on the one-point compactification $S^n$ of 
$\RR^n$, and $\xi\bul$ is the $(TOP_n, S^n)$-bundle associated with 
$\xi$. Furthermore, the fixed point $\infty$ of the $TOP_n$-action on 
$S^n$ yields a section of $\xi\bul$.} 
\rm 
 
\m So, $\xi\bul$ can be regarded as an $(S^n,*)$-fibration over $B$. 
In particular, $(\gga^n_{TOP})\bul$ can be regarded as an 
$(S^n,*)$-fibration over $BTOP_n$. So, there is a classifying morphism  
$$ 
\omega^{TOP}_{F}(n):\gga^n_{TOP} \to \gga^n_{F}. 
$$ 
We denote by $\ga^{TOP}_{F}(n): BTOP_n \to BF_n$ the base of 
$\omega^{TOP}_{F}(n)$. 
 
\m Finally, we note that an $n$-dimensional vector bundle over a polyhedron 
$X$ has a canonical structure of PL $\RR^n$-bundle over $X$. Similarly to 
above, this gives us  a (forgetful) map  
$$  
\ga^O_{PL}(n): BO_n \to BPL_n. 
$$ 
 
\m So, we have a sequence of forgetful maps 
\begin{equation}\label{alpha(n)} 
BO_n \xrightarrow{\ga'} BPL_n \xrightarrow{\ga''} BTOP_n \xrightarrow{\ga'''}BF_n 
\end{equation} 
where $\ga'=\ga^O_{PL}(n)$, etc. 
 
\begin{constructions}\rm 
1. Given an $F$-bundle $\xi=\{p: E \to B\}$ and an $F'$-bundle  
$\xi'=\{p': E' \to B'\}$, we define the product $\xi\times \xi'$ to be  
the $F \times F'$-bundle  
$$ 
p\times p': E \times E' \to B \times B'. 
$$ 
 
\m 2. Given an $F$-bundle $\xi=\{p: E \to B\}$ with a 
section $s: B \to E$ and an $F'$-bundle $\xi'=\{p': E' \to B'\}$ with 
a section $s': B' \to E'$, we define the smash product 
$\xi\wedge\xi'$ to be the $F \wedge F'$-bundle as follows. The map 
$p\times p': E \times E' \to B \times B'$ passes through the quotient 
map $q: E \times E' \to E \times E'/(E\times s(B') \cup E'\times s(B)$, 
and we set 
$$ 
\xi \wedge \eta =\{\pi: E \times E'/(E\times s(B') \cup E'\times s(B)  
\to B \times B', 
$$ 
where $\pi$ is the unique map with $p\times p'=\pi q$. Finally, the section $s$ 
and $s'$ yield an obvious section of $\pi$. 
 
\m 3. Given an $\RR^m$-bundle $\xi$ and an $\RR^n$-bundle 
 $\eta$ over the same space $X$, the {\it Whitney sum} 
 of $\xi$ and $\eta$ is the $\RR^{m+n}$-bundle $\xi \oplus 
\eta=d^*(\xi \times \eta)$ where $d: X \to X \times X $ is the 
diagonal. 
 
Note that if $\xi$ and $\eta$ are a PL $\RR^m$ and PL 
$\RR^n$-bundle, respectively, then $\xi \times \eta$ is a PL $\RR^{m+n}$-bundle.  
 
\m 4. Given a sectioned  $S^m$-bundle $\xi$ and sectioned 
$S^n$-bundle $\eta$ over the same space $X$, we set $\xi \dagger 
\eta=d^*(\xi \wedge \eta)$. 
\end{constructions}

We denote by $r_n=r_n^{TOP}: BTOP_n \to BTOP_{n+1}$ the map which classifies 
$\gga^n_{TOP} \oplus  \theta^1_{BTOP_n}$. The maps $r_n^{PL}: BPL_n \to 
BPL_{n+1}$ and  $r_n^{O}: BO_n \to BO_{n+1}$ are defined in a similar way. 
 
\m 
We can also regard the above map $r_n: BTOP_n \to BTOP_{n+1}$ as a map 
induced by the standard inclusion $TOP_n\subset TOP_{n+1}$. Using this  
approach, we define $r_n^F: BF_n \to BF_{n+1}$ as the map induced by the 
standard inclusion $F_n\subset F_{n+1}$, see \cite[p. 45]{MM}. 
 
\begin{remarks} \rm 
1. Regarding $\RR^m$ as the bundle over the point, we see that 
$(\RR^m)\bul=(S^m)$ and, moreover, 
$$ 
(\RR^m \times \RR^n)\bul = S^m\wedge S^n, \text{ i.e. } (\RR^m \oplus 
\RR^n)\bul = S^m\dagger S^n. 
$$ 
Therefore $(\xi \oplus\eta)\bul=\xi\bul \dagger \eta\bul$ for every 
$\RR^m$-bundle $\xi$ and $\RR^n$-bundle $\eta$. 
 
\m 2. Generally, the smash product of (sectioned) fibrations is not a 
fibrations. But we apply it to bundles only and so do not have any 
troubles. On the other hand, there is an operation $\wedge ^h$, the 
{\it homotopy smash product}, such that $\xi \wedge^h \eta$ is the $(F 
\wedge G)$-fibration over $X \times Y$ if $\xi$ is an $F$-fibration 
over $X$ and $\eta$ is an $G$-fibration over $Y$, see \cite{Rud}. In 
particular, one can use it in order to define an analog of Whitney sum 
for spherical fibrations and then use this one in order to construct the 
map $BF_n \to BF_{n+1}$. 
\end{remarks} 
 
\m Now we consider the classifying spaces $BO_n, BPL_n, BTOP_n$ and $BF_n$ as $n\to {\infty}$. In greater
detail,  we do the following. 
 
\m Choose classifying spaces $B'F_n$ for $(S^n,*)$-fibrations (i.e., in the 
weak homotopy type $BF_n$) and consider the maps $r_n^F: B'F_n \to B'F_{n+1}$ 
as above. We can assume that every $B'F_n$ is a $CW$-complex and every $r_n$ 
is a cellular map.  We define $BF$ to be the telescope (homotopy direct limit) 
of the sequence  
$$  
\CD 
\cdots @>>>  B'F_n @>r_n >> B'F_{n+1} @>>> \cdots \,, 
\endCD 
$$ 
see e.g. \cite[Definition I.3.19]{Rud}. Furthermore, we define $BF_n$ to be the telescope of the 
finite sequence $$ 
\CD 
\cdots @>>>  B'F_{n-1} @>r_{n-1} >> B'F_{n}. 
\endCD 
$$ 
(Note that $BF_n \simeq B'F_n$.)  So, we have the sequence (filtration) 
$$ 
\CD 
\cdots \suset \,BF_n \suset \,BF_{n+1} \suset \cdots \,. 
\endCD 
$$ 
So, $BF=\bigcup BF_n$ and $BF_n$ is closed in $BF$. Moreover, $BF$ has the direct 
limit topology with respect to the filtration $\{BF_n\}$. Furthermore, if $f: K 
\to BF$ is a map of a compact space $K$ then there exists $n$ such that 
$f(K)\suset BF_n$. 
 
\m Now, for every $n$ consider a $CW$-space $B'TOP_n$ in the weak homotopy type 
$BTOP_n$ and define $B''TOP$ to be the telescope of the sequence  
$$ 
\CD 
\cdots @>>>  B'TOP_n @>r_n >> B'TOP_{n+1} @>>> \cdots \,. 
\endCD 
$$ 
 Furthermore, we define $B''TOP_n$ to be the telescope of the finite sequence 
$$ 
\CD 
\cdots @>>>  B'TOP_{n-1} @>r_{n-1} >> B'TOP_{n}. 
\endCD 
$$ 
So, we have the diagram  
\begin{equation}\label{eq:auxdiag}
\CD 
\cdots \, @.\suset @.  \ B''TOP_n \ @. \suset @. \ B''TOP_{n+1}\  @. \suset 
\cdots \suset @. \, B''TOP\\  
@. @. @VVV @. @VVV  @. @VV p V \\ 
\cdots \,  @.\suset @. BF_n @. \suset @. BF_{n+1} @. \suset \cdots \suset @. BF 
\endCD 
\end{equation} 
where the map $p$ is induced by maps $\ga^{TOP}_F(n)$. 
Now we apply the Serre construction and replace every vertical map in the diagram \eqref{eq:auxdiag} by its fibrational substitute. 
Namely, we set   
$$ 
BTOP=\{(x,\omega)\bigm | x\in B''TOP, \,\omega \in (BF)^I,\, p(x)=\omega(0)\} 
$$ 
and define $\ga^{TOP}_F: BTOP \to BF$ by setting 
$\ga^{TOP}_F(x,\omega)=\omega(1)$. Finally, we set   
$$ 
BTOP_n=\{(x,\omega)\in BTOP \bigm|x\in B''TOP_n,\, \omega \in (BF_n)^I \suset 
(B''TOP)^I\} 
$$ 
and get the commutative diagram 
\begin{equation*} 
\CD 
\cdots \, @.\suset @.  \ BTOP_n \ @. \suset @. \ BTOP_{n+1}\  @. \suset \cdots \suset @. \, BTOP\\ 
@. @. @VVV @. @VVV  @. @VV p V \\ 
\cdots \,  @.\suset @. BF_n @. \suset @. BF_{n+1} @. \suset \cdots \suset @. BF 
\endCD 
\end{equation*} 
where all the vertical maps are fibrations. 
 
\m Now it is clear how to proceed and get the diagram 
\begin{equation}\label{alpha} 
\CD 
\cdots \, @. \suset @. BO_n @. \suset @. BO_{n+1} @. \suset \cdots \suset  @. BO\\ 
 @. @. @VVV @. @VVV @. @VV \ga^O_{PL} V \\ 
\cdots \, @. \suset @. BPL_n @. \suset @. BPL_{n+1} @. \suset \cdots \suset @. BPL\\  
 @. @. @VVV @. @VVV  @. @VV \ga^{PL}_{TOP} V \\ 
\cdots \, @.\suset @.  \ BTOP_n \ @. \suset @. \ BTOP_{n+1}\   @. \suset \cdots \suset @. \, BTOP\\ 
@. @. @V \ga^{TOP}_F(n)VV @. @VVV  @. @VV \ga^{TOP}_{F}V \\ 
\cdots \,  @.\suset @. BF_n @. \suset @. BF_{n+1} @. \suset \cdots \suset @. BF 
\endCD 
\end{equation} 
where all the vertical maps are fibrations. Moreover, each of limit spaces has the direct limit topology with 
respect to the corresponding filtration, and every compact subspace of, say,
$BO$  is contained in some $BO_n$.    

\begin{conv}\rm
Let $\xi$ classify a map $f_n: X \to BF_n$ (say). It is convenient for us to speak about $n=\infty$ and write that a map $f: X\to BF$ classify $\xi$ if $f$ can be expressed as 
\[
\CD
f:X@>f_n>> BF_n @>\subset >>BF.
\endCD
\]
\end{conv}

\m  Take a point $b\in BTOP$, put $(TOP/PL)_b:=(\ga^{PL}_{TOP})^{-1}(b)$ to  be the fiber of $\ga=\ga^{PL}_{TOP}$, and put 
\[
\gb=\gb_b: (TOP/PL)_b\to BPL
\]
 to be  the inclusion of the fiber. In further we allow us to omit the subscript $b$ and write the fibration $\ga$ as
\begin{equation}\label{albet}
\CD
TOP/PL@>\beta^{PL}_{TOP}>>BPL@>\ga^{PL}_{TOP}>> BTOP.
\endCD
\end{equation}
This will not lead to confusions because, if we choose another point $b'\in BTOP$, then the maps $\gb_b$ and $\gb_{b'}$ occur to be homotopy equivalent. We also use the notation $TOP/PL$ for the homotopy fiber of the map $\ga: BPL \to BTOP$.  

The homotopy fiber of $\ga^O_{PL} : BO\to BPL$ is denoted by $PL/O$, the fiber 
of $\ga^{F}_{TOP}$ is denoted by $F/TOP$, etc. Similarly, the homotopy fiber of the 
composition, say, 
 $$ 
\ga^{PL}_{F}:=\ga^{TOP}_{F}\sirc\ga^{PL}_{TOP} : BPL \to BF 
$$ 
is denoted by $F/PL$. In particular, we have a fibration 
\begin{equation}\label{a-b}
\CD 
TOP/PL @>a>>  F/PL @>b>>   F/TOP. 
\endCD 
\end{equation} 
Finally, note that $F/TOP = \bigcup F_n/TOP_n$ where $F_n/TOP_n$ denotes the fiber of 
the fibration $BTOP_n \to BF_n$, and $F/TOP$ has the direct limit topology 
with respect to the filtration $\{F_n/TOP_n\}$. The same holds for other 
``homogeneous spaces'' $F/PL, TOP/PL$, etc.  
 
Because of well-known results of Milnor \cite{Mi1}, all these 
``homogeneous spaces'' have the homotopy type of $CW$-spaces. 
Furthermore, all the spaces $BO, BPL, BTOP, BF, F/PL, TOP/PL$, etc. 
are infinite loop spaces, and the maps like in \eqref{albet}   \eqref{a-b} are infinite loop maps, see~\cite{BV}. In particular, the classifying spaces $BO$, etc. are homotopy associative and invertible $H$-spaces, and the fibrations \eqref{albet}, \eqref{a-b}, etc. are principal fibrations. 
 
\m We mention also the following useful fact. 

\begin{lemma}\label{stability} 
Let $Z$ denote on of the symbols $O, PL, F$. The above described  
 map $BZ_n \to BZ_{n+1}$ induces an isomorphism of homotopy groups in 
 dimensions $\leqslant n-1$ and an epimorphism in dimension $n$.  
\end{lemma} 
 
\p For $Z=O$ and $Z=F$ it is well known, see e.g. \cite{Br}, for $Z=PL$ 
it can be found in \cite{HW}. \qed 

\begin{remark}\rm
An analog of \lemref{stability} holds for $TOP$ as well, see \remref{r:stabletop}.
\end{remark} 
                                   
\begin{remark}\label{g} 
\rm Let $G_n$ denote the topological monoid of homotopy self-equivalences 
$S^{n-1} \to S^{n-1}$. Then the classifying space $BG_n$ of $G_n$ 
classifies $S^{n-1}$-fibrations (non-sectioned). Every $h\in TOP_n$ induces a map 
$\RR^n\setminus\{0\} \to \RR^n\setminus\{0\}$ which, in turn, yields a 
self-map  
$$ 
\pi_h: S^{n-1} \to S^{n-1}, \quad \pi_h(x)=h(x)/||h(x)||. 
$$  
So, we have a map $TOP_n \to G_n$ which, in turn, induces a map  
$$ 
BTOP_n \longrightarrow BG_n 
$$  
of classifying spaces. In the language of bundles, 
this map converts a topological $\RR^n$-bundle into a (non-sectioned) spherical fibration via 
deletion of the section. 
 
We can also consider the space $BG$ by tending $n$ to $\infty$. In 
particular, we have the spaces $G/PL$ and $G/TOP$. 
 
There is an obvious forgetful map $F_n \to G_{n+1}$ (ignore sections), and it turns out that 
the induced map $BF \to BG$ (as $n\to \infty$) is a homotopy 
equivalence, see e.g. \cite[Chapter 3]{MM}. In particular, $F/PL 
\simeq G/PL$ and $F/TOP \simeq G/TOP$. 
\end{remark}

\section{Structures on Manifolds and Bundles}\label{structures} 
 
A {\it PL atlas} on a topological manifold is an atlas such that all the
transition  maps are PL ones. We define a PL manifold as a topological
manifold with a  maximal PL atlas. Furthermore, given two PL manifolds $M$ and
$N$, we say that  a homeomorphism $h: M \to N$ a PL homeomorphism if $h$ is a PL
map. (One can  prove that in this case $h^{-1}$ is a PL map as well,
\cite{Hudson}.)   

\m 
\begin{definition} \label{str} 
\rm (a) We define a $\partial_{PL}$-manifold to be a topological 
manifold whose boundary $\partial M$ is a PL manifold. In particular,  
every closed topological manifold is a $\partial_{PL}$-manifold.  
Furthermore, every PL manifold can be canonically regarded as a  
$\partial_{PL}$-manifold.  
 
(b) Let $M$ be a $\partial_{PL}$-manifold. A {\it PL structuralization} 
of $M$ is a homeomorphism $h: V \to M$ such that $V$ is a PL manifold 
and $h$ induces a PL homeomorphism $\partial 
V \to \partial M$ of boundaries (or, equivalently, PL homeomorphism of 
corresponding collars). Two PL structuralizations $h_i: V_i \to M, 
i=0,1$ are {\it concordant} if there exist a PL homeomorphism $\gf: V_0 
\to V_1$ and a homeomorphism $H: V_0\times I \to M \times I$ such that 
$H|_{V\times \{0\}}=h_0$ and $H|_{V_0\times \{1\}}=h_1\gf$ and, 
moreover, $H: \partial V_0 \times I \to \partial M \times I$ coincides 
with $h_0\times 1_I$. Any concordance class of PL structuralizations 
is called a {\it PL structure on $M$}. We denote by  $\tpl(M)$ the set 
of all PL structures on $M$. 
 
(c) If $M$ on its own is a PL manifold then $\tpl(M)$ contains the 
distinguished element: the concordance class of $1_M$. We call it the 
{\it trivial} element of $\tpl(M)$. 
\end{definition} 
 
\begin{remarks}\label{conc} 
\rm 1. Clearly, every PL structuralization of $M$ 
equips $M$ with a certain PL atlas. Conversely, if we equip $M$ with a 
certain PL atlas then the identity map can be regarded as a PL 
structuralization of $M$. 
 
 2. If $M$ by itself is a PL manifold then the concordance class of 
any PL homeomorphism $h: V \to M$ is the trivial element of $\tpl (M)$.  Indeed, to prove this, we must find a homeomorphism $H: V \ts I \to M\ts I$ and a PL homeomorphism  $\gf: V \to M$ such that $H|_{V\times \{0\}}=h$ and $H|_{V\times \{1\}}=1_M\gf=\gf$.  But this is easy: put $\gf=h$ and $H(v,t)=(h(v),t)$.
 
3. Recall that two homeomorphism $h_0,h_1: X \to Y$ 
are {\it isotopic} if there exists a homeomorphism $H: X \times I \to 
Y \times I$ (isotopy) such that $p_2H: X \times I \to Y \times I \to 
I$ coincides with $p_2: X \times I \to I$. Given $A\subset X$, we say 
that $h_0$ and $h_1$ are isotopic $\rel A$ if there exists an isotopy 
$H$ such that $H(a,t)=h_0(a)$ for every $a\in A$ and every $t\in I$. 
In particular, if two PL structuralization $h_0,h_1: V \to M$ are 
isotopic $\rel\ \partial V$ then they are concordant. 
 
4. Given two PL structuralizations $h_i: V_i \to M, i=0,1$, they are 
not necessarily concordant if $V_0$ and $V_1$ are PL homeomorphic. We 
are not able to give such examples here, but we do it later, see 
Remark \ref{smooth}(2) and Example \ref{non-conc}. 
\end{remarks} 
 
\begin{definition} [cf. \cite{Br, Rud}]
\label{d:str} \rm
Given a topological $\RR^n$-bundle $\xi$, define a {\it PL structuralization} of $\xi$ to be a  topological 
$\RR^n$-morphism $\gf: \xi \to \gga^n_{PL}$. We define a {\it PL structure on $\xi$} to be a homotopy class of PL structuralizations of $\xi$.
\end{definition} 
 
Let $f_n: X \to BTOP_n$ classify a topological $\RR^n$-bundle $\xi$, and assume that there is an $\ga^{PL}_{TOP}(n)$-lifting
$$ 
g_n: X \to BPL_n 
$$  
of $f_n$. Take the $g_n$-adjoint classifying morphism 
\[
\mathfrak I:\mathfrak I_{g_n} : g_n^*\gga^n\to\gga^n.
\]
and consider the morphism
\begin{equation*} 
\xi \cong  f_n^*\gga^n_{TOP} =g_n^*\ga(n)^*\gga^n_{TOP}=g_n^*\gga^n_{PL} 
\xrightarrow{\mathfrak I} \gamma^n_{PL}, \quad \ga(n):=\ga^{PL}_{TOP}(n).  
\end{equation*}  
This morphism $\xi \to \gga^n_{PL}$ is a PL 
structuralization of $\xi$. It is easy to see that  in this way we have a 
correspondence   
\begin{equation}\label{corr} 
[\Lift_{\ga(n)} f_n] \longrightarrow \{ \text{PL structures on } \xi \}. 
\end{equation} 

\begin{theorem} 
The correspondence $\eqref{corr}$ is a bijection. 
\end{theorem} 
 
\p See~\cite[Theorem IV.2.3]{Rud}, cf. 
also~\cite[Chapter II]{Br}. 
\qed 
 
\m Consider now the map  
$$ 
\CD 
 f: X @>f_n >> BTOP_n \ \subset \ BTOP  
\endCD 
$$ 
and the map $\ga=\ga^{PL}_{TOP}: BPL \to BTOP$ as in \eqref{alpha}. Every 
$\ga(n)$-lifting $g_n: X\to BPL_n$ of $f_n$ gives us the $\ga$-lifting
\[
\CD
X @>g_n>> BPL_n @>>> BPL
\endCD
\] 
of $f$. So, we have a correspondence  
\begin{equation}\label{tpl-stab} 
u_{\xi}: \{ \text{PL structures on } \xi \}\longrightarrow [\Lift_{\ga(n)} 
f_n]\longrightarrow [\Lift_{\ga}f]  
\end{equation} 
where the first map is the inverse to \eqref{corr}. Furthermore, there is a 
canonical map   
\begin{equation*} 
v_{\xi}:\{ \text{PL structures on } \xi \}\longrightarrow 
 \{ \text{PL structures on } \xi\oplus \theta^1 \},  
\end{equation*} 
and these maps respect the maps $u_{\xi}$, i.e. $u_{\xi \oplus \theta^1}= 
v_{\xi}u_{\xi}$. So, we have the map  
\begin{equation}\label{lim} 
\lim_{n\to\infty}\,\{ \text{PL structures on } \xi \oplus \theta^n\} 
\longrightarrow [\Lift_{\ga} f]  
\end{equation} 
where $\lim$ means the direct limit of the sequence of sets. 
 
\begin{prop}\label{p:lim} 
If $X$ is a finite $CW$-space then the map $\eqref{lim}$ is a bijection. 
\end{prop} 
 
\p The surjectivity follows since every compact subset of $BTOP$ is contained 
in some $BTOP_n$. Similarly, every map $X \times I \to BPL$ passes through 
some $BPL_n$, and therefore the injectivity holds.  
\qed  
 
 \m The space $TOP/PL$ is a homotopy associative and homotopy invertible $H$-space, and hence the set $[X, TOP/PL]$ has a natural group structure. Here the neutral element is the homotopy class of inessential map $X \to TOP/PL$. Now, consider a principal $F$-fibration $F \to E \to B$ as in \defref{d:prin} and apply it to the case 
 \[
 \CD
 TOP/PL @>\gb >> BPL @>\ga >> BTOP. 
 \endCD
 \]
 Then for every map $f: X \to BTOP$ we have a right action
 \[
r: [\Lift_{\ga}f]\times [X,TOP/PL]\longrightarrow  [\Lift_{\ga}f]
 \]
 
 \begin{prop}\label{p:top/pl} 
 Suppose that the map $f: X \to BTOP$ lifts to $BPL$. Then the action $r$ is transitive. 
 Furthermore, for every $\ga$-lifting $g$ of $f$ the map
 \[
 [X, TOP/PL]\longrightarrow [\Lift_{\ga}f], \quad \gf \mapsto r(g,\gf)
 \]
 is a bijection.
 \end{prop}
 \p See \theoref{lift}.
 \qed
 
\m Note that, in view of Propositions \ref{lim} and \ref{p:top/pl}, if a topological bundle $\xi$ admits a PL structure then  the bijection \eqref{lim} turns into the bijection 
\begin{equation}\label{eq:tpl-lim} 
\lim_{n\to\infty}\,\{\text{PL structures on } \xi \oplus 
\theta^n\}\longrightarrow [X,TOP/PL]
\end{equation} 
 provided $X$ is a finite CW space.
 
\begin{definition} \rm \label{d:hstruct}
Let $M$ be a $\partial_{PL}$-manifold. A {\it 
homotopy PL structuralization} of $M$ is a homotopy equivalence $h: (V, \pa V)  
\to (M, \pa M)$ such that $V$ is a PL manifold and $h|_{\pa V}: \pa V \to \pa M$ is a PL homeomorphism. Two homotopy PL structuralizations $h_i: V_i  
\to M, i=0,1$ are {\it equivalent} if there exists a PL homeomorphism  
$\gf: V_0\to V_1$ such that $h_1\gf$ is homotopic to $h_0$ rel $\pa V$. In detail, there is a homotopy $H: V_0\times I \to M$ such that 
$H|_{V\times \{0\}}=h_0$ and $H|_{V_0\times \{1\}}=h_1\gf$ and,  
moreover, $H|_{V\times  \{t\}}: 
\partial V_0 \to \partial M$ coincides with $h_0$. Any equivalence 
class of homotopy PL structuralizations is called a {\it homotopy PL 
structure on $X$}. We denote by  $\spl(X)$ the set of all homotopy PL 
structures on $X$. 
 
\m If $M$ itself is a PL manifold, we define the {\it trivial} element of $\spl(M)$ as the equivalence class of $1_M: M \to M$. 
\end{definition} 

\m Pay attention to the map
\begin{equation}\label{forg} 
\CD 
\tpl(M) @>\phi>>  \spl(M)
\endCD 
\end{equation} 
that regards a PL structuralization as the homotopy PL 
structuralization.
 
\begin{definition} \label{d:bulstruct}
\rm 
Given an $(S^n,*)$-fibration $\xi$ over $X$, a {\it homotopy PL 
structuralization of $\xi$} is an $(S^n,*)$-morphism $\gf: \xi \to 
(\gga^n_{PL})\bul$. We say that two PL structuralizations $\gf_0, 
\gf_1: \xi \to (\gga^n_{PL})\bul$ are {\it equivalent} if there exists a 
morphism $\Phi: \xi\times 1_I \to (\gga^n_{PL})\bul$ of 
$(S^n,*)$-fibrations such that $\Phi|_{\xi\times 1_{\{i\}}}=\gf_i, 
i=0,1$. 
Every such an equivalence class is called a {\it homotopy PL structure 
on $\xi$}. 
\end{definition} 
 
\m 
Now, similarly to \eqref{eq:tpl-lim}, for a finite $CW$-space  
$X$ we have a bijection 
\begin{equation}\label{fpl-lim} 
\lim_{n\to\infty}\,\{ \text{homotopy PL structures on } \xi \oplus 
\theta^n\}\longrightarrow [X,F/PL]. \end{equation} 
However, here we can say more. 
 
\begin{prop}\label{fpl-stab}  
The sequence  
$$ 
\left\{\{ \text{\rm homotopy PL structures on } \xi \oplus 
\theta^n\}\right\}_{n=1}^{\infty} $$ 
stabilizes. In particular, the map  
$$ 
\{ \text{\rm homotopy PL structures on } \xi \oplus \theta^n\} \to [F/PL] 
$$ 
is a bijection if $\dim \xi >\!\!> \dim X$ 
\end{prop} 
 
\p This follows from \ref{stability}. 
\qed 
 
\m 
Thus, for every $\RR^N$-bundle $\xi$ that admits a PL structure we have a commutative diagram 
\begin{equation}\label{eq:fpl}
\CD 
\{\text{PL structures 
on $\xi$}\} @>>>  [X,TOP/PL]  \\ 
 @VVV @Va_*VV \\  
\{\text{homotopy PL structures 
on $\xi\bul$}\}@>>> [X,F/PL] 
\endCD 
\end{equation}
 
\m Here the right vertical map $a$ in \eqref{a-b} induces the map $a_*: [X, TOP/PL] \to [X, 
F/PL]$. The left vertical arrow converts a morphism of $\RR^N$ -bundles into
a morphism of $(S^N,*)$-bundles and regards the last one  as a morphism of
$(S^N,*)$-fibrations.   

For a finite $CW$-space $X$, the horizontal arrows turn into bijections if 
we stabilize the picture. i.e. pass to the limit as in \eqref{eq:tpl-lim}. 
Furthermore, the bottom arrow is an isomorphism if $N>\!\!>\dim X$. 
 
\begin{remark}\rm \label{r:stabletop}
Actually, following the proof of the Main Theorem, one can prove that 
$TOP_m/PL_m=K(\ZZ/2,3)$ for $m\geqslant 5$, see \cite[Essay V, \S 5]{KS2}. So, an 
obvious analog of \ref{stability} holds for $TOP$ also, and therefore the top 
map of the above diagram is a bijection for $N$ large enough. But, of course, 
we are not allowed to use these a posteriori arguments here, until we accomplish the proof of the Main Theorem.   
\end{remark}  
 
\begin{remark}\label{smooth} 
\rm
  We can also consider {\it smooth} (= differentiable $C^\infty$) 
structures on topological manifolds. To do this, we must replace the 
words ``PL'' in Definition \ref{d:str} by the word ``smooth''. The 
related set of smooth concordance classes is denoted by $\td(M)$. 
 
The set $\sd(M)$ of homotopy smooth structures is defined in a similar way: replace the 
words ``PL'' in Definition \ref{d:hstruct} by the word ``smooth''.

Moreover, every smooth manifold can be canonically 
converted into a PL manifold (S.~Cairns and J.~Whitehead Theorem~\cite{Cai, W}, 
see e.g. \cite{HM}). So, we can define the set $\pd(M)$ of smooth 
structures on a PL manifold $M$. To do this, we must modify definition 
\ref{str} as follows: $M$ is a PL manifold with a compatible smooth 
boundary, $V_i$ are smooth manifolds, $h_i$ and $H$ are PL 
isomorphisms. 
\end{remark}
 
 \m For convenience of references, we fix here the following theorem of 
Smale \cite{Sma}. Actually, Smale proved it for smooth manifolds, a good 
proof can also be found in Milnor \cite{Mi4}. However, the proof can  
be transmitted  to the PL case, see Stallings \cite[8.3, Theorem 
A]{Stallings}. 

\begin{theorem}\label{Sma}  
Let $M$ be a closed PL manifold that is homotopy equivalent 
to the sphere $S^n, n\geqslant 5$. Then $M$ is PL homeomorphic to $S^n$.  
\qed 
\end{theorem}

\begin{example}\rm
Now we construct an example of two smooth structuralizations 
$h_i: V \to S^n, i=1,2$ that are not concordant. First, note that 
there is a bijective correspondence between $\sd(S^n)$ and the 
Kervaire--Milnor group $\Theta_n$ of homotopy spheres, \cite{KM}. 
Indeed, $\Theta_n$ consists of equivalence classes of oriented 
homotopy spheres: two oriented homotopy spheres are equivalent if they 
are orientably diffeomorphic (= $h$-cobordant). 
 Now, given a homotopy smooth structuralization $h: \Sigma^n 
\to S^n$, we orient $\Sigma^n$ so that $h$ has degree 1. In this way we 
get a well-defined map $u:\sd (S^n) \to \Theta_n$. Conversely, given a 
homotopy sphere $\Sigma^n$, consider a homotopy equivalence $h: \Sigma^n \to 
S^n$ of degree 1. In this way we get a well-defined map $\Theta_n \to 
\sd(S^n)$ which is inverse to $u$. 
 
Note that, because of the Smale Theorem, every smooth homotopy sphere 
$\Sigma^n, n\geqslant 5$, possesses a smooth function with exactly two critical 
points. Thus, $\sd(S^n)=\td(S^n)=\pd(S^n)$ for $n\geqslant 5$. Kervaire and 
Milnor \cite{KM} proved that $\Theta_7=\ZZ/28$, i.e., because of what 
we said above, $\sd(S^7)=\td(S^7)=\pd(S^7)$ consists of 28 elements. 
 
On the other hand, there are only 15 classes of diffeomorphism of smooth manifolds which are 
homeomorphic (and PL homeomorphic, and homotopy equivalent) to $S^7$ but mutually
non-diffeomorphic. Indeed, if an oriented smooth 
7-dimensional manifold $\Sigma$ is homeomorphic to $S^7$ then $\Sigma$ 
bounds a parallelizable manifold $W_{\Sigma}$, \cite{KM}. We equip $W$ an 
orientation which is compatible with $\Sigma$ and set 
$$ 
a(\Sigma)=\frac{\sigma(W_{\Sigma})} 8 \mod 28 
$$ 
where $\sigma(W)$ is the signature of $W$. Kervaire and Milnor 
\cite{KM} proved that the correspondence 
$$ 
\Theta_7 \to \ZZ/28, \quad \Sigma \mapsto a(W_{\Sigma}) 
$$ 
is a well-defined bijection. 
 
However, if $a(\Sigma_1)=-a(\Sigma_2)$ then $\Sigma_1$ and $\Sigma_2$ 
are diffeomorphic: namely, $\Sigma_2$ is merely the $\Sigma_1$ with the 
opposite orientation. So, there are only 15 smooth manifolds which are 
homeomorphic (and homotopy equivalent, and PL homeomorphic) to $S^7$ but mutually 
non-diffeomorphic. 
 
In terms of structures, it can be expressed as follows. Let $\rho: S^n 
\to S^n$ be a diffeomorphism of degree -1. Then the smooth 
structuralizations $h: 
\Sigma^7 \to S^7$ and $\rho h: \Sigma^7 \to S^7$ are not concordant, 
if $a(\Sigma^7) \neq 0,14$. 
\end{example}

\section{From Manifolds to Bundles}\label{from} 
 
Recall that, for every topological manifold $M^n$, its tangent bundle 
$\tau_M$ and normal (with respect to an embedding $M \in \RR^{N+n}, N\gg n$) PL $\RR^N$-bundle $\nu_M$ are defined. Here $\tau_M$ 
is a topological $\RR^n$-bundle, and we can regard $\nu_N$ as a 
topological $\RR^N$-bundle. Furthermore, if $M$ is a PL 
manifold then $\tau_M$ and $\nu_M$ turns into PL bundles in a 
canonical way, see \cite{KS2, Rud}. 
 
 \m {\footnotesize Concerning tangent and normal (micro)bundles and their properties, see~Milnor~\cite{Mi3} for the topological category and Haefliger--Wall~\cite{HW} for the PL category.}
 
\begin{construction}\label{jtop} 
\rm 
 Consider a manifold $M$ (possibly with boundary) and a PL structuralization $h: V \to M$. 
 Let $g=h^{-1}: M \to V$. Since $g$ is a homeomorphism, it yields a 
topological  morphism $\gl_g: \tau_M \to \tau_V$ where $\tau_V, \tau_M$ denote the tangent bundles to $V, M$ respectively,  and so we have the correcting topological morphism  $c(\gl_g): \tau_M \to \gl^*\tau_V$. Now, let $\nu=\nu_M^N$ be a normal bundle of  $M$ in  $\RR^{N+n}$ with $N$ large enough. Consider the topological morphism  
$$  
\CD 
\theta_M^{N+n}=\tau_M\oplus \nu^{N}_M @>>> g^*\tau_V \oplus \nu^N_M 
@>\text{classif} >> \gga^{N+n}_{PL} \endCD 
$$ 
and regard it as a PL structuralization of $\theta^{N+n}$.
 It is easy to see that 
in this way we have the correspondence 
\begin{equation}\label{jtop1}
\jtop: \tpl(M) \longrightarrow \lim_{n\to\infty}\,\{\text{PL structures on } 
\theta_M^N\} \longrightarrow [M, TOP/PL] 
\end{equation}
where the last map comes from \eqref{eq:tpl-lim}.  

\m
Moreover, since $g: \pa M \to \pa V$ is a PL homeomorphism, we have a commutative diagram
\[
\CD
\tpl(M) @>\jtop >> [M,TOP/PL]@=[M, TOP/PL]\\
@VVV @VVV @VVV\\
\tpl(\pa M) @>>> \pt @>>> [\pa M, TOP/PL]
\endCD
\]
 cf. \remref{conc}(2).
So, we can (and sometimes shall) regard the map $\jtop$ from \eqref{jtop1} as the map 
\begin{equation}\label{jtop2}
\jtop: \tpl(M) \longrightarrow [(M, \partial M), (TOP/PL, *)]. 
\end{equation}
\end{construction} 
 
\m Now we construct a map $j_F: \spl (M) \to [M,  F/PL]$, 
a ``homotopy analogue'' of $\jtop$. This construction is more 
delicate, and we treat here the case of closed manifolds only.   
So, let $M$ be a connected closed PL manifold.

\begin{definition}\rm 
Given an $(S^n,*)$-fibration $\xi=\{E \to B\}$ with a section $s: B 
\to E$, we define its {\it Thom space} $T\xi$ as the quotient space 
$E/s(B)$, Given a topological $\RR^N$-bundle $\eta$, we define the Thom space 
$T\eta$ as $T\eta:=T(\eta\bul)$. 
\end{definition} 
 
Given a morphism $\gf: \xi \to \eta$ of $(S^n,*)$-fibrations, we define $T\gf: 
T\xi \to T\eta$ to be the unique map such that the diagram 
$$ 
\CD 
E @>>> E'\\ 
@VVV @VVV \\ 
T\xi @> T\gf >> T\eta 
\endCD 
$$  
commutes. Here $E'$ is the total space of $\eta$. 
 
\begin{df} 
\label{df-red} 
\rm A pointed space $X$ is called {\it reducible} if there is a pointed map $f: 
S^m \to X$ such that $f_*: \wt H_i(S^m) \to \wt H_i(X)$ is an isomorphism for 
$i\geqslant m$. Every such map $f$ (as well as its homotopy class or its stable 
homotopy class) is called a {\it reducibility} for $X$.  
\end{df} 
 
\m We embed $M^n$ in $\mathbb R^{N+n}, N\gg n$ and let $\nu=\nu_M, 
\dim \nu=N$ be a normal bundle of this embedding. Recall that $\nu$ is a 
PL bundle $E \to M$ whose total space $E$ is PL homeomorphic to a (tubular) 
neighborhood $U$ of $M$ in $\RR^{N+n}$. We choose such isomorphism and denote 
it by $\gf: U \to E$.

\begin{construdef}\label{collapse} 
\rm 
Let $T\nu$ be the Thom space of $\nu$. Then there is a unique map  
$$ 
\psi: \RR^{N+n}/(\RR^{N+n}\setminus U) \to T\nu 
$$ 
such that $\psi|_U=\gf$. We define the {\it collapse map}  $\iota: S^{N+n} \to 
T\nu_M$ (the {\it Browder--Novikov map}, cf.~\cite{Br, N1}) to be the composition  
$$ 
\CD 
\iota: S^{N+n} @> \text{quotient}>> S^{N+n}/(S^{N+n}\setminus U) = 
\RR^{N+n}/(\RR^{N+n}\setminus U) @>\psi >>  T\nu.  
\endCD 
$$ 
\end{construdef} 
 
\m It is well known and easy to see that $\iota$ is a reducibility for $T\nu$, see 
Corollary \ref{reducibility} below.

\m It turns out that, for $N$ large enough, the normal bundle of a given 
embedding $M \to \RR^{N+n}$ exists and is unique up to isomorphism. For detailed definitions and proofs, 
see \cite{HW, KL, LR}. The uniqueness gives us the following important fact. 
Let $\nu'=\{ E' \to M\}$ be another normal bundle and $\gf': U' \to E'$ be 
another PL homeomorphism. Let $\iota: S^{N+n} \to T\nu$ and $\iota': S^{N+n} \to 
T\nu'$ be the corresponding Browder--Novikov maps. Then there is a morphism 
$\nu \to \nu'$ of PL bundles which carries $\iota$ to a map homotopic to 
$\iota'$.

\begin{theorem}\label{t:bn} 
Consider a PL $\mathbb R^N$-bundle $\eta$ over $M$ such that $T\eta$ is 
reducible. Let $\ga\in \pi_{N+n}(T \eta)$ be an arbitrary reducibility for 
$T\eta$. Then there exist an $(S^N,*)$-equivalence $\mu: \nu_M\bul \to \eta\bul $ 
such that $(T\mu)_*(\iota)=\ga$, and such a $\mu$ is unique up to fiberwise 
homotopy over $M$. 
\end{theorem} 
 
The Theorem is a version of the Spivak Theorem~\cite[Theorem A]{Spi}, cf. also~\cite[I.4.19]{Br}. Note that our version does not require the simply-connectedness of $M$. We postpone the proof to the next Chapter, see \secref{BN}.

\m Given a homotopy equivalence $h:V \to M$ of closed connected PL manifolds, 
let $\nu_V$ be a normal bundle  of a certain embedding $V\subset \RR^{N+n}$, 
and let  $u\in\pi_{N+n}(T\nu_V)$ be the homotopy class of the collapsing map 
$S^{N+n} \to T\nu_V$. Let $g: M \to V$ be homotopy inverse to $h$ and set 
$\eta=g^*\nu_V$. The $g$-adjoint morphism  
$$ 
\mathfrak I=\mathfrak I_{g} : \eta =g^*\nu V  \to \nu_V 
$$ 
yields the map $ T\mathfrak I: T \eta \to T\nu_V$. It is easy to see that $T\\mathfrak I$ is a 
homotopy equivalence, and so there exists a unique $\ga \in 
\pi_{N+n}(T\eta)$ with $(T\mathfrak I)_*(\ga)=u$. Since $u$ is a reducibility for 
$T\nu_V$, we conclude that $\ga$ is a reducibility for $T\eta$. By 
Theorem \ref{t:bn}, we get an $(S^N,*)$-equivalence $\mu: 
\nu_M\bul \to \eta\bul$ with $(T\mu)_*(\iota)=\ga$. Now, the 
morphism    
\begin{equation} \label{eq:nu_to_gamma}
\CD
(\nu_M)\bul  @>\mu>> 
\eta\bul @>\text{ classif }>> \gga^{N}_F 
\endCD 
\end{equation} 
is a homotopy PL structuralization of $\nu_M$. Because of the uniqueness of 
the normal bundle, the homotopy class of this structuralization is well 
defined. So, in this way we have the function 
\begin{equation*} 
\jf: \spl(M) \longrightarrow  \{\text{homotopy PL structures 
on $\nu_M$}\} 
\cong   [M,F/PL] 
\end{equation*} 
where the last bijection comes from \ref{fpl-stab}. 
 
\begin{df}\rm\label{d:norminv}
The function $\jf$ is called the {\it normal invariant}, and its 
value on a homotopy PL structure (as well as on its PL structuralization) is called the normal invariant of 
this structure (structuralization). 
\end{df}

\section{Homotopy PL Structures on $T^k \times D^n$} 
 
Below $T^k$ denotes the $k$-dimensional torus. 
 
\begin{theorem}\label{t:norminv} 
Assume that $k+n\geqslant 5$. If $x\in \spl(T^k\times S^n)$ can be  
represented by a homeomorphism $M \to T^k \times S^n$ then $j_F(x)=0$.  
\qed
\end{theorem} 
 
This is a special case of the Sullivan Normal Invariant Homeomorphism 
Theorem. We prove \ref{t:norminv} (in fact, a little bit 
general result) in the next chapter. 
 
 \m
{\footnotesize We also prove the Sullivan Theorem in full generality in 
Chapter~\ref{c:norminv}, Section~\ref{secnorm}.  We must do this loop (repetition) since the proof 
in Chapter~\ref{c:norminv} uses Main Theorem and hence \theoref{t:norminv}.}
 
\begin{construdef} 
\rm 
\m Let $x\in \spl(M)$ be represented by a map $h: V \to M$, and 
let $p: \wt M \to M$ be a covering. Then we have a commutative diagram 
$$ 
\CD 
\wt V @>\wt h>> \wt M\\ 
@VqVV @VVpV\\ 
V @ >h>> M 
\endCD 
$$ 
where $q$ is the induced covering. Since $\wt h$ is defined uniquely 
up to deck transformations, the concordance class of $\wt h$ is well 
defined. So, we have a well-defined map 
$$ 
p^*: \spl(M) \to \spl (\wt M) 
$$ 
where $p^*(x)$ is the concordance class of $\wt h$. 

If $p$ is a finite covering, we say that a class 
$p^*(x)\in \spl(\wt M)$ {\it finitely covers} the class $x$. 
\end{construdef} 

\begin{theorem}\label{HS,W3} 
Let $k+n\geqslant 5$ Then the following holds: 
\par{\rm (i)} if $n>3$ then the set $\spl(T^k \times D^n)$ consists of 
precisely one $($trivial$)$ element; 
\par{\rm (ii)} if $n<3$ then every element of $\spl(T^k \times D^n)$ 
can be finitely covered by the trivial element; 
\par{\rm (iii)} the set  $\spl(T^k \times D^3)$ contains at most one 
element which cannot be finitely covered by the trivial element. 
\end{theorem} 
 
{\it Some words about the proof.} \rm First, we mention the proof given by 
Wall, \cite{Wall2} and \cite[Section 15 A]{Wall3}. Wall proved the bijection $w:
\spl(T^k \times D^n)  \to H^{3-n}(T^k)$. Moreover, he also proved that finite
coverings  respect this bijection, i.e. if $p: T^k \times D^n \to T^k \times
D^n$  is a finite covering then there is the commutative diagram 
$$ 
\CD 
\spl(T^k \times D^n) @>w>> H^{3-n}(T^k;\ZZ/2)\\ 
@Ap^*AA @AAp^*A\\ 
\spl(T^k \times D^n)  @>w>> H^{3-n}(T^k;\ZZ/2)\,. 
\endCD 
$$ 
Certainly, this result implies all the claims(i)--(iii). Wall's proof uses
difficult algebraic calculations.

\m   
Another proof of the theorem can be found in~\cite[Theorem C]{HS}. Minding the complaint of Novikov concerning 
Sullivan's results (see Preface), we mention that the nice paper of Hsiang and Shaneson~\cite{HS} use a Sullivan's result. Namely,  Hsiang and Shaneson consider the so-called surgery exact sequence 
$$ 
\CD 
@>\partial >> \spl(S^k \times T^n) @>\jf>>[S^k \times T^n, F/PL] @>>> 
\cdots 
\endCD 
$$ 
and write (page 42, Section 10): 
\begin{quote} 
By [44], every homomorphism $h: M \to S^k \times T^n, k=n \geq 5$, 
represents an element in the image of $\partial$. 
\end{quote} 
 
Here the item [44] of the citation is our bibliographical item ~\cite{Sul1}.  So, in fact, Hsiang and Shaneson use Theorem~\ref{t:norminv}. 
As I already said, we prove~\ref{t:norminv} in next Chapter.
\qed

\section[The Product Structure Theorem]{The Product Structure Theorem, or from Bundles to Manifolds} 
 
Let $M$ be an $n$-dimensional $\partial_{PL}$-manifold. Then every PL 
structuralization  $h: V \to M$ yields a PL structuralization   
$$ 
h\times 1: V \times \RR^k \to M \times \RR^k. 
$$  
Thus, we have a well-defined map  
$$ 
e: \tpl(M) \to \tpl(M\times \RR^k). 
$$ 
 
\begin{theorem}[\rm The Product Structure Theorem] 
\label{product} 
For every $n\geqslant 5$ and every $k\geqslant 0$, the map $e: \tpl(M) \to \tpl(M\times 
\RR^k)$ is a bijection.  
\end{theorem} 
 
\sm In particular, if $\tpl(M\times \RR^k) \ne \emptyset$ then $\tpl(M) \ne \emptyset$. 
 
\m Kirby and Siebenmann made the breakthrough for 6.1~\cite{K1, KS1,KS2}. Quinn~\cite{Q2} gave a nice short proof of 6.1 by developing his theory of ends of maps,~\cite{Q1}.
 
 \begin{corollary}[The Classification Theorem] 
\label{classif} 
If $\dim M \geqslant 5$ and $M$ admits a PL structure, then the map 
$$ 
\jtop: 
\tpl(M) \to [(M, \partial M), (TOP/PL, *)] 
$$ 
is a bijection. 
\end{corollary} 

\p 
We construct a map 
\begin{equation}\label{sigma} 
\sigma: [(M, \partial M), (TOP/PL, *)] \to \tpl(M) 
\end{equation} 
which is inverse to $\jtop$. Take an element 
\[
a\in [(M, \partial M), (TOP/PL, *)]
\]
 and, using 
\eqref{eq:tpl-lim}, interpret it as a homotopy class of a topological $\RR^N$-morphism $\gf: 
\theta^N_M \to \gga^N_{PL}$ such that $\gf|_{\pa M}$ is a PL $\RR^N$-morphism. The morphism $\gf$ yields a correcting isomorphism $\theta^N_M 
\to b^*\gga^N_{PL}$ of topological $\RR^N$-bundles over $M$, where $b: M \to 
BPL$ is the base of the morphism $\gf$. So, we have the commutative diagram 
$$ 
\CD 
M \times \RR^N @>h>> W\\ 
@VVV @VVV\\ 
M @=M 
\endCD 
$$ 
where $h$ is a fiberwise homeomorphism and $W \to M$ is a PL 
$\RR^N$-bundle $b^*\gamma^N_{PL}$. In particular, $W$ is a PL manifold. 
Regarding $h^{-1}: W \to M \times \RR^N$ as a PL structuralization of 
$M\times \RR^N$, we conclude that, by the Product Structure Theorem 
\ref{product}, $h^{-1}$ is concordant to a map $g\times 1$ for some PL 
structuralization $g: V \to M$. We define $\sigma(a)\in \tpl(M)$ to be 
the concordance class of $g$. One can check that $\sigma$ is a 
well-defined map which is inverse to $\jtop$. Cf. \cite[Essay IV Theorem 4.1]{KS2}. 
\qed 
 
\begin{corollary}[The Existence Theorem] 
\label{exist} 
A topological manifold $M$ with $\dim M \geqslant 5$ admits a PL  
structure if and only if the tangent bundle of $M$ admits a PL  
structure.  
\end{corollary} 
 
\p Only claim ``if'' needs a proof. Let $\tau=\{\pi: D \to M\}$ be 
the tangent bundle of $M$, and let $\nu=\{r: E \to M\}$ be a stable normal 
bundle of $M$, $\dim \nu = N$. Then $E$ is homeomorphic to an open subset of 
$\RR^{N+n}$, and therefore we can (and shall) regard $E$ as a PL manifold. 
Since $\tau$ is a PL bundle, we conclude that $r^*\tau$ is a PL bundle over 
$E$. In particular, the total space $M\times \RR^{N+n}$ of $r^*\tau$ turns 
out to be a PL manifold, cf.~\cite{Mi3}. Now, because of the Product Structure Theorem 
\ref{product}, $M$ admits a PL structure.  Cf. \cite[Essay IV Theorem 3.1]{KS2}
\qed 
 
\m Let $f:M \to BTOP$ classify the stable tangent bundle of a 
closed topological manifold $M$, $\dim M \geqslant 5$. 
 
\begin{cory}\label{c:exist} 
The following conditions are  equivalent: 
\par{\rm (i)} $M$ admits a PL structure; 
\par{\rm (ii)} $\tau$ admits a PL structure; 
\par{\rm (iii)} there exists $k$ such that $\tau\oplus \theta^k$ admits a PL 
structure;  
\par{\rm (iv)} the map $f$ admits an $\ga^{PL}_{TOP}$-lifting to $BPL$. 
\end{cory} 
 
\p It suffices to prove that (iv) $\Longrightarrow $ (iii) $\Longrightarrow $ 
(i). The implication (iii) $\Longrightarrow $ (i) can be proved similarly to 
\ref{exist}. Furthermore, since $M$ is compact, we conclude that $f(M) \suset 
BTOP_m$ for some $m$. So, if (iv) holds then $f$ lifts to $BPL_m$, i.e. 
$\tau \oplus \theta^{m-k}$ admits a PL structure. 
\qed

\begin{remark} 
\rm  It follows from \ref{lift}, \ref{exist} and 
\ref{classif} that the set $\tpl(M)$ of PL structures on $M$ is in a bijective 
correspondence with the set of fiber homotopy classes of 
$\ga^{PL}_{TOP}$-liftings of $f$.  
\end{remark} 
 
\begin{remark}\label{ker-coker} 
\rm 
It is well known that $\jf$ is not a bijection in general. The ``kernel''' 
and ``cokernel'' of $\jf$ can be described in terms of so-called Wall groups, 
\cite{Wall3}. (For $M$ simply-connected, see also \theoref{bns}.) On the 
other hand, the bijectivity of $\jtop$ (the Classification Theorem) follows 
from the Product Structure Theorem. So, informally speaking, kernel and 
cokernel of $\jf$ play the role of obstructions to splitting of structures. 
It seems interesting to develop and clarify these naive arguments.  
\end{remark}  
 
\begin{remark} \label{pl/o} 
\rm Since tangent and normal bundles of smooth manifolds turn out to be 
vector bundles, one can construct a map   
$$ 
k: \pd(M)  \to [M, PL/O] 
$$ 
which is an obvious analog of $\jtop$. 
Moreover, the obvious analog of the Product Structure Theorem (as well as of 
the Classification and Existence Theorems) holds without any dimensional 
restriction. In particular, $k$ is a bijection for every smooth manifold, 
\cite{HM}.    
 
It is well known (although difficult to prove) that $\pi_i(PL/O)=0$ for $i 
\leqslant 6$. (See \cite[IV.4.27(iv)]{Rud} for the references.) Thus, every PL 
manifold $M$ of dimension $\leqslant 7$ admits a smooth structure, and this 
structure is unique if $\dim  M \leqslant 6$.   
\end{remark}

\section{Non-contractibility of $TOP/PL$} 

\begin{theorem}[{\rm Rokhlin Signature Theorem}]\label{rohlin} 
Let $M$ be a closed $4$-dimensional PL manifold with 
$w_1(M)=0=w_2(M)$. Then the signature of $M$ is divisible by $16$. 
\end{theorem}

\p 
See~\cite{MK}, \cite[XI]{K2}, or the original work \cite{Ro}. In fact, Rokhlin  
proved the result for smooth manifolds, but the proof works for PL manifolds 
as well. On the other hand, in view of \ref{pl/o}, there is no 
difference between smooth and PL manifolds in dimension 4. 
\qed 
 
 \begin{theorem}[{\rm Freedman's Example}]\label{t:F} 
There exists a closed simply-connected topological $4$-dimensional  
manifold $V$ with $w_2(V)=0$ and such that $E_8$ is the matrix of the intersection 
form $H^2(V)\otimes H^2(V)\to \ZZ$. In particular, the signature of $V$ 
is equal to $8$. Furthermore, such a manifold $V$ is unique up to homeomorphism.
\end{theorem} 
 
\p 
See~\cite{FU},~\cite{FQ}, or the original work~\cite{F}. 
\qed 

\begin{comment}\rm \label{plumb}
Some words on constructing of $V$. Take the manifold $W$ (plumbing) from Browder~\cite[Complement V.2.6]{Br}. This is a smooth 4-dimensional simply-connected parallelizable manifold whose boundary $\partial W$ is a homology sphere. Furthermore, $E_8$ is the matrix of the intersection form $H_2(W)\otimes H_2(W)\to \ZZ$. A key (and very difficult) result of Freedman~\cite{F} claims that $\pa W$ bounds a compact contractible topological manifold $P$. Now, put $V=W\cup_{\pa W}P$. 
\end{comment}

\begin{corollary}\label{c:F} 
The topological manifolds $V$ and $V\times T^k, k\geqslant 1$ do not admit any PL 
structure. 
\end{corollary} 
 
\p 
The claim about $V$ follows from \ref{rohlin}. (Note that $w_1(V)=0$ because $V$ is simply-connected.)
Suppose that $V\times  T^k$ has a PL structure. Then $V\times \RR^k$ has a PL structure. So, because 
of the Product Structure Theorem \ref{product}, $V\times \RR$ has a PL 
structure. Hence, by \ref{pl/o}, it possesses a smooth structure.  Then the 
projection $p_2: V \times \RR \to \RR$ can be $C^0$-approximated by a map $f: 
V \times \RR \to \RR$ which coincides with $p_2$ on $V\times (-\infty,0)$ and 
is smooth on $V\times (1,\infty)$. Take a regular value $a\in (0,\infty)$ of 
$f$ (which exists because of the Sard Theorem) and set $U=f^{-1}(a)$. Then $U$ 
is a smooth manifold (by the Implicit Function Theorem), and it is easy to see 
that $w_1(U)=0=w_2(U)$ (because it holds for both manifolds $\RR$ and $V\times 
\RR$). On the other hand,  both manifolds $V$ and $U$ cut the ``tube''  $V
\times \RR$. So, they are (topologically) bordant, and therefore $U$ has 
signature $8$. But this contradicts the Rokhlin Theorem \ref{rohlin}. 
 
\qed 
 
\begin{corollary}\label{noncontr} 
The space $TOP/PL$ is not contractible. 
\end{corollary} 
 
\p Indeed, suppose that  $TOP/PL$ is contractible. Then every map $X \to BTOP$ 
lifts to $BPL$, and so, by \ref{exist}, every closed topological manifold of 
dimension greater than 4 admits a PL structure. But this contradicts~\ref{c:F}. 
\qed 
 
\begin{remark}\rm 
Kirby and Siebenmann \cite[Annex C]{KS2} constructed the original  example of a 
topological manifold which does not admit a PL structure. Namely, they considered the space
\[
X^4=T^4 \# \text{ cone of }W
\]
with $W$ as in \ref{plumb} and proved that $X\ts S^1$ is a topological manifold. If we assume that $X\times S^1$ and argue as
in the end of the proof of \ref{c:F}, we construct a manifold $Y$ (an analog of $U$ in \ref{c:F}) with $w_1(Y)=0=w_2(Y)$ and $\gs(Y)=8$. Thus, the 5-manifold $X\ts S^1$ does not admit PL structure.
\end{remark}

\section{Homotopy groups of $TOP/PL$} 
 
Let $M$ be a compact topological manifold equipped with a metric 
$\rho$. Then the space $\HH$ of self-homeomorphisms $M\to M$  gets a metric $d$ with 
$d(f,g)=\sup\{x\in M\bigm|\rho(f(x), g(x))\}$. 
 
\begin{theorem}\label{chern} 
The space $\HH$ is locally contractible. 
\end{theorem} 
 
\p See \cite{Ch,EK}. 
\qed 
 
\begin{corollary}\label{eps} There exists $\eps>0$ 
such that every homeomorphism $h\in \HH$ with $d(H, 1_M)<\eps$ is  
isotopic to $1_M$. 
\qed 
\end{corollary} 
 
\begin{construction} \label{useful} 
\rm 
We regard the torus $T^k$ as a commutative Lie group (multiplicative) 
and equip it with the invariant metric $\rho$. Consider the map 
$p_{\lambda}: T^k \to T^k, p_{\lambda}(a)=a^{\lambda}, \lambda \in 
\NN$. Then $p_{\lambda}$ is a $\lambda^k$-sheeted covering. It is 
also clear that all the deck transformations of the covering torus are 
isometries. Hence the diameter of each of (isometric) 
fundamental domain for $p_{\gl}$ tends to zero as $\lambda \to 
\infty$. 
  
\end{construction} 
 
\begin{lemma} 
Let $h: T^k \times D^n \to T^k \times D^n$ is a self-homeomorphism 
which is homotopic $\rel \partial(T^k \times D^n)$ to the identity. 
Then there exist $\lambda \in 
\NN$ and a commutative diagram 
$$ 
\CD 
T^k \times D^n @>\wt h>> T^k \times D^n \\ @Vp_{\lambda}VV 
@VVp_{\lambda}V\\ T^k 
\times D^n @>h>>   T^k \times D^n 
\endCD 
$$ 
where the lifting $\wt h$ of $h$ is  isotopic $\rel \partial(T^k 
\times D^n)$ to the identity. 
\end{lemma} 
 
\p (Cf. \cite[Essay V]{KS2}.) First, consider the case $n=0$. Without 
loss of generality we can assume that $h(e)=e$ where $e$ is the neutral 
element of $T^k$. Consider a covering $p_{\lambda}: T^k \to T^k$ 
as in \ref{useful} and take a covering $\wt h: T^k \to T^k, 
p_{\lambda}\wt h =\wt h p_{\lambda}$ of $h$ such that $\wt h(e)=e$. 
In 
order to distinguish the domain and the range of $p_{\lambda}$, we 
denote the domain of $p_{\lambda}$ by $\wt T$ and the range of 
$p_{\lambda}$ by $T$. Since $h$ is homotopic to $1_T$, we conclude 
that every point of the lattice $L:=p_{\lambda}^{-1}(e)$ is fixed 
under $\wt h$. 
 
Given $\eps>0$,  choose $\delta$ such that $\rho(\wt h(x), \wt 
h(y))<\eps/2$ whenever $\rho(x,y)<\gd$. Furthermore, choose $\lambda$ 
so large that the diameter of any closed fundamental domain is less 
then $\min\{\eps/2,\gd\}$. Now, given $x\in \wt T$, choose $a\in L$ 
such that $a$ and $x$ belong to the same closed fundamental domain. 
Now, 
$$ 
\rho(x,\wt h(x)\leqslant \rho(x,a)+\rho(a, \wt h(x))=\rho(x,a)+ 
\rho(\wt h(a), \wt h(x))< \frac \eps 2  + \frac \eps 2 =\eps. 
$$ 
So, for every $\eps>0$ there exists $\lambda$ such that $d(\wt h, 
1_{\wt T})<\eps$. Thus, by \ref{eps}, $\wt h$ is isotopic to $1_{\wt 
T}$ for $\gl$ large enough. 
 
\m 
The proof for $n>0$ is similar but a bit more technical. Let 
$D_{\eta}\subset D^n$ be the disk centered at $0$ and having the radius 
$\eta$. We can always assume that $h$ coincides with identity outside 
of $T^k\times D_{\eta}$. Now, asserting as for $n=0$, take a covering 
$p_{\lambda}$ as above and choose $\gl$ and $\eta$ so small that the 
diameter of every fundamental domain in $\wt T \times D_{\eta}$ is 
small enough. Then 
$$ 
\wt h: \wt T \times D_{\eta} \to  \wt T \times D_{\eta} 
$$ 
is isotopic to the identity and coincides with identity 
outside $\wt T \times D_{\eta}$. This isotopy is not an isotopy rel 
$\wt T \times \partial D_{\eta}$. Nevertheless, we can easily extend 
it to the whole $\wt T \times D^n$ so that this extended isotopy is an 
isotopy $\rel \partial (\wt T \times D^n)$. 
 
If you want an explicit formula, do the following. Given $a=(b,c)\in \wt T \times 
D_{\eta}$, set $|a|=|c|$. Consider an isotopy 
$$ 
\gf: \wt T \times D_{\eta} \times I \to \wt T \times D_{\eta} 
\times I, \quad \gf(a,0)=a, \ \gf(a,1)=\wt h(a), \quad a\in \wt T 
\times D_{\eta}. 
$$ 
Define $\ov \gf : \wt T \times D_{\eta}\times I \to \wt T \times D_{\eta} = 
\times I$ by setting 
\begin{equation*} 
\ov \gf(a,t)= 
\left\{ 
\begin{array}{lcl} 
\gf(a,t) &\text{if } |a|\leqslant \eta, \\ \, 
\gf(a,  \frac{|a|-1} {\eta-1}t) & \text{if } |a|\geqslant \eta. 
\end{array} 
\right. 
\end{equation*} 
Then $\ov\gf$ is the desired isotopy $\rel \partial (\wt T \times D^n)$. 
\qed 
 
\begin{corollary}\label{cover} 
Let $\phi: \tpl( T^k \times D^n) \to \spl( T^k \times D^n)$ be the 
forgetful map as in $\eqref{forg}$. If $\phi(x)=\phi(y)$ then there exists a 
finite covering $p: T^k \times D^n \to T^k \times D^n$ such that 
$p^*(x)=p^*(y)$. 
\qed 
\end{corollary} 
 
\m Consider the map 
\begin{equation*} 
\begin{array}{l l l} 
\psi : &\pi_n(TOP/PL) = [(D^n, \partial D^n), (TOP/PL, *)]  
\xrightarrow{\ p_2^*\ } \\ 
 & [(T^k\times D^n, \partial ( T^k \times D^n)),  
(TOP/PL,*)] 
 \xrightarrow{\sigma}\tpl( T^k \times D^n) 
\end{array} 
\end{equation*} 
where $\sigma$ is the map from \eqref{sigma} (the inverse map to $\jtop$). 
 
\begin{lemma}\label{l:red} 
The map $\psi$ is injective. Moreover, if $p^*\psi(x)=p^*\psi(y)$ for 
some finite covering $p:T^k \times D^n \to T^k \times D^n$ then $x=y$.  
 
In particular, if $p^*\psi(x)$ is the trivial element of $\tpl(T^k 
\times D^n)$ then $x=0$. 
\end{lemma} 
 
\p The injectivity of $\psi$ follows from the injectivity of $p_2^*$ 
and $\sigma$. Furthermore, for every finite covering $p:T^k \times D^n 
\to T^k \times D^n$ we have the commutative diagram 
$$ 
\CD 
\pi_n(TOP/PL) @>\psi>> \tpl (T^k \times D^n)\\ 
@| @AAp^*A\\ 
\pi_n(TOP/PL) @>\psi>> \tpl (T^k \times D^n) 
\endCD 
$$ 
Therefore $x=y$ whenever $p^*\psi(x)=p^*\psi(y)$. Finally, if 
$p^*\psi(x)$ is trivial element then $p^*\psi(x)=p^*\psi(0)$, and thus $x=0$. 
\qed 
 
\m Consider the map 
$$ 
\CD 
\gf: \pi_n(TOP/PL) @>\psi>>  \tpl((T^k \times D^n) @>\phi >> \spl((T^k 
\times D^n)
\endCD 
$$ 
where $\phi$ is the forgetful map described in \eqref{forg}.

\begin{theorem}[The Reduction Theorem, cf.~\cite{K1}] 
\label{t:reduction} 
The map $\gf$ is injective. \\ 
Moreover, if $p^*\gf(x)=p^*\gf(y)$ for 
some finite covering  
$$ 
p:T^k \times D^n \to T^k \times D^n 
$$  
then $x=y$.  
 
In particular, if $p^*\gf(x)$ is the trivial element of $\tpl(T^k 
\times D^n)$ then $x=0$. 
\end{theorem} 
 
\m We call it the Reduction Theorem because it {\it reduces} the calculation 
of the group $\pi_i(TOP/PL)$ to the calculation of the sets $\spl(T^k\times 
D^n)$.

\p If $\gf(x)=\gf(y)$ then $\phi(\psi(x))=\phi(\psi(y))$. Hence, 
by Corollary \ref{cover}, there exists a finite covering $\pi:T^k 
\times D^n \to T^k \times D^n$ such that $\pi^*\psi(x)=\pi^*\psi(y)$.  So,
by  Lemma \ref{l:red}, $x=y$, i.e. $\gf$ is injective. 
 
Now, suppose that $p^*\gf(x)=p^*\gf(y)$ for some finite 
covering $p:T^k \times D^n \to T^k \times D^n$. Then 
$\phi (p^*\psi(x))=\phi (p^*\psi(y))$.  
Now, by Corollary \ref{cover}, there exists a finite covering  
$$ 
q:T^k \times D^n 
\to T^k \times D^n 
$$  
such that $q^*p^*\psi(x)=q^*p^*\psi(y)$, i.e. 
$(pq)^*\psi(x)=(pq)^*\psi(y)$. Thus, by Lemma \ref{l:red}, $x=y$. 
\qed 
 
\begin{corollary}[The Main Theorem] 
\label{main} 
$\pi_i(TOP/PL)=0$ for $i\ne 3$. Furthermore, $\pi_3(TOP/PL)=\ZZ/2$. 
Thus, $TOP/PL=K(\ZZ/2.3)$. 
\end{corollary} 
 
\p The equality $\pi_i(TOP/PL)=0$ for $i\ne 3$ follows from 
Theorem \ref{HS}(i,ii) and Theorem \ref{t:reduction}. Furthermore, again 
because of \ref{HS} and \ref{t:reduction}, we conclude that 
$\pi_3(TOP/PL)$ has at most two elements. In other words, 
$TOP/PL=K(\pi,3)$ where $\pi=\ZZ/2$ or $\pi=0$. Finally, by  
Corollary 
\ref{noncontr}, the space $TOP/PL$ is not contractible. Thus, 
$TOP/PL=K(\ZZ/2,3)$. 
\qed

{\chapter{Normal Invariant}} \label{c:accom}

\m 
The goal of this chapter is to prove Theorem \ref{t:norminv}. The proof uses  the 
Sullivan's result on the homotopy type 
of $F/PL$,~\cite{Sul1, Sul2}. Note that Madsen and Milgram~\cite{MM} gave a detailed proof of 
those Sullivan result. 

We also need \theoref{t:bn}. We prove it in \secref{BN} in the form that is suitable for our aim.

\section{Stable equivalences of spherical bundles} 
 
\m Given a sectioned spherical bundle $\xi$ over a finite $CW$-space $X$, let 
$\aut \xi$ denote the group of fiberwise homotopy classes of 
self-equivalences $\xi \to \xi$ over $X$, where we assume 
that all homotopies preserve the section.  

\m We denote by $\gs^k=\gs^k_X$ the trivial $S^k$-bundle over $X$ with a 
fixed section. In another words, $\gs^k=(\theta^k)\bul$. 
 
\begin{prop}\label{aut} 
There is a natural bijection   
$$ 
\aut \gs^k = [X, F_k]. 
$$ 
\end{prop} 
 
\p Because of the exponential law, every map $X \to F_k$ yields 
a section-preserving map $X \times S^k \to X \times S^k$ over $X$, and vice versa.  
Cf. \cite[Prop. I.4.7]{Br}. \qed 
 
\m Consider the map 
$$ 
\mu: F_k \times F_k \to F_{2k}, \quad \mu(a,b)=a \wedge b: S^{2k}= S^k 
\wedge S^k \to S^k \wedge S^k =S^{2k} 
$$ 
where we regard $a,b \in F_k$ as pointed maps $S^k \to S^k$. Let $T: F_k  
\times F_k \to F_k \times F_k$ be the transpose map, $T(a,b)=(b,a)$. 
 
\begin{lemma}\label{transp} 
The maps $\mu: F_k \times F_k \to F_{2k}$ and $\mu T : F_k \times F_k \to  
F_{2k}, k >0$ are homotopic. 
\end{lemma} 
\p  Consider the map  
$$ 
\tau: S^{2k} =S^k \wedge S^k \to S^k \wedge S^k = S^{2k}, \quad 
\tau(u,v)=(v,u) 
$$  
and note that, for every $a,b \in F_k$, we have 
$$ 
(\mu \sirc T)(a,b)=\tau \sirc \mu(a,b) \sirc \tau. 
$$ 
 
First, consider the case of $k$ odd. Then there is a pointed homotopy  
$h_t$ between $\tau $ and $1_{S^{2k}}$. Now, the pointed homotopy $h_t \sirc  
\mu(a,b)\sirc h_t$ is a pointed homotopy between $(\mu \sirc T)(a,b)$ and $\mu (a,b)$  
which yields a homotopy between $\mu T$ and $\mu$. 
 
\m Now consider the case of $k$ even. We regard $S^{2k}$ as $\RR^{2k} \cup  
\infty$ with $\RR^{2k}=\{(x_1, \ldots, x_{2k})\}$ and define $\tau', \tau'':  
S^{2k} \to S^{2k}$ by setting 
\begin{eqnarray*} 
\tau'(x_1, x_2, x_3, \ldots, x_{2k})&=&(x_2, x_1, x_3, \ldots, x_{2k}),\\  
\tau''(x_1, \ldots, x_{2k-2}, x_{2k-1},x_{2k} )&=& (x_1, \ldots,  
x_{2k-2},x_{2k},x_{2k-1}),  
\end{eqnarray*} 
(i.e. $\tau'$ permutes the first two coordinates and $\tau''$ permutes the  
last two coordinates). Since $k$ is even, we conclude that $\tau' \simeq \tau  
\simeq \tau''$. Furthermore, $\tau''\tau' \simeq 1_{S^{2k}}$. If we fix such  
pointed homotopies then we get the pointed homotopies 
\begin{eqnarray*} 
(\mu \sirc T)(a,b) &=& \tau \sirc \mu(a,b) \tau \simeq \tau'' \sirc \mu(a,b)  
\tau'=\tau'' \sirc (a \wedge b) \tau'\\ 
&=&\tau'' \sirc (a \wedge 1)\sirc (1\wedge  
b) \tau'=(a \wedge \tau'')\sirc (\tau' \wedge b)\\ 
&=& (a\wedge 1)\sirc (\tau'' \tau') \sirc (1\wedge b) \simeq a\wedge b =\mu(a,b) 
\end{eqnarray*} 
which yield the homotopy $\mu \sirc T \simeq \mu$. 
\qed  
 
\begin{corollary}\label{commut}  
Let $\gf, \psi: \sigma^k \to \sigma^k$ be two automorphisms of $\sigma^k$.  
Then the automorphisms $\gf \dagger \psi$ and $\psi \dagger \gf$ of $\gs^{2k}$ 
are fiberwise homotopic.  
\qed 
\end{corollary} 
 
\m Given two spherical bundles $\xi$ and $\eta$ over $X$, consider the bundle 
$\xi \wedge \eta$ over $X\times X$. We denote by 
$\Delta : X \to X \times X$ the diagonal and consider  
the $\Delta$-adjoint bundle morphism   
$$  
J:=\mathfrak I_{\Delta, \xi \wedge \eta}: \xi \dagger \eta \to \xi \wedge \eta. 
$$ 

\begin{prop} \label{commut-j1} 
For every automorphism $\gf: \eta \to \eta$ the diagram 
$$ 
\CD 
\xi \dagger \eta @>J>> \xi \wedge \eta \\ 
@V1 \dagger \gf VV @VV1 \wedge \gf V\\ 
\xi \dagger \eta @>J>> \xi \wedge \eta  
\endCD 
$$ 
commutes 
\qed 
\end{prop} 
 
\begin{cory}\label{commut-j2} 
The diagram 
$$ 
\CD 
\xi \dagger \eta \dagger \eta @>J>> (\xi \dagger \eta)\wedge \eta @>(1\dagger 
1)\wedge \gf >> (\xi \dagger \eta)\wedge \eta \\  
@| @. @|\\ 
\xi \dagger \eta \dagger \eta @>J>> (\xi \dagger \eta)\wedge \eta @>(1\dagger 
\gf)\wedge 1 >> (\xi \dagger \eta)\wedge \eta   
\endCD 
$$ 
commutes up to homotopy. 
\qed 
\end{cory}

\section[Proof of Theorem 4.5 ]{Proof of Theorem~\ref{t:bn}}\label{BN} 
 
We need certain preliminaries on stable duality~\cite{Spa1}. Given a pointed map 
$f: X\to Y$, let $Sf: SX \to SY$ denote the reduced suspension over $f$. 
So, we have a well-defined map 
\[
S:[X,Y]\bul \to [SX,SY]\bul, \quad [f]\to [Sf]
\] 
 
\begin{prop}\label{freuden} 
Suppose that $\pi_i(Y)=0$ for $i<n$ and that $X$ is a $CW$-space with $\dim X 
< 2n-1$. Then the map $S: [X,Y]\bul \to [SX, SY]\bul$ is a bijection. 
\end{prop} 
 
\p This is the famous Freudenthal Suspension theorem, see e.g. \cite{FFG, H, Spa2, Sw} 
\qed 
 
\m 
Given two pointed spaces $X,Y$, we define $\{X,Y\}$ to be the direct limit of 
the sequence $$ 
\CD 
[X,Y]\bul @>S>> [SX,SY]\bul @>S>> \cdots @>>>[S^nX,S^nY]\bul @>S>>\ldots . 
\endCD 
$$ 
In particular, we have the obvious maps  
$$  
\Sigma: (Y,*)^{(X,*)} \longrightarrow  [X,Y]\bul \longrightarrow  \{X,Y\}. 
$$  
Given a pointed map $f: X \to Y$, the element $\Sigma(f) \in \{X,Y\}$ is called the {\it stable 
homotopy class of $f$}. The standard notation for this one is $\{f\}$, but, as 
usual, in several cases we use the same notation $f$ for $f$, $[f]$ and $\{f\}$.  
 
\m It is well known that, for $n\geqslant 2$, the set $[S^nX,S^nY]\bul$ has a natural 
structure of the abelian group, and the corresponding maps $S$ are 
homomorphisms, \cite{Sw}. So, $\{X,Y\}$ turns out to be a group. Furthermore, 
by Theorem \ref{freuden}, if $X$ is a finite $CW$-space then the map  
$$ 
[S^NX, S^NY]\bul\to \{S^NX,S^NY\} 
$$ is a bijection for $N$ large enough.      
 
\begin{df} 
\rm A map $f: S^d \to A\wedge A^{\perp}$ is called {\it a 
{\rm(}stable{\rm)} $d$-duality} if, for every space $E$,  the maps   
$$ 
u_E: \{A,E\} \to \{S, E \wedge A^{\perp}\} ,\quad 
u_E(\varphi)=(\varphi \wedge 1_{A^{\perp}})u 
$$   
and   
$$ 
u ^E: \{A^{\perp},E\} \to \{S,A \wedge E\},\quad u^E(\varphi)=(1_A \wedge 
\varphi)u 
$$ 
are isomorphisms. 
\end{df} 
 
\begin{propcon}\label{hom} 
Let $u: S^d \to A \wedge A^{\perp}$ be a $d$-duality between two finite  
$CW$-spaces. Then, for all $i$ and $\pi$, the map $u$ yields an isomorphism  
$$ 
H_i(u;\pi): \wt H^i(A^{\perp};\pi) \to \wt H_{d-i}(A,\pi) .
$$ 
\end{propcon} 
 
\p Recall that  
$$ 
H^n(A^{\perp};\pi)=[A, K(\pi,n)]\bul=[S^NA, K(\pi, N+n)]\bul 
$$  
where $K(\pi,i)$ is the Eilenberg--MacLane space. Because of 
Theorem \ref{freuden}, the last group coincides with $\{S^NA, K(\pi, N+n)\}$ 
for $N$ large enough, and therefore  
$$  
H^n(A^{\perp};\pi)=\{S^NA, K(\pi, 
N+n)\} \text{ for $N$ large enough}.  
$$ 
 
Furthermore, let $\eps_n: SK(\pi, n) \to K(\pi, n+1)$ be the adjoint map to the standard homotopy equivalence $K(\pi,n) \to \Omega K(\pi, 
n+1)$, see e.g. \cite{Sw}.  G. Whitehead \cite{Wh} noticed that  
$$  
\wt H_n(A;\pi)=\varinjlim \, [S^{N+n}, K(\pi, N)\wedge A]\bul.  
$$ 
Here $\varinjlim$ is the direct limit of the sequence 
\begin{eqnarray*} 
[S^{N+n}, K(\pi, N)\wedge A]\bul &\longrightarrow &[S^{N+n+1}, SK(\pi, N)\wedge A]\bul \\ 
&\xrightarrow{\ \eps_*\ \ } & [S^{N+n+1}, K(\pi, N)\wedge A]\bul 
\end{eqnarray*} 
(see \cite[Ch 18]{Gray} or \cite[II.3.24]{Rud} for greater details). 
Since $\eps_n$ is an $n$-equivalence, and because of Theorem \ref{freuden}, 
we conclude that  $$ 
\wt H_n(A;\pi)=[S^{N+n}, K(\pi,N) \wedge A]\text{ for $N$ large enough}. 
$$  
So, again because of Theorem \ref{freuden}, 
$$ 
\wt H_n(A;\pi)=\{S^{N+n}, K(\pi,N) \wedge A\} 
$$ 
for $N$ large enough.  
 
\m Now, consider a $d$-duality $u: S^d \to A \wedge A^{\perp}$. Fix $i$ and 
choose $N$ large enough such that  
\begin{eqnarray*} 
\wt H^i(A^{\perp};\pi)&=&\{S^NA^{\perp}, K(\pi, N+i)\}, \\ 
\wt H_{d-i}(A;\pi)  &=& \{S^{N+d}, K(\pi, N+i)\wedge A\}. 
\end{eqnarray*} 
Put $K=K(\pi, N+i)$. By suspending the domain and the range, we get a duality (denoted also by $u$) 
$$ 
u: S^{N+d} \to A \wedge S^NA^{\perp}. 
$$ 
 
This duality yields the isomorphism  
\begin{equation*} u^{K}: \wt H^i(A^{\perp};\pi) =\{S^NA^{\perp}, K\} \to  \{S^{N+d}, K\wedge A\} =\wt H_{d-i}(A;\pi),
\end{equation*} 
and we set $H_i(u;\pi) :=u^{K}$.
\qed 

\forget
\begin{eqnarray*} 
H_i(u;\pi):  &=& u^{K(\pi,N+i)}: \wt H^i(A^{\perp};\pi) =\{S^NA^{\perp}, K(\pi, N+i)\}\\ 
 &\to & \{S^{N+d}, K(\pi, N+i)\wedge A\} =\wt H_{d-i}(A;\pi). 
\end{eqnarray*} 
\forgotten
 
\begin{df}\label{df-cored}\rm 
We dualize \ref{df-red} and say that a pointed map $a: A \to S^k$ (or its 
stable homotopy class $a\in \{A,S^k\}$) is a {\it coreducibility} if the induced 
map $$ a^* : \wt H^i(S^k) \to \wt H^i(A) 
$$ 
is an isomorphism for $i\leqslant k$. 
\end{df} 
 
\begin{prop}\label{red-cored} 
Let $u: S^d \to A \wedge A^{\perp}$ be a $d$-duality between two finite  
$CW$-spaces, and let $k \leqslant d$. A class $\ga\in \{A^{\perp},S^k\}$ is a  
coreducibility if and only if the class $\gb:=u^{S^k}\ga \in \{ S^{d-k},A\}$ is a  
reducibility. 
\end{prop} 
 
\p Let $H_i(u): \wt H^ i(A^{\perp}) \to \wt H_{d-i}(A)$ be the isomorphism as  
in \ref{hom}. Note that the standard homeomorphism $v: S^d \to S^k \wedge S^{d-k}$ 
is a $d$-duality. It is easy to see that the diagram 
$$ 
\CD 
\wt H^ i(A^{\perp}) @>H_i(u)>> \wt H_{d-i}(A)\\  
@A \ga^* AA @ AA\gb_* A\\ 
\wt H^i(S^k) @> H_i(v) >> \wt H_{d-i}(S^{d-k}) 
\endCD 
$$ 
commutes. In particular, the left vertical arrow is an isomorphism if and
only if the   right one is. 
\qed  
  
\m Consider a closed connected $n$-dimensional PL manifold $M$ and embed it 
in $\RR^{N+n+k}$ with $N$ large enough. Let $\iota: S^{N+n+k}\to T\nu^{N+k}$ 
be a collapse map as in \ref{collapse}, and let 
$$ 
J: (\nu^{N+k})\bul = (\nu^N)\bul \dagger \gs^k \longrightarrow (\nu^N)\bul 
\wedge \gs^k   
$$  
be the morphism as in \ref{commut-j1},  where $\gs^k=\gs^k_M$. 
 
\begin{theorem}\label{N-duality} 
The map  
$$ 
\CD 
S^{N+n+k} @>\iota>> T\nu^{N+k} @>TJ>> T\nu^N \wedge \sigma^k 
\endCD 
$$ 
is an $(N+n+k)$-duality map. 
\end{theorem} 
 
\p This is actually proved in \cite{DP}. For greater detail, see 
\cite[V.2.3(i)]{Rud}. \qed 
 
\begin{cory}\label{reducibility} 
 If the manifold $M$ is orientable then the collapse map $\iota: S^{N+n} \to T\nu^N$ is a reducibility.  
\end{cory} 
 
\p Recall that $T\gs^k=(M\times S^k)/M=S^k(M^+)$. Consider a surjective map 
$e:M^+ \to S^0$ and define $\eps=S^ke: T\gs^k \to S^k$. The map  
$$ 
S^k\iota: S^kS^{N+n}=S^{N+n+k} \to T\nu^{N+k}=S^kT\nu 
$$ 
can be written as 
\begin{equation*} 
S^{N+n+k} \xrightarrow{S^k\iota} T\nu^{N+k}=T\nu\sp N \wedge T \gs^k 
\xrightarrow{1\wedge \eps} T\nu^N \wedge S^k=S^kT\nu,  
\end{equation*} 
where the composition of first two maps is the duality from~\ref{N-duality}.
Hence $\iota$ is dual to $\eps$ with respect to duality~\eqref{D}. 
Clearly, $\eps$ is a coreducibility because $M$ is orientable. Thus, the result follows from 
\ref{red-cored}. \qed 
 
\m 
For technical reasons, it will be convenient for us to consider the duality  
\begin{equation}\label{duality-k} 
S^{N+n+2k} \to T\nu^{N+2k} \xrightarrow{TJ} T\nu^{N+k}\wedge T\gs^k. 
\end{equation} 

This duality yields an isomorphism 
\begin{equation}\label{D} 
\begin{aligned}
D:=u^{S^k}:\{T\gs^k, S^k\} &\to \{S^{N+n+2k}, S^k\wedge 
T\nu^{N+k}\}\\
&=\{S^{N+n+k}, T\nu^{N+k}\}.  
\end{aligned}
\end{equation} 
 
\begin{prop}\label{dual-thom} 
For every  automorphism $\gf: \gs^k \to \gs^k$ the following diagram commutes 
up to homotopy:  
$$ 
\CD 
S^{N+n+2k} @>\iota >> T\nu^{N+2k} @>TJ>> T\nu^{N+k} \wedge T\sigma^k 
@>T(1\dagger \gf)\wedge 1>> T\nu^{N+k} \wedge T\sigma^k\\ @| @| @. @|\\ 
S^{N+n+2k} @>>> T\nu^{N+2k} @>TJ>> T\nu^{N+k} \wedge T\sigma^k @>1\wedge T\gf 
>> T\nu^{N+k} \wedge T\sigma^k  
\endCD 
$$ 
\end{prop} 
 
\p This follows from \ref{commut-j2} .
\qed 
 
\m Every automorphism $\gf: \gs^k \to \gs^k$ yields a homotopy equivalence 
$$ 
\CD 
T(1\dagger\gf): T\nu^{N+k}=T(\nu^N \dagger \gs^k)@>>> T(\nu^N \dagger \gs^k)= 
T\nu^{N+k} \endCD 
$$ 
and hence an isomorphism  
$$ 
T(1\dagger\gf)_*: \{S^{N+n+k},T\nu^{N+k}\}\to  \{S^{N+n+k},T\nu^{N+k}\}. 
$$  
So, we have the $\aut \gs^k$-action  
\begin{eqnarray*} 
a_{\nu}: \aut \sigma^k \times  \{S^{N+n+k},T\nu^{N+k}\} &\to & 
 \{S^{N+n+k},T\nu^{N+k}\}, \\  
a_{\nu}(\gf, \ga) &=& T(1\dagger\gf)_*(\ga). 
\end{eqnarray*} 
Similarly, every automorphism $\gf$ of $\gs^k$ induces a homotopy equivalence 
$T\gs^k \to T\gs^k$, and therefore we have the action $$ 
a_{\gs}: \aut \gs^k \times \{T\gs^k, S^k\} \to   \{T\gs^k, S^k\}. 
$$ 
 
\begin{theorem}\label{actions}  
The diagram 
\begin{equation*} 
\CD 
\aut \sigma^k \times \{T\gs^k, S^k\} @>a_{\sigma}>> \{T\gs^k, S^k\} \\ 
@V1\times D VV @ VV D V\\ 
\aut \sigma^k \times \{S^{N+n+k},T\nu ^{N+k}\} @>a_{\nu}>>  
\{S^{N+n+k},T\nu ^{N+k}\} 
\endCD 
\end{equation*} 
commutes.  
\end{theorem} 
 
\p This follows from \ref{dual-thom} and the definition of $D, a_{\nu}$ and $a_{\gs}$. 
\qed 
 
\m Because of Theorem \ref{freuden}, for $k$ large enough we have   
$$ 
\{T\gs^k, S^k\}=\pi^k(T\gs^k) \text{ and } \{S^{N+n+k}, T\nu^{N+k}\}=\pi_{N+n+k}(T\nu^{N+k}). 
$$  
Then we can rewrite the diagram from Theorem \ref{actions} as 
\begin{equation}\label{non-stab} 
\CD 
\aut \sigma^k \times \pi^k\left(T\gs^k\right) @>a_{\sigma}>> \pi^k\left(T\gs^k\right)\\ 
@V1\times D VV @ VV D V\\ 
\aut \sigma^k \times \pi_{N+n+k}\left(T\nu ^{N+k}\right) @>a_{\nu}>>  
\pi_{N+n+k}\left(T\nu ^{N+k}\right) 
\endCD 
\end{equation} 

\m Let ${\C R} \in \pi_{N+n+k}(T\nu^{N+k})$ be the set of reducibilities, and let  
${\C C} \in \pi^k(T\gs^k)$ be the set of coreducibilities. Then, clearly, $a_{\nu}(\C  
R)\subset \C R$ and $a_{\sigma}(\C C)\subset \C C$. Therefore, in view of 
\propref{red-cored}, the diagram \eqref{non-stab} yields the 
commutative diagram 
\begin{equation} 
\CD 
\aut \sigma^k \times \C C @>a_{\sigma}>> \C C \\ 
@V1\times D VV @ VV D V\\ 
\aut \sigma^k \times\C R @>a_{\nu}>>\C R 
\endCD 
\end{equation} 
 
\m 
\begin{theorem}\label{bn-prefinish} 
For every $\ga, \gb \in \C C$ there exists an automorphism $\gf$ of 
$\sigma^k$ such that  $a_{\sigma}(\gf, \ga)=\gb$. Moreover, this $\gf$ is 
unique up to fiberwise homotopy. In other words, the action $a_{\sigma}: 
\aut  \sigma^k \times \C C \to \C C$ is free and transitive. 
\end{theorem} 
 
\p Recall that $T\gs^k = (M\times S^k)/M$. So, for every $m\in M$, a pointed  
map $f: T\gs^k \to S^k$ yields a pointed map $f_m: S^k_m \to S^k$ where 
$S^k_m$ is the fiber over $m$. Furthermore, $f$ represents a coreducibility if 
and only if all maps $f_m$ belong to $F_k$. In other words, every coreducibility $T\gs^k \to S^k$ yields a homotopy class $M \to F_k]$, and in fact we have  a bijectioj $\C C \to  [M, F_k]$ there.  Moreover, it is 
easy to see that, in view of \propref{aut}, the action $a_{\gs}$ 
coincides with the map $$ [M, F_k]\times  [M, F_k] \to [M, F_k] $$ induced by 
the product in $F_k$, and the result follows. \qed 
 
Since $D$ is an isomorphism, Theorem \ref{bn-prefinish} yields the following 
corollary.  
 
\begin{cory}\label{bn-finish} 
The action $a_{\nu}: \aut  
\sigma^k \times \C R \to \C R$ is free and transitive.  
\qed 
\end{cory} 
 
\m Now we can finish the proof of Theorem \ref{t:bn}. Assuming $\dim 
\eta=N+k$ to be large enough, we conclude that $\nu\bul$ and $\eta\bul$ are 
homotopy equivalent over $M$, see Atiyah~\cite[Prop. 3.5]{At}. (Note that 
Atiyah works with non-sectioned bundles, but there is no problem to adapt the 
proof for sectioned ones.) Choose any such $F_{N+k}$-equivalence $\gf: \eta\bul 
\to \nu\bul$ and consider the induced homotopy equivalence $T\gf: T\eta \to 
T\nu$. Clearly, the composition   
$$   
\CD  
\gb: S^{N+n+k}  @>\ga >> T\eta @> T\gf >> T\nu   
\endCD  
$$  
is a reducibility. So, by \ref{bn-finish}, there exists an $F_{N+k}$-equivalence 
$\gl: \nu\bul \to \eta\bul$ over $M$ with $(T\gl)_*(\gb)=\iota$. Now, we 
define $\mu : \nu\bul \to \eta\bul$ to be the fiber homotopy inverse to 
$\gl\gf$. (The existence of an inverse equivalence can be proved following 
Dold~\cite{Dold}, cf. also~\cite{May}). Clearly, $\mu_*\iota=\ga$. This proves the 
existence of the required equivalence $\mu$. 
 
\m Furthermore, if there exists another equivalence $\mu': \eta\bul \to 
\nu\bul$, then $\mu'\sirc \mu^{-1}(\iota)=\iota$, and so $\mu$ and $\mu'$ are 
homotopic over $M$. This proves the uniqueness of $\mu$. Thus, 
Theorem \ref{t:bn} is proved.    
\qed

 \section{Normal Morphisms, Normal bordisms,  and $F/PL$} 
 
Throughout the section we fix a closed orientable  $n$-dimensional PL manifold $M$. 
 
\begin{df}[\cite{Br}] \rm  \label{d:normmap}
A {\it normal morphism at 
$M$} is a PL $\RR^N$-morphism 
$\gf: \nu_V\to \xi$ where $\xi$ is a PL $\RR^N$-bundle over $M$, $V$ is a closed PL manifold,  and $\nu_V$ is PL $\RR^N$-bundle of $V$ in $\RR^{N+n}$.

\begin{example}\label{ex:norm}\rm
Let $h: V \to M$ be a homotopy equivalence and $g: M \to V$ a homotopy inverse map to $h$. Consider the normal bundle $\nu$ of $V$ and put $\xi:=g^*(\nu)$. Then $h^*(\xi)=h^*g^*\nu=\nu$. The correcting morphism $\nu =h^*(\xi) \to \xi$ is a normal morphism.
\end{example}

\m A normal morphism is called {\it reducible} if the map 
$$ 
\CD 
S^{N+n} @>\text{collapse}>> T\nu_V @>T\gf >> T\xi 
\endCD 
$$ 
is a reducibility. 
\end{df} 
 
\m Because of the Thom Isomorphism Theorem, a normal morphism is reducible whenever  
its base $V\to M$ is a map of degree 1.  (One can prove that $\xi$ is orientable if $V$ and
$M$ are.)  

\m We denote the set of all reducible normal morphisms at $M$ by $\Nor(M)$. (For persons who aks whether $\Nor(M)$ is a set, we take the space $\RR^{\infty}$ and assume that all the spaces in \ref{d:normmap} are contained in $\RR^{\infty}$.)

\begin{construdef}\rm \label{condef:ass}
Represent a map (homotopy class) $f: M \to F/PL$ by an 
$(S^N,*)$-morphism $\gf: \nu_M\bul \to (\gga^N_{PL})\bul$ with $N$ large enough, 
see \eqref{fpl-stab}. Set $\xi=(\bs \gf)^*\gga^N_{PL}$. Then the correcting $(S^N,*)$-morphism $\nu_M\bul \to 
\xi\bul$ is a commutative diagram  
\begin{equation}\label{uu}  
\CD U\bul @>g>> {U'}\bul\\ @VqVV @VVpV\\ 
M @= M 
\endCD 
\end{equation} 
where $\nu_M=\{q: U \to M\}$, $\xi=\{p: U' \to M\}$ are PL $\RR^N$-bundles, and $U\bul, {U'}\bul$ are 
fiberwise one-point compactifications of $U$ and $U'$, respectively. 
 
We regard $M$ as the zero section of $\xi$, $M \subset U'$ and deform 
$g$ to a map $t: U\bul \to {U'}\bul$ which is transverse to $M$. Set 
$V=t^{-1}(M)$ and $b=t|_V:V \to M$. We can assume that $V\subset U$. So, we get the $b$-adjoint PL $\RR^N$-morphism 
\begin{equation}\label{ass}
\mathfrak I_b: b^*\xi \to \xi, \quad \bs(\mathfrak I_b)=b: V \to M. 
\end{equation}
Note that $b^*\xi$ is the normal bundle of $V$ in $U$, and therefore it is the normal bundle $\nu_V$ ot $V$ in $\RR^{N+n}$ because $U$ is the open subset of $\RR^{N+n}$. In other words, the morphism \eqref{ass} 
is a normal morphism at  $M$. We say that the normal morphism $\eqref{ass}$ is {\it 
associated with a map {\rm(}homotopy class{\rm)} $f: M \to F/PL$}. 
 
Clearly, there are many normal morphisms that are associated with a given map $f: M \to F/PL$. 
\end{construdef} 
 
\begin{construdef}\rm \label{condef:norm-ass}
Let  
\begin{equation}\label{norm-ass} 
\gf: \nu_V \to \xi  
\end{equation} 
be a reducible normal morphism at $M$ and assume that $\dim \nu_V$ is large. 
Consider a collapse map (homotopy class) $\iota: S^{N+n} \to T\nu_M$ as in 
\ref{collapse}. Since the map   
$$ 
\CD  
\ga: S^{N+n} @>\text{collapse}>> T\nu_V @>T\gf >> T\xi  
\endCD 
$$ 
is a reducibility, there exists, by Theorem \ref{t:bn}, a unique $(S^N,*)$-morphism 
$\mu: \nu_M\bul \to \xi\bul$ with $\mu_*(\iota)=\ga$. Now, the morphism  
\begin{equation*} 
\CD 
\nu_M\bul @>\mu_* >> \xi\bul @>\text{classif} >> (\gga^N_{PL})\bul 
\endCD 
\end{equation*} 
is a homotopy PL structuralization of $\nu_M$. Thus, in view of~\eqref{fpl-stab}, we get a homotopy class  
in $[M, F/PL]$. We denote by  $f_{\gf}: M \to F/PL$ any representative of this class.  
\end{construdef} 
                               
\begin{prop}\label{p:ass} 
The normal morphism $\eqref{norm-ass}$ is associated with the map $f_{\gf}: M \to 
F/PL$.  
\qed 
\end{prop} 
 
\begin{df}\label{d:nor}
 The function 
 \begin{equation}\label{eq:ass}  
\Gamma=\Gamma_M:  \Nor(M)\longrightarrow [M, F/PL], \quad \gf \mapsto [f_{\gf}].
 \end{equation}
is called the {\it normal invariant} for $M$.
\end{df}.

Probably, a reader noticed that we already defined normal invariant $\jf$ in \defref{d:norminv}. Now we show that these two definitions (for homotopy structures and for normal morphisms) are quite close to each other. The relation between $\Gamma$ and $\jf$ appears in the commutative diagram
 \[
 \CD
{\text{ he}}(M) @>>>\Nor M \\
@VVV @VV \Gamma V  \\
\spl(M)@> \jf >> [M, F/PL]
\endCD
 \]
  where he$(M)$ is the set of homotopy equivalences whose targets is $M$. The horizontal top map is explained in \exref{ex:norm}, the left vertical map send a homotopy equivalence to its equivalence class as \defref{d:hstruct}.
  
 \begin{df}[\cite{Br}] \rm \label{df:normbord}
 A {\it normal bordism} between two normal morphisms $\gf_i: \nu_{V_i}\to M, i=0,1$ at $M$ is a PL $\RR^N$-morphism $\Phi: \nu_W\to \xi$ where $W$ is a compact PL manifold with $\pa W=V_0\sqcup V_1$ and $\Phi|_{V_i}=\xi, i=0,1$. Furthermore, $\nu_M$ is the PL normal $\RR^N$-bundle of $W$. 
 \[
 \CD
 D @>\wh c >> E'\\ 
@VVV @VVV\\ 
W @>c>> M 
 \endCD
 \]
 where $D\to W$ is a normal bundle of $W$ and $E\to M$ is the same bundle as in \defref{d:normmap}. Furthermore, $W$ is a compact manifold with $\pa W=V_0\sqcup V_1$ and $(\wh c, c)|_{V_i}=(\wh b_i, b_i), i=0,1$. 
 \end{df}
 
 We say that two normal morphisms are {\it normally bordant} if there exists a normal bordism between these two normal morphisms. Clearly, ``to be normally bordant'' is an equivalence relation. The equivalence classes are called the {\it normal bordism classes}.  We denote by $[\Nor M]$ the set of normal bordism classes at $M$.
 
 \begin{thm}\label{t:normbord}
If $\gf_0, \gf_1$ are two normally bordant normal morphisms at $M$ then $f_{\gf_0}\cong f_{\gf_1}: \Nor M \to F/PL$. So, the map $\Gamma$ yields a map
\[
\wt \Gamma: [\Nor M]\to F/PL, \quad \wt \Gamma[\gf]=[\Gamma (\gf)].
\]
Moreover, the map $\wt \Gamma$ is a bijection;
\end{thm}

\m
{\it Proof} (sketch).  This is a version of the Pontryagin--Thom theorem. We give a sketch and leave the detail to the reader. Let $\Phi$ be a normal bordism as in \defref{df:normbord}. Follow \ref{condef:norm-ass} and construct a map $F_{\phi}: M\to F/PL$. Then $F_{\Phi}|_{M\times \{i\}}=f_{\gf_i}$ for $i=0,1$. So, the above mentioned map $\wt\Gamma$ is well-defined.

To construct an inverse map $\Delta$ to $\wt \Gamma$, take a map $f: M\to F/PL$ and put $\Delta[f]$ to be the normal bordism class that is associated to $f$. Check that this normal  bordism class is well-defined and that $\wt \Gamma$ and $\Delta$ are inverse to each other.
\qed

\m Recall that a closed manifold is called {\it almost parallelizable} if it  
becomes parallelizable after deleting of a point. Note that every almost  
parallelizable manifold is orientable (e.g., because its first 
Stiefel--Whitney  class is equal to zero). 
 
\begin{prop}\label{real} 
For every $V^m$ be an almost parallelizable PL manifold $V^m$ there exists a reduciible normal morphism with a base $V\to M$.
\end{prop}
 
\p We regard $S^m=\{(x_1, \ldots, x_{m+1}) \bigm| \sum x_i^2=1\}$ as the union of  
two discs, $S^m = D_+\cup D_-$, where 
$$ 
D_+= \{x\in S^m|x_{m+1} \geqslant 0\}, \quad D_-= \{x\in  
S^m|x_{m+1} \leqslant 0\}. 
$$  
Take a map $b: V \to M$ of degree 1. We can assume that there is a small closed 
disk $D_0$ in $V$ such that  $b_+:=b|_{D_0} :D_0 \to D_+$ is a PL homeomorphism. We 
set $W=V\setminus  (\operatorname {Int} D_0)$. Since $W$ is parallelizable, 
there exists a PL morphism  $\gf_-: \nu_V|_W \to \theta_{D_-}$  such 
that $b|_W: W \to D_-$ is the base of  $\gf$. Furthermore, since $b_+$ is a PL 
homeomorphism, there exists a morphism  $\gf_+: \nu_V|_{D_0} \to \theta |_{D_+}$ over 
$b_+$ such that $\gf_+$ and $\gf_-$  coincide over $b|_{\partial W}: \partial W 
\to S^{m-1}$. Together $\gf_+$ and  $\gf_-$ give us a PL morphism $\gf: \nu_V \to 
\xi$ where $\xi$ is a PL bundle over  $S^m$. Clearly, $\gf$ is a normal morphism 
with the base $b$, and it is reducible  because $\deg b=1$. 
\qed

\section {The Sullivan Map $s:[M, F/PL] \to P_{\dim M}$} 
 
We define the groups $P_i$ by setting 
\begin{equation*} 
P_i=\left\{ 
\begin{array}{lcl} 
\mathbb Z &&\text{if $i=4k$},\\ 
\mathbb Z/2 &&\text{if $i=4k+2$},\\ 
0 &&\text{if $i=2k+1$} 
\end{array} 
\right. 
\end{equation*} 
where $k \in \mathbb N$. 
 
\m Given a closed connected $n$-dimensional PL manifold $M$ (which is 
assumed to be orientable for $n=4k$), we define a map 
\begin{equation}\label{sul} 
s: [M,F/PL] \to P_n 
\end{equation} 
as follows. Given a homotopy class $f: M \to F/PL$, consider a normal morphism 
\[
\mathfrak I_b: b^*\xi \to \xi, \quad \bs(\mathfrak I)=b: V \to M
\]
associated with $f$, see \eqref{ass}. 
 
For $n=4k$, let $\psi$ be the symmetric bilinear intersection form on 
$$ 
\Ker\{b_*: H_{2k}(V; \mathbb Q) \to H_{2k}(M;\mathbb Q)\}. 
$$ 
We define $s(u)=\sigma(\psi)/ 8$ where $\sigma(\psi)$ is the 
signature of $\psi$. It is well known that $\gs(\psi)$ is divisible by 8, 
(see e.g.~\cite{Br}), and so $s(u)\in \ZZ$. 
 
Also, it is easy to see that  $\gs(\psi)= \gs(V)-\gs(M)$, and so 
$$ 
s(u)=\frac{\gs(V)-\gs(M)}8 
$$ 
where $\gs(V), \gs(M)$ is the signature of the manifold 
$V, M$, respectively.    
 
\m For $n=4k+2$, we define $s(u)$ to be the Kervaire 
invariant of the normal morphism $\II_b$, see e.g.~\cite{Br}. 

\m The routine arguments 
show that $s$ is well-defined, i.e. it does not depend on the choice of the 
associated normal morphism. See \cite[Ch. III, \S 4]{Br} or \cite{N1} for 
details. 
 
\m In particular, if $b$ is a homotopy equivalence then $s(u)=0$. 
 
\m One can prove that, for all $M$, the map $s$ is a homomorphism of abelian groups, where the abelian group 
structure on $[M,F/PL]$ is given by the $H$-space structure on $F/PL$.   
 
\m Given a map $f: M \to F/PL$, it is useful 
to introduce the notation $s(M,f):=s([f])$ where $[f]$ is the homotopy 
class of $f$. 

\begin{theorem}\label{image}  
{\rm (i)} The map $s: [S^{4i}, F/PL] \to \ZZ$ is surjective for all $i>1$,  
  
\par {\rm (ii)} The map $s: [S^{4i-2}, F/PL] \to \ZZ/2$ is surjective for all
$i>0$ .  

\par {\rm (iii)}  The image of the map $s: [S^4, F/PL] \to \ZZ$ is
the  subgroup of index $2$.   
\end{theorem} 
  
\p (i) For every $k>1$ Milnor constructed a parallelizable $4k$-dimensional  
smooth manifold $W^{4k}$ of signature 8 and such that $\partial W$ is a homotopy 
sphere, see \cite[V.2.9]{Br}. Since, by Theorem \ref{Sma}, every 
homotopy sphere of dimension $\geqslant 5$ is PL homeomorphic to the standard one, we 
can form a closed PL manifold   
$$  
V:=W\cup_{S^{4k-1}} D^{4k} 
$$ 
of the signature 8. Because of Proposition \ref{real}, there exists a  
reducible normal morphism with the base $V^{4k} \to S^{4k}$. Because of Proposition 
\ref{p:ass}, this normal morphism is associated with a certain map (homotopy 
class) $f: S^{4k} \to F/PL$. Thus, 
\[
s(S^{4k},f)=\frac{\gs(V^{4k})-gs(S^{4k})}{8}=1.
\] 
 
\m (ii) The proof is similar to that of (i), but we must use 
(4k+2)-dimensional parallelizable Kervaire manifolds $W$, $\partial W= 
S^{4k+1}$ of the Kervaire invariant one, see \cite[V.2.11]{Br}. 
 
\m (iii) The Kummer algebraic surface \cite{K2} gives us an example of  
4-dimensional almost parallelizable smooth manifold of the signature 16. So,  
$\IM s \supset 2\ZZ$.  
 
Now suppose that there exists $f: S^4 \to F/PL$ with $s(S^4,f)=1$. Then  
there exists a normal morphism with the base $V^4 \to S^4$ and such that $V$ has  
signature 8. Since normal bundle of $V$ is induced from a bundle over $S^4$, 
we  conclude that  $w_1(V)=0=w_2(V)$. But this contradicts the Rokhlin Theorem 
\ref{rohlin}. 
\qed 
 
\begin{theorem}[Sullivan~\cite{Sul2}]\label{bns} 
For any closed simply-connected PL manifold $M$ of dimension $\geqslant 5$, the 
sequence 
$$ 
\CD 
0 @>>> \spl (M) @>\jf >> [M,F/PL] @>s>> P_{\dim M} 
\endCD 
$$ 
is exact, i.e. $\jf$ is injective and $\IM \jf =s^{-1}(0)$. 
\end{theorem} 
 
\p See~\cite[II.4.10 and II.4.11]{Br}. Note that the map $\omega$ in 
loc. cit. is the zero map because, by Theorem \ref{Sma}, 
every homotopy sphere of dimension $\geqslant 5$ is PL homeomorphic to the standard 
sphere.  
\qed

\begin{corollary} 
\label{c:gpl} 
 We have $\pi_{4i}(F/PL)=\mathbb Z$, $\pi_{4i-2}(F/PL)=\mathbb Z/2$, and 
$\pi_{2i-1}(F/PL)=0$ for every $i>0$. Moreover, the homomorphism 
$$ 
s: [S^k, F/PL] \to P_k 
$$  
is an isomorphism for $k\ne 4$, while for $k=4$ it has the form  
$$ 
\CD 
\ZZ =\pi_4(F/PL) @>s>> P_4=\ZZ, \quad a \mapsto 2a. 
\endCD 
$$   
\end{corollary} 
 
\p First, if $k>4$ then, because of the Smale Theorem \ref{Sma}, 
$\spl(S^k)$ is the one-point set. Now the result follows from \theoref{bns} and 
\ref{image}. 
 
 If $k\leqslant 4$ then $\pi_k(PL/O)=0$, cf. Remark \ref{pl/o}. So, 
$\pi_k(F/PL)=\pi_k(F/O)$. Moreover, the forgetful map $\pi_k(BO) \to 
\pi_k(BF)$ coincides with the Whitehead $J$-homomorphism. So, we have the long 
exact sequence  
$$  
\CD \cdots \to \pi_k(F/O) \to \pi_k(BO) @>J>> \pi_k(BF) \to 
\pi_{k-1}(F/O) \to \cdots .  
\endCD 
$$ 
For $k\leqslant 5$ all the groups $\pi_k(BO)$ and $\pi_k(BF)$ are known (note 
that $\pi_k(BF)$ is the stable homotopy group $\pi_{k+N-1}(S^N)$), and it is 
also known that $J$ is an epimorphism for $k=1,2,4,5$, see e.g. \cite{Adams}. 
Thus, $\pi_k(F/O)\cong P_i$ for $k\leqslant 4$. 

The last claim follows from \theoref{image}. 
\qed

\section{The Homotopy Type of $F/PL[2]$} 
 
Recall that, given a space $X$ and an abelian group $\pi$, we allow us to ignore the distinction between lements of $H^n(X;\pi)$ and maps (homotopy classes) $X \to K(\pi,n)$.

 \forget For example, regarding a Steenrod cohomology operation 
$Sq^2$ as an element $Sq^2\in H^{k+2}(K(\ZZ/2, k);\ZZ/2)$, we can treat  it
as a map $Sq^2: K(\ZZ/2,k) \to K(\ZZ/2,k+2)$. 

\forgotten
 
\begin{notation}\label{loc}\rm
 Given a prime $p$, let $\ZZ[p]$ be the subring of $\QQ$ consisting 
 of all irreducible fractions with denominators relatively prime 
to $p$, and let $\ZZ[1/p]$ be the subgroup of $\QQ$ consisting of the 
fractions $m/p^k, m\in \ZZ$. Given a simply-connected space $X$, we denote by 
$X[p]$ and $X[1/p]$ the $\ZZ[p]$- and $\ZZ[1/p]$-localization of $X$, 
respectively. Furthermore, we denote by $X[0]$ the $\QQ$\,-localization of $X$. 
For the definitions, see \cite{HMR}. 
\end{notation}  

\begin{prop}[Sullivan \cite{Sul1, Sul2}]\label{t:cl}
For every $i>0$ there are cohomology classes 
$$
K_{4i}\in H^{4i}(F/PL;\ZZ[2]), \  K_{4i-2}\in H^{4i-2}(F/PL;\ZZ/2)
$$
such that
$$
s(M^{4i},f) = \langle f^*K_{4i},\, [M] \rangle
$$
for every closed connected oriented PL manifold $M$,
and                                    
$$
s(N^{4i-2},f) = \langle f^*K_{4i-2},\, [N]_2 \rangle.
$$
for every closed connected manifold $N$.
Here $[M]\in H^{4i}(M)$ is the fundamental class of $M$, $[N]_2\in 
H^{4i-2}(N;\ZZ/2)$ is the modulo $2$ fundamental class of $N$, and $\langle 
-,-\rangle$ is the Kronecker pairing. 
\end{prop} 
 
\p Let $MSO_*(-)$ denote the oriented bordism theory, see e.g. \cite{Rud}. 
Recall that if two maps $f: M^{4i} \to F/PL$ and  $g: N^{4i} \to F/PL$ are 
bordant (as oriented singular manifolds) then $s(M,f)=s(N,g)$. Thus, $s$ 
defines a homomorphism  
$$  
\wt s: MSO_{4i}(F/PL) \to \ZZ. 
$$ 
 
\m It is well known that the Steenrod--Thom map   
$$ 
t: MSO_*(-)\otimes \ZZ[2] \to H_*(-; \ZZ[2]) 
$$  
splits, i.e. there is a natural map $v: H_*(-; \ZZ[2]) \to MSO_*(-)\otimes 
\ZZ[2]$ such that $tv=1$ (a theorem of Wall \cite{Wall1}, see also
\cite{Stong, Rud, Astey}. In  particular, we have a natural homomorphism   
$$  
\CD  
\wh s:  H_{4i}(F/PL; \ZZ[2]) @>v>>  MSO_{4i}(F/PL)\otimes 
\ZZ[2] @>\wt s >> \ZZ.  
\endCD 
$$ 
Since the evaluation map 
\begin{equation*}
\begin{aligned}
\ev: H^*(X; \ZZ[2]) &\to \Hom (H_*(X;\ZZ[2]),\ZZ[2]), \\
(\ev(u)(v) =\langle u,v\rangle, \quad & u\in  H^*(X; \ZZ[2]), \quad u\in   H_*(X;\ZZ[2]) 
\end{aligned}
\end{equation*}
is surjective for all $X$, there exists a class $K_{4i}\in H^{4i}(F/PL; 
\ZZ[2])$ such that $\ev(K_{4i})=\wh s$. Now $$ 
s(M,f)=\wh s(f_*[M]) =\langle K_{4i}, f_*[M]\rangle = \langle f^*K_{4i}, [M]\rangle. 
$$ 
So, we constructed the desired classes $K_{4i}$. 
 
\m The construction of classes $K_{4i-2}$ is similar. Let $MO_*(-)$ denoted 
the non-oriented bordism theory. Then the map $s$ yields a homomorphism  
$$ 
\wt s: MO_{4i-2}(F/PL) \to \ZZ/2. 
$$ 
Furthermore, there exists a natural map $H_*(-;\ZZ/2) \to MO_*(-)$ which 
splits the Steenrod--Thom homomorphism, and so we have a homomorphism $$ 
\CD 
\wh s:  H_{4i-2}(F/PL; \ZZ/2) @>>>  MO_{4i-2}(F/PL)\otimes \ZZ[2] @>\wt s >> \ZZ/2 
\endCD 
$$ 
with $\wh s(f_*([M]_2)=s(M,f)$. Now we can complete the proof similarly to 
the case of classes $K_{4i}$.  \qed  
 
\m We set  
\begin{equation}\label{pi}
\Pi:=\prod_{i>1}\left( K(\mathbb Z[2], 4i)\times K(\mathbb Z/2, 
4i-2)\right). 
\end{equation} 
Together the classes $K_{4i}: F/PL \to K(\ZZ[2],4i), i>1$ and $K_{4i-2}: F/PL 
\to K(\ZZ[/2,4i-2), i>1$ yield a map 
\begin{equation} \label{bigpi}
K: F/PL \to \Pi
\end{equation}
 such that for each $i>1$ the map
\begin{equation*}
 \CD
F/PL @>K>> \Pi @>\text{ projection}>>  K(\ZZ[2],4i)
 \endCD
 \end{equation*}
 coincides with $K_{4i}$ and the map
 \begin{equation*}
 \CD
F/PL @>K>> \Pi @>\text{ projection}>>  K(\ZZ/2,4i-2)
 \endCD
 \end{equation*}
  coincides with $K_{4i-2}$.
  
\begin{lemma}\label{l:gpl} 
The map 
\begin{equation}\label{k2} 
K[2]: F/PL[2] \to \Pi 
\end{equation} 
induced an isomorphism of homotopy groups in dimensions $\geq 5$. 
\end{lemma} 
 
\p This follows from \theoref{image} and \corref{c:gpl}. 
\qed 
 
 \m Let $Y$ be the  Postnikov 4-stage of  $F/PL$. So, we have a map 
 \begin{equation}\label{psi}
\psi: F/PL \lr Y 
\end{equation}
that induces an isomorphism of homotopy groups in dimension $\leqslant 4$.  Consider the map 
\[
\phi: F/PL[2] \to Y[2] \times \Pi, \quad \phi(x)=(\psi[2](x), K[2](x)).
\]
 
\begin{theorem}\label{t:gpl} 
The map  
\begin{equation*}
\phi: F/PL[2] \to Y[2] \times \Pi 
\end{equation*} 
is a homotopy equivalence. 
\end{theorem} 
 
\p The maps 
\[
\phi_*: \pi_i(F/PL[2] \to \pi_i(Y[2]\times\Pi)
\]
ari isomorphism for all $i$. Indeed, for $i\leqslant 4$ the holds since $\psi$ is the Postnikov 4-approximation of $F/PL$,  for $\geqslant 5$ it  follows from $\ref{k2}$. Thus, $\phi$ is a homotopy equivalence by the Whitehead Theorem.  
\qed 

\m Now we discuss the space $Y$ in greater detail. We have $\pi_2(Y)=\ZZ/2$, $\pi_4(Y)=\ZZ$, and $\pi_i(Y)=0$ otherwise. So, we have a $K(\ZZ,4)$-fibration
\begin{equation}\label{y-fib}
\CD
K(\ZZ,4) @>i>> Y @>p>> K(\ZZ/2.2)
\endCD
\end{equation}
whose characteristic class is the Postnikov invariant 
\[
\varkappa\in H^5(K(\ZZ/2, 2)):=\varkappa\in H^5(K(\ZZ/2, 2); \ZZ)
\]
of $Y$. We shall see in \theoref{postn}  below that $\varkappa=\gd Sq^2\iota_2\ne 0$. Hence, $\varkappa$ is also the first non-trivial Postnikov invariant of $F/PL$.

\begin{lemma}\label{cp2}\rm 
There exists a map $g:\CC \PP^2\to F/PL$ such that $s(\CC \PP^2,g)=1$.
\end{lemma}
  
\p Let $\eta$ denote the canonical 
complex line bundle over $\CC \PP^2$. First, we prove that $24\eta$ is fiberwise 
homotopy trivial. Let $H: S^3\to S2$ be the Hopf map. Consider the Puppe sequence
\[
\CD
S^3 @>H>> S^2@>>> \CC \PP^2 @>>> S^4
\endCD
\]
and the induced exact sequence 
\[
\CD
[S^3, BF]  @<H^*<< [S^2, BF]@<<< [\CC \PP^2, BF] @<<< [S^4,BF] 
\endCD
\]
 Let $\Pi_n^S$ denotes the $n$-th stable homotopy group $\pi_{n+N}(S^N)$, $N$ large. Note that $[S^n, BF]=\Pi_{n-1}^S$. We have $[S^4, BF]=\Pi_3^S=\ZZ/24$,~\cite{H}. Furthermore, the homomorphism
\[
\CD
\ZZ/2=[S^3, BF]  @<H^*<< [S^2, BF]=\ZZ/2 
\endCD
\]
is an isomorphism, because the suspension $S^N H:S^{N+3} \to S^{N+2}$ is the generator of $\pi_{N+3}S^{N+2}=\Pi_n=\ZZ/2$ for $N$ large. Hence, $[\cp2, BF]$ is a quotient group of $\ZZ/24$. In particular, $24\eta$ is fiberwise 
homotopy trivial.

So, the classifying map $\cp2\to BO\to BPL$ for 24$\eta$ lifts to $F/PL$, In other words,  there exist a map $g: \CC \PP^2 \to F/PL$ such that the map  
$$ 
\CD 
\CC \PP^2 @>g>> F/PL @>>> BPL 
\endCD 
$$ 
classifies $24\eta$. 
Since $\la p_1(\eta) , [\CC \PP^2]\ra=1$, we have $\la(p_1(24\eta), [\CC \PP^2]\ra=24$, and therefore 
$\la L_1(24\eta),[\CC \PP^2]\ra=8$ (here $p_1$ and $L_1$ denote the first Pontryagin class and first Hirzebruch class, respectively), see \cite{MS}. Thus, $s(\CC \PP^2,g)=8/8=1$, and therefore $\langle K_4, g_*[\CC \PP^2] \rangle =1$.        
\qed

\m Let $h: \pi_4(F/PL)\to H_4(F/PL)$ be the Hurewicz homomorphism. Let $\tors$ denotes the torsion subgroup of $H_4(F/PL)$.

\begin{lemma}\label{l:sur}
 The map
\begin{equation*}
\begin{aligned}
\CD
a:\ZZ=\pi_4(F/PL) @>h>> H_4(F/PL)
@>\text{\rm quotient} >>  H_4(F/PL)/\tors=\ZZ
\endCD
\end{aligned}
\end{equation*}
is not surjective.
\end{lemma}

\p  Consider the Leray--Serre spectral sequence of the fibration \eqref{y-fib} and note that $H_4(Y)/\tors =\ZZ$, because $H_4(K(\ZZ,4))=\ZZ$ and all the groups $\wt H_i(K(\ZZ/2,2))$ are finite. Furthermore, $H_4(F/PL)\cong H_4(Y)$ since $Y$ is a Postnikov 4-stage of $F/PL$. Thus,  $H_4(F/PL)/\tors=\ZZ$.  

Because of \ref{c:gpl} and \ref{t:cl}, the subgroup $\la K^4, \IM a\ra $ of $\ZZ$ consist of even numbers. On the other hand, $\langle K_4, g_*[\CC \PP^2] \rangle =1$ by \lemref{cp2}. Thus, the image of  $g_*[\CC \PP^2]$ in $H_4(F/PL)/\tors$ does not belong to $\IM a$.
\qed

\m 
Consider the short exact sequence 
$$ 
\CD 
0 @>>> \ZZ @>2>>\ZZ[2] @>\rho >> \ZZ/2 @ >>> 0 
\endCD 
$$ 
where 2 over the arrow means multiplication by 2 and $\rho$ is the 
modulo 2 reduction. This exact sequence yields the Bockstein exact 
sequence 
\begin{equation} 
\begin{split} 
\cdots \longrightarrow H^n(X;\ZZ) & \xrightarrow {\ \ 2\ \ } H^n(X;\ZZ) 
\xrightarrow {\ \ \rho_*\ \ } H^n(X; \ZZ/2)\\ 
\label{bockstein} 
& \xrightarrow {\ \ \delta\ \ } H^{n+1}(X;\ZZ) \longrightarrow \cdots\,. 
\end{split} 
\end{equation} 
Put $X=K(\ZZ/2,n)$ and consider the fundamental class  
$$ 
\iota_n\in  H^n(K(\ZZ/2, n);\ZZ/2). 
$$  
Then we have the class $\gd :=\gd(\iota_n)\in  H^{n+1}(K(\ZZ/2,n),\ZZ)$.
According to what we said above, we  regard $\gd$ as a map $\gd: K(\ZZ/2,n)
\to K(\ZZ,n+1)$ and/or the cohomology operation 
\[
\gd: H^n(-;\ZZ/2)\to H^{n+1}(-,\ZZ).
\]

\begin{lemma}\label{l:N}
We have: $H^{n+3}(K(\ZZ/2,n))=\ZZ/2=\gd Sq^2\iota_n$ for all $n\geqslant 4$.
\end{lemma}

\p  Put $\iota=\iota_n$. We have 
\[
H^{n+3}(K(\ZZ/2,n), \ZZ/2)=\ZZ/2\oplus \ZZ/2=\{Sq^3\iota, Sq^2Sq^1\iota\},
\]~\cite[Ch. 9]{MT}. 
Let
\[
\gb:=\rho_*\gd: H^i(-;\ZZ/2)\to H^{i+1}(-;\ZZ/2)
\]
be the Bockstein homomorphism. Since $\gb(Sq^3\iota)\ne 0$, we conclude that
\[
Sq^3\iota\notin\IM\{\rho_*: H^{n+3}(K(\ZZ/2,n)) \to H^{n+3}(K(\ZZ/2,n), \ZZ/2)\},
\]
see~\cite[Ch. 11]{MT}. Since the homomorphism 
\[
\rho_*\otimes 1: H_*(-)\otimes \ZZ/2\to H_*(-;\ZZ/2)
\]
 is injective, we conclude that the group $H^{n+3}(K(\ZZ/2,n))$ is cyclic, and this group is 2-primary in view of Serre Class Theory,~\cite[Ch. 10]{MT}. Hence, the group  $H^{n+3}(K(\ZZ/2,n))$ generates $\gd Sq^2\iota$ since $\rho_* \gd Sq^2\iota=Sq^3\iota$ generates $H^{n+3}(K(\ZZ/2,n); \ZZ/2)$. Finally, $\gd Sq^2\iota$ has the order 2 since $2\gd=0$.
\qed

\m Recall that we denote by $\varkappa \in H^5(K(\ZZ/2,2))$ the characteristic class of the fibration \eqref{y-fib}.

\begin{thm}\label{postn} 
We have $\varkappa =\gd Sq^2\iota_2$, and it is an element of order $2$  in the group $H^5(K(\ZZ/2,2))=\ZZ/4$.
So,  $\varkappa$ is the first non-trivial Postnikov invariant of $F/PL$.
\end{thm}
 
\p Note that  $\varkappa \neq 0$ because of \ref{l:sur}. Indeed, otherwise $Y\cong K(\ZZ,4)\ts K(\ZZ/2,2)$. But this contradicts \lemref{l:sur}

Let $\Omega$ be the loop functor on category of topological spaces and maps. 
 Since $F/PL$ is an infinite loop space, the Postnikov invariant $\varkappa$ of $F/PL$ can be written as
$\Omega^N a_N$ for all $N$ and suitable 
\[
a_N \in [K(\ZZ/2,N+2), K(\ZZ. N+5)]= H^{N+5}(K(\ZZ/2, N+2); \ZZ).
\]
 By \lemref{l:N}, for $N>5$  the last group is equal to $\ZZ/2$. Thus, $\varkappa$ has the order 2. It is easy to see that 
 \[
  H^5(K(\ZZ/2,2))=\ZZ/4=\{x\}
  \]
   with $2x=\gd Sq^2\iota_2$, see e.g.~\cite[Lemma VI.2.7]{Rud}. 
Thus, $\varkappa =\gd Sq^2\iota_2$.  
\qed 
 
\begin{lemma}\label{l:free} 
Let $X$ be a finite $CW$-space such that the group $H_*(X)$ is torsion 
free. Let $Z$ be an infinite loop space such the groups $\pi_i(Z)$ have no odd torsion for all $i$. Then the group $[X,Z]$ is torsion free. In particular, the group $[X,F/PL[1/2]]$ is torsion free. 
\end{lemma} 
 
\p It suffices to prove  that $[X, Z[p]]$ is torsion free for 
every odd prime $p$. Note that $Z[p]$ is an infinite loop space since 
$Z$ is. So, there exists a connected $p$-local spectrum $E$ such that  $$ 
\wt E^0(-)=[-, Z[p]] =[-, Z\otimes \ZZ[p]. 
$$  
Moreover, $E^{-i}(\pt)=\pi_i(E)=\pi_i(Z)\otimes \ZZ[p]$,  
  So, because of the isomorphism $\wt E^0(X) 
\cong [X, Z[p]]$, it suffices to prove that $E^*(X)$ is torsion free. 
Consider the Atiyah--Hirzebruch spectral sequence for $E^*(X)$. Its 
initial term is torsion free because $E^*(\pt)$ and $H^*(X)$ are torsion free. Hence, the 
spectral sequence degenerates, and thus the group $E^*(X)$ is torsion free. 
\qed 
 
\begin{proposition}\label{p:null} 
Let $X$ be a finite $CW$-space such that the group $H_*(X)$ is torsion 
free. Let $f: X \to F/PL$ be a map such that $f^*K_{4i}=0$ and 
$f^*K_{4i+2}=0$ for all $i\geqslant 1$. Then $f$ is null-homotopic. 
\end{proposition} 
 
\p Consider the commutative square 
$$ 
\CD 
F/PL @>l_1>> F/PL[2]\\ 
@Vl_2VV @VVl_3V\\ 
F/PL[1/2] @>l_4>> F/PL[0] 
\endCD 
$$ 
where the horizontal maps are the $\mathbb Z[2]$-localizations and the 
vertical maps are the $\mathbb Z[1/2]$-localizations. Because of 
\ref{t:gpl},  $[X, F/PL]$ is a finitely generated abelian group, and so 
it suffices to prove that both $l_1\sirc f$ and $l_2\sirc f$ are 
null-homotopic. First, we remark that $l_2\sirc f$ is null-homotopic 
whenever $l_1\sirc f$ is. Indeed, since $H_*(X)$ is torsion free, the 
group $[X, F/PL[1/2]]$ is torsion free by \ref{l:free}. Now, if 
$l_1\sirc f$ is null-homotopic then $l_3\sirc l_1\sirc f$ is 
null-homotopic, and hence $l_4 \sirc l_2\sirc f$ is null-homotopic. 
Thus, $l_2\sirc f$ is null-homotopic since $[X, F/PL[1/2]]$ is torsion 
free. 
\par So, it remains to prove that $l_1\sirc f$ is null-homotopic.

Clearly, the equalities $f^*K_{4i}=0$ and $f^*K_{4i-2}=0$, $i>1$, 
imply that the map 
$$ 
\CD 
X @>>> F/PL @>l_1>> F/PL[2] \simeq Y[2] \times \Pi @>p_2>> \Pi 
\endCD 
$$ 
is null-homotopic. So, it remains to prove that the map 
$$ 
\CD 
g: X @>f>> F/PL @>l_1>> F/PL[2] \simeq Y[2] \times \Pi @>p_1>> Y[2] 
\endCD 
$$ 
is null-homotopic.

It is easy to see that $H^4(Y[2];\ZZ[2]) = \ZZ[2]$, see e.g.~\cite[VI.2.9(i)]{Rud}. Let $u \in
H^4(Y;\ZZ[2])$ be a free generator of the free
$\ZZ[2]$-module $H^4(Y;\ZZ[2])$. The fibration \eqref{y-fib} gives us the
following diagram with  the exact row: 
$$
\CD
H^4(X; \mathbb Z[2]) @>i_*>>
[X,Y[2]] @>p_*>> H^2(X;\mathbb  Z/2)\\ 
@. @VVu_*V @.\\ 
@. H^4(X;\mathbb Z[2]) @. 
\endCD 
$$ 

Note that  
$$ 
u_*i_* : \mathbb Z[2] \to \mathbb Z[2] 
$$ 
is the multiplication by $2\eps$ where $\eps$ is an invertible element of 
the ring $\mathbb Z[2]$, see e.g.~\cite[VI.2.9(ii)]{Rud}. Since $f^* K_2=0$, we conclude 
that $p_*(g)=0$, and so $g=i_*(a)$ for some $a\in H^4(X;\mathbb Z[2])$. Now, 
$$ 
0=u_*(g)=u_*(i_*a)=2a\eps. 
$$ 
 But $H^*(X; \mathbb Z[2])$ is torsion free, and thus $a=0$. 
\qed 

\m For completeness, we mention also that $F/PL[1/2] \simeq BO[1/2]$. So, 
there is a Cartesian square (see \cite{MM, Sul2})  
$$ 
\CD F/PL @>>> \Pi \times Y \\ 
@VVV @VVV \\ 
BO[1/2] @>\operatorname{ph}>> \prod K(\QQ, 4i) 
\endCD 
$$ 
where $\operatorname{ph}$ is the Pontryagin character.

 \section{Splitting  Theorems} 
\begin{definition}\label{def:split}\rm 
 Let $ A^{n+k}$ and $W^{n+k}$ be two  connected PL 
 manifolds (without boundaries), and let $M^n$ be a closed PL submanifold of 
$A$.  We say that a map $g: W^{n+k} \to  A^{n+k}$ {\it splits along $M^n$}  
if there exists a homotopy  
$$ 
g_t: W^{n+k} \to A^{n+k},\quad t\in I 
$$ 
such that: 
\par (i) $g_0=g$; 
\par (ii) there is a compact subset $K$ of $W$ such that 
$g_t|_{W\setminus K}=g|_{W\setminus K}$ for all $t\in 
I$; 
\par (iii) the map $g_1$ is transverse to  $M$ (and hence $g_1^{-1}(M)$ is a closed PL submanifold of $M_1)$, and the map 
$b:=g_1|_{g_1^{-1}(M)}: g_1^{-1}(M) \to M$ is a homotopy equivalence. 

We also say that the homotopy $G: W\ts I \to A, G(w,t)=g_t(w)$ {\it realizes the splitting of $g$}.
\end{definition} 
 
An important special case is when $A^{n+k}=M^n \times B^k$ for some 
connected manifold $B^k$. In this case we can regard $M$ as submanifold 
$M\times \{b_0\}, b_0\in B$ of $A$ and say that $g: W \to A$ splits along $M$ 
if it splits along $M\times \{b_0\}$. Clearly, this does not depend on the 
choice of $\{b_0\}$, i.e. $g$ splits along $M \times \{b_0\}$ if and only if 
it splits along $M\times\{b\}$ with any other $b\in B$.     
 
\m 
Recall that a map $f$ is called {\it proper} if $f^{-1}(C)$ is compact 
whenever $C$ is compact. A map $f: X \to Y$ is called a {\it proper 
homotopy equivalence} if there exist a map $g: Y \to X$ and
homotopies $F: gf \simeq 1_X$, $G: fg\simeq 1_Y$ such that all the four 
maps $f,g, F: X \times I \to X$, and $G: Y \times I \to Y$ are proper. 
 
\begin{theorem}\label{t:nov} 
 Let $M^n, n\geq 5$ be a closed connected $n$-dimensional 
 PL manifold such that $\pi_1(M)$ is a free abelian 
 group. Then every proper homotopy equivalence $h: W^{n+1} 
 \to M^n \times \mathbb R$ splits along $M^n$. 
\end{theorem} 
 
\p Because of the Thom transversality theorem, there is a 
homotopy $h_t: W \to M \times \RR$ 
 satisfies condition (ii) of \ref{def:split} and such that $h_1$ is transversal to $M$ . 
 We let $f=h_1$. Because of a crucial theorem of Novikov ~\cite[Theorem 3]{N2}, 
 there is a sequence of surgeries of the inclusion
  $f^{-1}(M) \subset W$ in $W$ such that the final result 
  of these surgeries is a homotopy equivalence $V \subset W$. Using the 
Pontryagin--Thom construction, we can realize this sequence of 
surgeries via a homotopy $f_t$ such that $f_t$ satisfies conditions 
(i)--(iii) of \ref{def:split} with $f_1^{-1}(M)=V$. 
\qed 
 
\begin{theorem}\label{t:fh} 
Let $M^n$ be a manifold as in $\ref{t:nov}$. Then every homotopy 
equivalence $W^{n+1} \to M^n \times S^1$ splits along $M^n$. 
\end{theorem} 
 
\p See~\cite{Fa}, cf. also \cite{FH}. 
\qed 
 
\begin{corollary}\label{c:fh} 
Let $M^n$ be a manifold as in $\ref{t:nov}$. Let $T^k$ denote the 
$k$-dimensional torus. Then every homotopy equivalence $W^{n+k} \to 
M^n \times T^k$ splits along $M^n$. 
\end{corollary} 
 
\p This follows from \ref{t:fh} by induction. 
\qed 
 
\begin{theorem}\label{t:splitting} 
Let $M^n$ be a manifold as in $\ref{t:nov}$. Then every homeomorphism 
$
h: W^{n+k} \to M^n \times \RR^k
$
splits along $M^n$. 
\end{theorem} 
 
\p We use the Novikov's torus trick. The inclusion $T^{k-1} \times 
\mathbb R \subset \mathbb R^k$ yields the inclusion 
$$
M \times T^{k-1} \times 
\mathbb R \subset M \times \RR^k.
$$ 
We set $W_1:= h^{-1}(M \times T^{k-1} \times R)$. Note that $W_1$ is 
a smooth manifold since it is an open subset of $W$. Now, set 
\begin{equation}\label{u-split}
u=h|_{W_1}: W_1 \to M \times T^{k-1} \times \mathbb R.
\end{equation}
 Then, by  \ref{t:nov}, $u$ splits along $M \times T^{k-1}$, i.e. there is a 
homotopy $u_t$ as in \ref{def:split}. We set $f:=u_1, V:=f^{-1}(M 
\times T^{k-1})$, and $g:=f|_V$. Because of \ref{c:fh}, $g: V \to M 
\times T^{k-1}$ splits along $M$. Hence, $f$ splits along $M$, and 
therefore $u$ in \eqref{u-split} splits along $M$. 
This splitting yields a homotopy $\overline u_t$ with $\overline u_0=u$ as in \defref{def:split}. Now, we define 
the homotopy 
$$
k_t: W \to M \times \mathbb R^k
$$
 by setting $k_t|_{W_1}:= 
\overline u_t|_{W_1}$ and $k_t|_{W\setminus W_1}:=h|_{W\setminus W_1}$. Note 
that $\{k_t\}$ is a well-defined and continuous family since the family $\{ 
\overline u_t\}$ satisfies \ref{def:split}(ii). It is clear that $k_t$ 
satisfies the conditions (i)---(iii) of \ref{def:split} and that $k_1$ 
extends $f$ on the whole $W$, i.e. $h$ splits along $M$. Thus, $h$ splits along $M$. 
\qed 
 
\begin{remarks}\rm 
1. The above used Theorems \ref{t:nov},\ref{t:fh}. and \ref{t:splitting} of Novikov and Farrell--Hsiang were 
originally proved for smooth manifolds, but they are valid for PL 
manifolds as well, because there is an analog of the Thom Transversality 
Theorem for PL manifolds, \cite{Wil}. 

2. In the above mentioned theorems we require the spaces to have free abelian fundamental groups. For arbitrary fundamental groups, there are obstructions to the splittings that involves algebraic $K$-theory of the fundamental group $\pi$. In fact, in \theoref{t:fh} there is an obstruction that is  in an element of the Whitehead group $\Wh(\pi)$ of  $\pi$. For \theoref{t:nov}, there are two obstructions: in $K^0(\pi)$ and in $\Wh(\pi)$. 
\end{remarks} 
 
\begin{lemma}\label{xi} 
Suppose that a map $g: W \to A$ splits along a submanifold $M$ of $A$. Let 
$\xi=\{E \to A\}$ be a PL bundle over $A$, let $g^*\xi =\{D \to W\}$, and let 
$\mathfrak I_g:g^*\xi \to \xi$ be the $g$-adjoint bundle morphism. Finally, let 
$l: D \to E$ be the map of the total spaces induced by $k$. Then $l$ splits 
over $M$. $($Here we regard $A$ as the zero section of $\xi$, and so $M$ 
turns out to be a submanifold of $E)$.    
\end{lemma} 

\p Let $G: W \times I \to A$ be a homotopy which realizes the splitting of $g$. 
Recall that $g^*\xi \times I$ is equivalent to $G^*\xi$. Now, the morphism  
$$ 
\CD 
g^*\xi \times I \cong G^*\xi @>\mathfrak I_g>> \xi 
\endCD 
$$  
gives us the homotopy $ D \times I \to E$ which realizes the
splitting of $l$. 
\qed  
 
\begin{lemma}\label{stab-split} 
Let $M$ be a manifold as in $\ref{t:nov}$. Consider two PL $\RR^N$-bundles 
$\xi=\{U \to M\}$ and $\eta=\{E\to M\}$ over $M$ and a topological morphism 
$\gf: \xi \to \eta$ over $M$ of the form  
$$  
\CD 
U @>g>> E\\ 
@VVV @VVV\\ 
M @= M . 
\endCD 
$$ 
Then there exists $k$ such that the map 
$$ 
g\times 1: U \times \RR^k \to E \times \RR^k 
$$  
splits along $M$, where $M$ is regarded as the zero section of $\eta$.   
\end{lemma} 
 
\p Take a PL $\RR^m$-bundle $\zeta$ such that $\eta \oplus \zeta=\RR60 ^{N+m}$ 
and let $W$ be the total space of $\xi \oplus \zeta$. Then the morphism   
$$ 
\gf\oplus 1: \xi \oplus \zeta \to \eta \oplus \zeta = \theta^{N+m} 
$$ 
yields a map   of the total spaces
\begin{equation}\label{WM} 
\Phi: W \to M \times \RR^{N+m}. 
\end{equation} 
Because of Theorem \ref{t:splitting}, the map $\Phi$ 
splits along $M$. Furthermore, the morphism   
$$ 
\gf \oplus 1 \oplus 1: \xi \oplus \zeta \oplus \eta \to  \eta \oplus \zeta \oplus \eta 
$$ 
yields a map  of the total spaces
$$ 
g\times 1: U \times \RR^{2N+m} \to E \times \RR^{2N+m}. 
$$ 
In view of Lemma \ref{xi}, this map splits over $M$ because $\Phi$ does. So, 
we can put $k=2N+m$.  
\qed 
 
\m Now, let $a: TOP/PL \to F/PL$ be a map as in \eqref{a-b}. 
 
\begin{theorem}\label{t:vanishing} 
Let $M$ be as in $\ref{t:nov}$. Then the composition 
$$ 
\CD 
[M,TOP/PL] @>a_*>> [M, F/PL] @>s>>P_{\dim M} 
\endCD 
$$ 
is trivial, i.e., $sa_*(v)=0$ for every $v\in [M, TOP/PL]$. In other 
words, $s(M, af)=0$ for every $f: M \to TOP/PL$. 
\end{theorem} 
 
\p In view of \eqref{eq:tpl-lim}, every element $v\in [M, TOP/PL]$ gives us 
a (class of a) topological morphism  
$$ 
\gf: \nu^N_M \longrightarrow \gga^N_{PL} 
$$  
of PL $\RR^N$-bundles. To map the class $v\in [M,TOP/PL]$ to the class 
$a_*v\in [M,F/PL]$, we must convert $\gf$ to the (equivalence class of the) $(S^N,*)$-morphism 
$\gf\bul: (\nu_M)\bul\to (\gga^N_{PL})\bul$. Now, we follow \ref{construdef:ass} and construct a commutative diagram 
\begin{equation} 
\CD 
U\bul @>g>>  {U'}\bul\\ 
@VqVV @VVpV\\ 
M @= M 
\endCD 
\end{equation} 
like \eqref{uu}. However, here $g$ is a homeomorphism. Thus, 
$g(U)=U'$, and so we get the diagram  
\begin{equation}\label{u-pl} 
\CD 
U @>g>> {U'}\\ 
@VqVV @VVpV\\ 
M @= M 
\endCD 
\end{equation} 
which is a topological morphism of PL bundles over $M$. 

We embed $M$ in $U'$ as the zero section of $p$. By the definition of the map $s$, we conclude that $s(M, 
a_*v)=0$ if the map $g: U \to U'$ splits along $M$ (because in this case the 
associated normal morphism is a map over a homotopy equivalence). Furthermore 
for any $k$,  the topological morphisms $\gf$ and  
$$  
\CD  
\nu_M \oplus \theta^k @>(\gf\oplus 1)>> \gga^N_{PL}\oplus \theta^k @>>> 
\gga^{N+k}_{PL}  \endCD  
$$ 
represent the same element of $[M, TOP/PL]$. Hence,  $s(M, a_*v)=0$ provided there 
exists {\it at least one} $k$ such that the map   
$$ 
g\times 1: U \times  \RR^k \to U' \times \RR\sp k 
$$ 
splits along $M$. But this follows from Lemma \ref{stab-split}, since \eqref{u-pl} is a {\it topological} morphism
\qed

\m Now we show that the condition $\dim M \geq 5$ in \ref{t:vanishing} 
is not necessary. 

\begin{corollary}\label{c:dim} 
Let $M$ be a closed connected PL manifold such that $\pi_1(M)$ is 
a free abelian group. Then $s(M, af)=0$ for every map $f: M \to TOP/PL$. 
\end{corollary} 
 
\p Let $\cp2$ denote the complex projective plane, and let 
$$ 
p_1: M \times \cp2 \to M 
$$  
be the projection on the first factor. Then 
$s(M \times \cp2, gp_1)=s(M,g)$ for every $g: M \to F/PL$, see 
\cite[Ch. III, \S 5]{Br}. In particular, for every map $f; M \to TOP/PL$ 
we have  
$$ 
s(M, af)= s(M \times \CC \PP^2, (af)p_1)=s(M \times \CC \PP^2, 
a(fp_1))=0 
$$ 
where the last equality follows from Theorem \ref{t:vanishing}.  
\qed

\section{Detecting Families} 
 
\m Recall the terminology: a singular smooth manifold in a space $X$ is a map 
$M \to X$ of a smooth manifold $M$. 
 
\m Given a $CW$-space $X$, consider a connected closed smooth singular 
manifold $\gamma: M \to X$ in $X$. Then, for every map $f: X \to 
F/PL$, the invariant $s(M, f \gamma)\in P_{\dim M}$ is defined. 
Clearly, if $f$ is null-homotopic then $s(M, f \gamma)=0$. 
 
\begin{definition}\rm 
 Let $\{\gamma_j: M_j \to X\}_{j\in J}$ be a family of closed connected  
smooth singular manifolds in $X$; here $J$ is an index set. We say that the family $\{\gamma_j: M_j \to 
X\}$ is a {\it detecting family} for $X$ if, for every map $f: X \to F/PL$, 
the validity of  all the equalities $s(M_j, f \gamma_j)=0, j\in J$ implies 
that $f$ is null-homotopic. 
\end{definition} 
 
\m Note that $F/PL$ is an $H$-space, and hence, for every detecting family  
$\{\gamma_j: M_j\to X\}$, the collection $\{s(M_j, f \gamma_j)\}$ determine a map 
$f: X \to F/PL$ uniquely up to homotopy.  
 
\m The concept of detecting family is related to Sullivan's ``characteristic 
variety'', but it is not precisely the same. If a family $\mathcal F$ of 
singular manifolds in $X$ contains a detecting family, then $\mathcal F$ 
on its own is a detecting family. On the contrary, the characteristic variety is 
in a sense ``minimal'' detecting family. 
 
\begin{lemma}\label{l:detecting} 
Let $X$ be a finite $CW$-space such that the group $H_*(X)$ 
is torsion free. Let $\{\gamma_j: M_j \to X\}$ be a family of smooth oriented closed 
connected singular manifolds in $X$ such that, for each $m$, the elements 
$(\gamma_j)_*[M^{2m}_j]$ generate the group $H_{2m}(X)$. Then $\{\gamma_j\}$ is a detecting family for $X$
\end{lemma} 
 
\p Consider a map $f: X \to F/PL$ such that $s_j(M_j, f 
\gamma_j)=0$ for all $j\in J$. We must prove that $f$ is null-homotopic. 

Because of \ref{p:null}, it suffices to prove that 
$f^*K_{i}=0$ and $f^*K_{4i-2}=0$ for all $i\geq 1$. Furthermore, 
$H^*(X)=\Hom (H_*(X), \mathbb Z)$ because $H_*(X)$ is torsion free. 
So, it suffices to prove that 
\begin{equation}\label{4i} 
\langle f^*K_{4i}, z \rangle =0 \text{ for all } z\in H_{4i}(X) 
\end{equation} 
and 
\begin{equation}\label{4i+2} 
\langle f^*K_{4i-2}, z\rangle =0 \text{ for all } z\in 
H_{4i-2}(X;\mathbb Z/2). 
\end{equation} 
First, we prove \eqref{4i}. Since the classes $(\gamma_j)_*[M_j], \dim 
M_j=4i$ generates the group $H_{4i}(X)$, it suffices to prove that 
$$ 
\left\langle f^*K_{4i},\, (\gamma_j)_*[M_j] \right\rangle =0 
\text{ whenever $\dim M_j= 4i$}. 
$$ 
But, because of \ref{t:cl}, for every $4i$-dimensional $M_i$ we 
have 
$$ 
0=s(M_j, f \gamma_j)= 
\left\langle (f \gamma_i)^*K_{4i},\, [M_j] \right\rangle = 
\left\langle f^*K_{4i},\, (\gamma_j)_*[M_j] \right\rangle. 
$$ 
This completes the proof of the equality \eqref{4i}. 

For the case $i=4k-2$, note that the group 
$H_{4i-2}(X; \mathbb Z/2)$ is generated by the elements 
$(\gamma_j)_*[M_j]_2, \dim M_j=4k-2$, since $H_*(X)$ is torsion free. 
Now the proof can be completed similarly to the case $i=4k$. 
\qed

 \begin{thm}\label{t:detect} 
Let $X$ be a connected finite $CW$-space such that the group $H_*(X)$ 
is torsion free. Then $X$ admits a detecting family
$\{\gamma_j: M_j \to X\}$ such that each $M_j$ is orientable. 
\end{thm} 
 
\p Since $H_*(X)$ is torsion free, every homology class in $H_*(X)$ 
can be realized by a closed connected smooth oriented singular 
manifold, see e.g. \cite[15.2]{Conner} or \cite[6.6 and 7.32]{Rud}. Now apply \lemref{l:detecting}.
\qed
 
 \begin{example}\label{tksn}\rm Let $X$ be the space $T^k\ts S^n$.  Clearly, $H_{2m}(X)$ is generated by fundamental classes of submanifolds $T^{2m}$ and $T^{2m-n}\ts S^n$ of $T^k\ts S^n$. Hence, $X$ has a detected family $\{\gamma_j: M_j \to X\}$ such that each $M_j$ is  either $T^r$ or $T^r\times S^n$. 
 \end{example}
 
\section[A Special Case]{Normal Invariant of a 
Homeomorphism: a Special Case} 
 
\begin{theorem}\label{t:weaknorm} 
 If the element $x\in \spl(T^k\ts S^n)$ can be represented by a {\em homeomorphism} $h: V \to 
M$, then $\jf(x)=0$
\end{theorem}

\p Put $M=T^k\ts S^n$. The maps $\jtop$ and $\jf$ from section~\ref{structures} can be included in the 
 commutative diagram
\begin{equation}\label{corr3} 
\CD 
\tpl(M) @>\jtop>> [M, TOP/PL]\\ 
@VVV @VVa_*V\\ 
\spl(M) @>\jf>> [M, F/PL] 
\endCD 
\end{equation} 
where the left arrow is the obvious forgetful map and $a_*$ is induced 
by $a$ as in \eqref{a-b}. 
 
Suppose that $x$ can be represented by a homeomorphism $h: V \to M$. Consider 
a map $f: M \to TOP/PL$ such that $\jtop (h)$ is 
homotopy class of $f$. Then, clearly, the class $\jf(x)\in [M,F/PL]$ is 
represented by the map 
$$ 
\CD 
M @>f>> TOP/PL @>a>> F/PL. 
\endCD 
$$

As we explained in Example~\ref{tksn}, $M$ possesses a detecting family 
$\{\gamma_j: M_j \to M\}$ where each $M_j$ is either $T^r$ or $T^r\ts S^n$. Hence, by \ref{t:vanishing} 
and \ref{c:dim}, $s(M_j, af \gamma_j)=0$ for all $j$. So, $af$ is null-homotopic since $\{\gamma_j\}$ is a detecting family . Thus, 
$\jf(x)=0$.  
\qed 
 \newpage

\vfill
 \pagebreak
 
 	 \chapter{{Applications and Consequences}\label{c:norminv}}

\section{The Space $F/TOP$} 
 
Because of the Main Theorem and results of Freedman~\cite{F} and Scharlemann~\cite{Schar}, the Transversality Theorem holds for {\it topological} manifolds and bundles. For the references, see Rudyak~\cite[IV.7.18]{Rud}.   
 
\m 
Since we have the topological transversality, we can define the maps 
$$ 
s': [M, F/TOP] \to P_{\dim M} 
$$  
where $M$ turns out to be a topological manifold. These map $s'$ are obvious analog of maps $s$ defined in \eqref{sul}: you need merely replace PL by TOP in Equation~\eqref{sul} and \theoref{t:cl}. We leave it to the reader.

\m The following proposition states the main difference between $F/PL$ and $F/TOP$. 
 
\begin{prop}\label{dim4} 
The map $s': \pi_4(F/TOP) \to \ZZ$ 
is a surjection. 
\end{prop}  
 
\p Note that the Freedman manifold $V$ from \theoref{t:F} is almost 
parallelizable and has signature 8. Now the proof can be completed just as 
\ref{image}(i). \qed 

\begin{remark}\rm Kirby and Siebenmann~\cite{KS2} used a {\it homology} 
4-manifold of signature 8 in order to prove \propref{dim4}. The paper of Freedman appeared later.
\end{remark}

\begin{thm}\label{ftop} 
{\rm (i)} For $i \ne 4$ the map $b: F/PL \to F/TOP$ induces an isomorphism 
$$ 
b_*: \pi_i(F/PL) \to \pi_i(F/TOP). 
$$ 
{\rm (ii)} The homomorphism 
$$ 
b_*: \ZZ = \pi_4(F/PL) \to  \pi_4(F/TOP) =\ZZ
$$ 
is the multiplication by $2$. 
\end{thm} 

\p 
(i) Recall that $TOP/PL=K(\ZZ/2,3)$ and $\pi_4(F/PL)=\ZZ$. So, the exactness of the homotopy sequence of the  
fibration 
\begin{equation*}  
\CD 
TOP/PL @>a>> F/PL @>b>> F/TOP 
\endCD 
\end{equation*} 
in \eqref{a-b} yields an isomorphism $b_*: \pi_i(F/PL) \cong \pi_i(F/TOP)$ for $i\ne 4$.  

(ii) We have the commutative diagram 
$$ 
\CD 
0 =\pi_4(TOP/PL) \\ 
@Va_*VV @.\\ 
\ZZ = \pi_4(F/PL) @>s>> \ZZ \\ 
 @V b_*VV  @. \\ 
 \pi_4(F/TOP) @>s'>> \ZZ \\ 
@VVV @. \\   
\ZZ/2= \pi_3(TOP/PL) \\ 
 @VVV @. \\ 
0= \pi_3(F/PL) \\ 
\endCD 
$$ 
 where the middle vertical line is a short exact sequence. Therefore $\pi_4(F/TOP)=\ZZ$ or $\pi_4(FTOP)=\ZZ \oplus\ZZ/2$. By \theoref{image}(iii),
$\IM s$ is the subgroup $2\ZZ$ of $\ZZ$, while $s'$ is a surjection by 
\ref{ftop}. Thus,  $\pi_4(F/TOP)=\ZZ$ and $b_*$ is the multiplication by 2.  
\qed 
 
\m Now, following \ref{t:cl}, we can introduce the classes  
$$ 
K'_{4i}\in H^{4i}(F/TOP,\ZZ[2]) \text{ and } K'_{4i-2}\in 
H^{4i-2}(F/TOP,\ZZ/2) 
$$  
such that 
$$ 
s'(M^{4i},f) = \langle f^*K'_{4i},\, [M] \rangle \text{ and } 
s'(N^{4i-2},f) = \langle f^*K'_{4i-2},\, [N]_2 \rangle. 
$$ 
However, here $M$ and $N$ are allowed to be topological (i.e. not 
necessarily PL) manifolds. 
 
\m Similarly to \eqref{bigpi}, together these classes yield the map 
$$ 
K': F/TOP \longrightarrow \prod_{i>0}\left(K(\ZZ[2], 4i) \times 
K(\ZZ/2,4i-2)\right). 
$$ 
 such that for each $i>0$ the map
\begin{equation*}
 \CD
F/TOP @>K'>> \Pi @>\text{ projection}>>  K(\ZZ/2,4i-2) \text{ (resp. } K(\ZZ[2], 4i)
 \endCD
 \end{equation*}
 coincides with $K'_{4i-2}$ (resp. $K'_{4i}$).
 
\begin{theorem} 
The map  
$$ 
K'[2]: K': F/TOP[2] \lr \prod_{i>0}\left(K(\ZZ[2], 4i) \times 
K(\ZZ/2,4i-2)\right) 
$$ 
is a homotopy equivalence. 
\end{theorem} 
 
\p Together \ref{image} and \ref{dim4} imply that the homomorphisms 
$$
s':  \pi_{2i}(F/TOP) \lr P_{2i}
$$ 
are surjective. Now, in view of \ref{ftop}, all 
the homomorphisms $s'$'s are isomorphisms, and the result follows.  
\qed 
 
\m So, the only difference between the spaces $F/PL$ and 
$F/TOP$ is that $F/TOP[2]$ has  trivial Postnikov invariants, while $F/PL[2]$ 
has exactly one non-trivial Postnikov  invariant $\gd Sq^2\iota_2 \in H^5(K(\ZZ /2,2) 
;\ZZ[2])$.

\m Now we discuss the groups $\pi_i(BTOP)$. Consider the map 
\[
\ga=\ga^{PL}_{TOP} : BPL\to BTOP
\] 
and the fibration 
\[
TOP/PL \lr BPL \stackrel{\ga}\lr BTOP
\]
as in \eqref{albet}.  Since $\pi_3(BPL)=0$ and $\pi_i(TOP/PL)=0$ for $i\neq 3$, we conclude that 
\[
\ga_*: \pi_i(BPL)\lr \pi_i(BTOP)
\]
is an isomorphism for $i\ne 4$. Furthermore, we have the exact sequence
\[
0\lr \pi_4(BPL)\stackrel{\ga_*}\lr \pi_4(BTOP)\lr \pi_3(TOP/PL)\lr
\]
where $\pi_4(BPL)=\ZZ$ and $\pi_3(TOP/PL)=\ZZ/2$. Hence we have that either $\pi_4(BTOP=\ZZ$ or $\pi_4(BTOP)=\ZZ\oplus \ZZ/2$.

Now, consider the diagram of fibrations
$$
\CD
F/PL @ >>> BPL @>>> BF\\
@VVV @VVV @|\\
F/TOP @ >>> BTOP @>>> BF
\endCD
$$

It is known that $J$-homomorphism 
$$
J: \ZZ=\pi_4(BPL)\lr \pi_4(BSF)=\ZZ/24
$$ 
is surjective~\cite{Adams, MK} (recall that $\pi_i(PL/O)=0$ for $i<7$, and so there is no difference between $\pi_i(BPL)$ and $\pi_i(BO)$ up to dimension 6). Furthermore, $\pi_5(BF)$ is finite and $\pi_3(F/PL)=\pi_3(F/TOP)=0$. Now, we apply $\pi_4$ and get the commutative diagram with exact rows 
$$
\CD
0@>>>\ZZ @>24 >> \ZZ @>>>\ZZ/24@>>>0\\
@. @V2VV @VVV @VV=V @.\\
0@>>>\ZZ @>>> \pi_4(BTOP) @ >>> \ZZ/24 @>>>0.
\endCD
$$
The assertion $\pi_4(BTOP)$ contradicts the commutativity of the diagram. Thus, $\pi_4(BTOP)=\ZZ\oplus\ZZ/2$. Cf. Milgram~\cite{Mil}.

\section{The Map $a: TOP/PL \to F/PL$}

 \m Recall that in \eqref{a-b} we described the fibration 
\begin{equation*}  
\CD 
TOP/PL @>a>> F/PL @>b>> F/TOP. 
\endCD 
\end{equation*} 

\begin{prop}\label{a-ess} 
The map $a: TOP/PL \to F/PL$ is essential. 
\end{prop} 
 
\p  For general reasons, the fibration 
$$ 
\CD 
TOP/PL @>a>> F/PL @>>> F/TOP 
\endCD 
$$ 
yields a fibration 
$$ 
\CD 
\Omega(F/TOP) @>u>> TOP/PL @>a>> F/PL. 
\endCD 
$$ 
 
If $a$ is inessential then there exists a map 
\[
v: TOP/PL \to \Omega(F/TOP)
\]
with $uv\simeq 1$. But this is impossible because $\pi_3(TOP/PL)=\ZZ/2$ while 
$\pi_3(\Omega(F/TOP))=\pi_4(F/TOP)=\ZZ$. 
\qed
 
\m  Let $\ell: F/PL \to F/PL[2]$  denote the localization map. Let  $\psi: F/PL \to Y$ be the Postnikov 4-approximation of $F/PL$ as in \eqref{psi}. Take an arbitrary map $f: X \to TOP/PL$.
 
\begin{prop}\label{ess} 
The following three conditions are equivalent: 
 
\par{\rm (i)} the map  
$$ 
\CD 
X @>f>> TOP/PL  @>a>> F/PL
\endCD 
$$  
is essential; 
 
\par{\rm (ii)} the map  
\begin{equation*} 
\CD 
X @>f>> TOP/PL  @>a>> F/PL @>\ell>> F/PL[2] 
\endCD 
\end{equation*} 
is essential; 
 
\par{\rm (iii)} the map  
$$ 
\CD 
X @>f>> TOP/PL @>a>> F/PL @>\ell>> F/PL[2] @>\psi[2]>> Y[2] 
\endCD 
$$ 
is essential. 
 
\end{prop} 
 
\p It suffices to prove that (i) $\Rightarrow $ (ii) $\Rightarrow $(iii). To 
prove the first implication, recall that a map $u: X \to F/PL$ is inessential 
if both localized maps  
$$  
\CD 
X@>u>> F/PL \to F/PL[2], \qquad  X@>u>> F/PL \to F/PL[1/2] 
\endCD 
$$ 
are inessential. Now,  (i) $\Rightarrow $ (ii) holds since $TOP/PL[1/2]$ is 
contractible. 
 
To prove the second implication, note that a map $v: X \to F/PL[2]$ is 
inessential if both maps  (we use notation as in \ref{t:gpl})  
$$ 
\CD 
X @>v>>  F/PL[2] @> K >> \Pi, \qquad X @>v>>  F/PL[2] @>>>  Y 
\endCD 
$$ 
are inessential. So, it suffices to prove that the map  
$$ 
\CD 
X @>\ell a f>>  F/PL[2] @>>> \Pi 
\endCD 
$$ 
is inessential. This holds, in turn, because the map $TOP/PL \to F/PL \to 
F/TOP$ is inessential and the diagram   
$$ 
\CD 
F/PL[2] @ > K[2] >> \Pi @= \Pi\\ 
 @VV b[2] V @.  @| \\ 
F/TOP[2] @>K'[2]>> \prod_{i>0} (K(\ZZ/2, 4i-2) \times K(\ZZ[2], 4i)) 
@>\text{proj}>> \Pi  
\endCD 
$$ 
commutes. 
\qed 
 
\m Consider the fibration  
$$ 
\CD 
K(\ZZ[2], 4) @>i>> Y [2]@ >>> K(\ZZ/2. 2) 
\endCD 
$$ 
that is the $\ZZ[2]$-localization of the fibration \eqref{y-fib}. 
 
\begin{lemma}\label{x-y} 
For every space $X$, the homomorphism 
$$ 
\CD 
H^4(X; \ZZ[2])=[X, K(\ZZ[2], 4)] @>i_* >> [X,Y[2]] 
\endCD 
$$ 
is injective. Moreover, $i_*$ is an isomorphism if $H^2(X;\ZZ/2)=0$. 
\end{lemma} 
 
\p The fibration \eqref{y-fib} yields the exact sequence (see e.g. \cite{MT}) 
\begin{equation}\label{exact-y} 
H^1(X;\ZZ/2) \xrightarrow{\gd Sq^2} H^4(X;\ZZ[2]) \xrightarrow{i_*} [X,Y[2]] 
\rightarrow H^2(X;\ZZ/2) \end{equation} 
where $\gd Sq^2(x)\equiv 0$ (because $\gd Sq^2(x)=0$ whenever $\deg x=1$). 
\qed 
 
\m Let $g: TOP/PL \to Y$ be the composition 
$$ 
\CD 
TOP/PL @>a>> F/PL @>\ell >> F/PL[2] @>\psi[2]>> Y[2]. 
\endCD 
$$ 
 
Note that $g$ is essential because of \ref{a-ess} and \ref{ess}. 
 
\begin{cory}\label{g=igd} 
The map 
$$ 
\CD 
TOP/PL =K(\ZZ/2,3) @>\gd >> K(\ZZ[2],4) @>i >> Y[2] 
\endCD 
$$ 
is homotopic to $g$, i.e. $g\simeq i\gd$. 
\end{cory} 
 
\p Because of \lemref{x-y} applied to $X=K(\ZZ/2,3$, the set $[K(\ZZ/2,3),Y[2]]$ has 
exactly two elements. Since both maps $g$ and $i\sirc \gd$ are essential (the 
last one because of Lemma \ref{x-y}), we conclude that $g\simeq i\gd$.  
\qed 
 
\begin{theorem}\label{gd-ess} 
Given a map $f: X \to TOP/PL$, the map  
\begin{equation*} 
\CD 
X @>f>> TOP/PL @>a>> F/PL 
\endCD 
\end{equation*} 
is essential if and only if the map 
\begin{equation*} 
\CD 
X @>f>> TOP/PL = K(\ZZ/2,3) @>\gd>> K(\ZZ[2],4) 
\endCD 
\end{equation*} 
is essential. 
\end{theorem} 
  
\p We have the chain of equivalences
\begin{equation*}
\begin{aligned}
af \text{ is essential }&\overset{\ref{ess}} {\Longleftarrow\! \Longrightarrow} gf
\text{ is essential }\stackrel{\ref{g=igd}}{\Longleftarrow\! \Longrightarrow}  i\gd \text{ is essential }\\
 &\stackrel{\ref{x-y}}
{\Longleftarrow\! \Longrightarrow}  \gd f\text{ is essential. }
\end{aligned}
\end{equation*}
\qed

\section{Normal Invariant of a Homeomorphism}\label{secnorm} 
 
\begin{lemma}\label{2-tor}  
Let $X$ be a finite $CW$-space such that $H_n(X)$ 
is $2$-torsion free. Then the homomorphism   
$$ 
\gd: H^n(X;\ZZ/2) \to H^{n+1}(X;\ZZ[2]) 
$$  
is zero. 
\end{lemma} 
 
\p Because of the exactness of the sequence \eqref{bockstein} 
$$ 
\CD 
H^n(X;\ZZ/2) @>\gd >> H^{n+1}(X;\ZZ[2]) @> 2 >> H^{n+1}(X;\ZZ[2]), 
\endCD 
$$  
it suffices to prove that $H^{n+1}(X;\ZZ[2])$ is 2-torsion free. Since 
$H_n(X)$ is 2-torsion free, we conclude that $\Ext (H_n(X),\ZZ[2])=0$. (Indeed, $\Ext(\ZZ/m, A)=A/mA$ for all $A$.) Thus,  because of the Universal Coefficient Theorem,
\begin{eqnarray*} 
H^{n+1}(X;\ZZ[2]) &=& \Hom (H_{n+1}(X;\ZZ[2]) \oplus \Ext (H_n(X);\ZZ[2])\\ 
&=& \Hom (H_{n+1}(X;\ZZ[2]), 
\end{eqnarray*} 
and the result follows. 
\qed 
 
\begin{theorem}\label{ni} 
Let $M$ be a closed PL manifold such that $H_3(M)$ is $2$-torsion free. Then 
the normal invariant of any homeomorphism $h: V \to M$ is trivial.  
\end{theorem}  
 
\p Since $h$ is a homeomorphism, the normal invariant $\jf(h)$ turns out to 
be the homotopy class of a map  
$$ 
\CD 
M @>f>> TOP/PL @>a>> F/PL 
\endCD 
$$ 
where the homotopy class of $f$ is $\jtop(h)$. Because of \ref{ess} and 
\ref{x-y}, it suffices to prove that the map  
$$ 
\CD 
M @>f>> TOP/PL =K(\ZZ/2,3) @>\gd >> K(\ZZ[2],4) 
\endCD 
$$ 
is inessential. But this follows from Lemma~\ref{ni}. 
\qed 
 
 \m Now we have the following version of the Hauptvermutung, cf.~\cite[Corollary on p.68]{Cas} and~\cite[Theorem H on p. 93]{Sul2}.
 
\begin{cory} Let $M, \dim M \geqslant 5$ be a closed simply-connected PL manifold 
such that $H_3(M)$ is $2$-torsion free. Then every homeomorphism $h: V \to M$ 
is homotopic to a PL homeomorphism. In particular, $V$ and $M$ are PL homeomorphic. 
\end{cory} 
 
\p This follows from \ref{bns} and \ref{ni}. 
\qed 
 
\begin{remark}\rm 
Rourke \cite{Rou} suggested another proof of \ref{ni}, using the technique of 
simplicial sets. \end{remark}

 \section{Kirby-Siebenmann and Casson-Sullivan Invariants}\label{ksks}
 
Recall some facts on obstruction theory~\cite{DK, FFG, H, MT, Spa2}. Let $F\to E \to B$ be a principal $F$-fibration such that $F$ is an Eilenberg-MacLane space $K(\pi,n)$, and assume that the $\pi_1(B)$-action on $\pi=\pi_n(F)$ is trivial. Let $\iota=\iota_n\in H^n(K(\pi,n);\pi)$ be the fundamental class of $F$. and let $\gk=\tau\iota\in H^{n+1}(B;\pi)$ be the characteristic class of the fibration $F\to E \to B$, where $\tau: H^n(F;\pi)\to H^{n+1}(B;\pi)$ is the transgression. This is well-known that the fibration $F\to E \to B$ admits a section if and only if $\gk=0$ and, if a section exists then the vertically homotopy class of sections of the fibration are in a bijective correspondence with elements of $H^n(B;\pi)$. Hence, given a map $f:X\to B$, the map $f$ can be lifted to $E$ iff $f^*(\gk)=0$, and the vertical homotopy classes of liftings of $f$ to $E$ are in a bijective correspondence with elements $H^n(X;\pi)$ provided such a lifting exists.
 
 \m Since $TOP/PL$ is the Eilenberg-MacLane space $K(\ZZ/2,3)$, we can apply previous arguments to the principal $TOP/PL$-fibration \eqref{albet}
 $$
 \CD
 TOP/PL @>>> BPL @>a^{PL}_{TOP}>> BTOP.
 \endCD
 $$
Then we get the characteristic class
\begin{equation}\label{eq:univksclass}
\gkk=\tau\iota\in H^4(BTOP;\pi_3(K(\ZZ/2,3))=H^4(BTOP;\ZZ/2)
\end{equation}
where $\iota\in H^3(K(\ZZ/2,3);\ZZ/2)$ is the fundamental class. We call $\gkk$ the {\it universal Kirby-Siebenmann class}.

\m Let $M$ be a topological manifold, and let $f: M\to BTOP$ classify the stable tangent bundle of $M$. Since $f$ is unique up to homotopy, the class $f^*(\gkk)\in H^*(M;\ZZ/2)$ is a well-defined invariant of $M$. We put
\begin{equation}\label{eq:ksclass}
\gkk(M):=f^*(\gkk)\in H^*(M;\ZZ/2)
\end{equation}
and call it the {\it Kirby-Siebenmann class of $M$}.

\begin{thm}\label{t:ks-exist} Let $M$ be a topological manifold. If $M$ admits a PL structure then $\gkk(M)=0$.  If $\dim M \geqslant 5$ and $\gkk(M)=0$ then $M$ admits a PL structure. In particular, if $\dim M \geqslant 5$ and $H^4(M;\ZZ/2)=0$ then $M$ admits a PL 
structure. 
\end{thm}

\p If $M$ admits a PL structure then the classifying map $f: M \to BTOP$ can be lifted to $BPL$, and hence $f^*(\gkk)=0$, i.e $\gkk(M)=0$. Conversely, if $\gkk(M)=0$ then $f$ can be lifted to $BPL$. Thus, in case $\dim M \geqslant 5$ the manifold $M$ admits a PL structure by \corref{exist}.
\qed         

\begin{thm}\label{ks-class}
If a topological manifold $M, \dim M \geqslant 5$ admits a PL structure then set of concordance classes of PL structure on $M$ is in bijective correspondence with $H^3(M;\ZZ/2$, i.e.
\[
\tpl(M)\cong H^3(M;\ZZ/2.
\]
In paticular, if $H^3(M;\ZZ/2)=0$ then the Hauptvermutung holds for $M$.
\end{thm}
 
 \p Because of \corref{classif}, we have a bijection
 \[
 \tpl(M)\cong [M, TOP/PL].
 \]
 Thus, because of the Main Theorem $TOP/PL \simeq K(\ZZ/2,3$,  we get
 \[
 \tpl(M)\cong [M, TOP/PL]\cong [M, K(\ZZ/2,3]\cong H^3(M;\ZZ/2). 
 \]
 \qed
 
\begin{definition}\rm Let $M$ be a PL manifold and $h: V\to M$ be a PL structuralization. In view of bijection from \theoref{ks-class}, the PL structure $h$ gives us a  cohomology class $\gkk(h)\in H^3(M;\ZZ/2)$. This class is called the {\it Casson-Sullivan invariant} of $h$, and it measures the difference between $h:V \to M$ and $1_M$. 
\end{definition}  
 
\m So, $\gkk (h)=0$ if and only if $h: V\to M$ is concordant to the identity map of $M$. It is also worthy to mention that, for every $a\in H^3(M;\ZZ/2)$ there exists a homeomorphism $h: V \to M$ with $a=\gkk(h)$. 
 
 \begin{remark} \rm We know that Hauptvermutung holds for $T^k\times S^n$ with $k+n\geqslant 5$ and $n\geqslant 3$,~\cite{HS}. In other words, if two PL manifolds $M_1, M_2$ are homomorphic to $T^k\times S^n$ then there are PL homeomorphic. On the other hand, the group $H^3(T^k\times S^n;\ZZ/2)$ is quite large for $k$ large enough, i.e. $T^k\times S^n$ has many different PL structure. Is it a contradiction? No, it is not. The explanation comes because, given a homeomorphism $h:T^k\times S^n\to T^k\times S^n$, there are many PL concordance classes $T^k\times S^n\to T^k\times S^n$ that are homotopic to $h$.
\end{remark}

\section{Several Examples}\label{s:ex} 
 
\begin{example}\label{hptv} 
There are two closed PL manifolds that are homeomorphic but not PL homeomorphic. 
\end{example} 
   
\m Let $\RR\PP^n$ denote the real projective space of dimension $n$ and assume that $n>4$.

\m Recall that $\jtop: \tpl(\RR\PP^n) \to [\RR\PP^n, TOP/PL]$ is a 
bijection.  Consider a homeomorphism $k: M \to \RR\PP^5$ such that   
$$ 
\jtop(k)\ne 0\in [\RR\PP^n, TOP/PL] =H^3(\RR\PP^n;\ZZ/2)=\ZZ/2. 
$$ 
Note that  the Bockstein homomorphism
$$ 
\gb: \ZZ/2=H^3(\RR\PP^n;\ZZ/2) \to H^4(\RR\PP^n;\ZZ/2)=\ZZ/2 
$$ 
is an isomorphism, and hence $\gd(\jtop(k))\ne 0$ for $\gd:H^3(\RR\PP^n;\ZZ/2) \to H^4(\RR\PP^n)$.  
So, by Theorem~\ref{gd-ess}, $a_*\jtop(k)\ne 0$.   
In view of commutativity of the diagram \eqref{corr3}, $\jf(k)=a_*\jtop(k)$,  
i.e. $\jf(k)\ne 0$. 

On the other hand, it  follows from the obstruction theory that every homotopy equivalence 
$h: \RR\PP^n \to \RR\PP^n$ is homotopic to the identity map. In partcular, $\jf(h)=0$.
Thus, $M$ is not PL homeomorphic to $\RR\PP^n$. 
 
\begin{example}\label{non-conc} 
For every $n>3$ there is a homeomorphism 
$$
h=h_n: S^3 \times S^n \to S^3\times S^n, n>3
$$
which is homotopic 
to a PL homeomorphism but is not concordant to any PL homeomorphism.  
\end{example} 
 
\m Take an arbitrary homeomorphism  $f: M \to S^3\times S^n, n>3$. Then $\jf(f)$ 
is trivial by \theoref{ni}. Thus, by \theoref{bns}, $f$ is homotopic to a PL 
homeomorphism. In particular, $M$ is PL homeomorphic to $S^3\times S^n$.   
 
\m Now, we refine the situation and take a homeomorphism 
\[
h: S^3\times S^n \to S^3\times S^n
\]
such that   
$$ 
\jtop(h) \ne 0\in \tpl (S^3\times S^n)=H^3(S^3 \times S^n; \ZZ/2)=\ZZ/2. 
$$ 
Such $h$ exists because $\jtop$ is a bijection. So, $h$ is not concordant to the identity map, and therefore $h$ is not concordant to a PL homeomorphism, see \remref{conc}(2). But, as we have already seen, $h$ is homotopic 
to a PL homeomorphism.  
 
\m Note that the maps $h$ and the identity map have the same domain while 
they are not concordant. So, this example serves also the \remref{conc}(3). 
 
\begin{examples}\label{non-triang} 
There are topological manifolds that do not admit any PL structure. 
\end{examples}    

See manifold $V \times T^n$ that are described in \corref{c:F}. 

\m In 1970 Siebenmann~\cite{Sieb} published a paper with the intriguing title: Are nontriangulable manifolds triangulable? The paper cerebrated about the following problem: Are there manifolds that can be triangulated as simplicial complexes but do not admit any PL structure? Later, people made a big progress related to this issue, \cite{AM, GaSt1, GaSt2, Mat, Man, Rand,  Sav}. Here we give only a brief survey notice on the issue because non-combinatorial triangulations are far from the main line of our concern. 

\m Recall that a homology $k$-sphere is defined to be a $k$-dimensional closed PL  manifold $\gS$ such that $H_*(\gS)\cong H_*(S^k)$.  
 
\begin{examples}\label{tri-nonpl} 
There are topological manifolds that can be triangulated as simplicial complexes but do not 
admit any PL structure.  
\end{examples}  
 
\begin{theorem}\label{t:triang}
Every orientable topological $5$-dimensional closed manifold can be triangulated as a simplicial complex.
\end{theorem}

\p Let us say that a homology 3-sphere $\Sigma$ is good if the double suspension $S^2\Sigma$ over $\Sigma$ is homeomorphic to $S^5$ and $\Sigma$ is bounded by a compact parallelizable manifold of signature 8.  Siebenmann~\cite[Assertion on p. 81]{Sieb} proved that every orientable topological $5$-dimensional closed manifold can be triangulated as a simplicial complex provided that there exists a good homology 3-sphere. Cannon~\cite{Cannon} proved that, for any homology 3-sphere $\Sigma$, the double suspension $S^2\Sigma$ is homeomorphic to $S^5$. Now, note that the homology 3-sphere $\pa W$ from~\comref{plumb} is good, and the theorem follows.
\qed
 
\m Take $M=V \times S^1$. Then $M$ does not admit any PL structure by \corref{c:F}. On the other 
hand, $M$ can be triangulated as a simplicial complex by \theoref{t:triang}. Because of this, for each $k\geqslant 1$ the manifold $V\times T^k$ also have these properties. This is remarkable that $V$ cannot be triangulated as a simplicial complex, see below.

\begin{examples}\label{non-tri} 
There are topological manifolds that cannot be triangulated as simplicial complexes.
\end{examples}  

First, note that if a 4-dimensional topological manifold $M$ can be triangulated as a simplicial complex then  $M$ admits a PL structure. In paticular, $V$ cannot be triangulated as a simplicial complex, (Casson), see~\cite{AM, Sav}. 

\m Now we pass to higher dimensions.

\m Define two oriented homology 3-spheres $\gS_1, \gS_2$  to be {\it equivalent} if there exists an oriented PL bordism $W$, $\pa W=\gS_1\sqcup \gS_2$ such that  $H_1(W)=0=H_2(W)$. Let $\Theta_3^H$ denote the abelian group obtained from the set of equivalence classes using the operation of connected sum. We define a homomorphism  $\mu: \Theta_3^H\to \ZZ/2$ as follows.

\m It is well known that every homology 3-sphere (in fact, every orientable 3-manifold) $\Sigma$ bounds a 4-dimensional parallelizable manifold $P$. By the Rokhlin \theoref{rohlin}, the signature $\gs(P)\mod 16$ is a is well-defined invariant of $\gS$, and 8 divides $\gs(\gS)$. Take $a\in \Theta_3^H$, let $\gS_a$ be a homology 3-sphere that represents $a$, and let $P_a$ be a  4-dimensional parallelizable manifold with $\pa P_a=\gS_a$. Now, put $\mu(a) =(\gs(P_a) \mod 16)/8$ and  get a well-defined homomorphism $\mu: \Theta_H^3\to \ZZ/2$.  

Consider the short exact sequence
\[
\CD
0@>>> \ker\mu @>\subset>> \Theta_3^H @>\mu>> \ZZ/2 @>>> 0
\endCD
\]
and let $\delta:H^4(-;\ZZ/2)\to H^{5}(-; \ker \mu)$ be the Bockstein homomorphims associated with this sequence.

\begin{thm}[Galewski--Stern~\cite{GaSt2}, Matumoto~\cite{Mat}]\label{t:crit} 
A topological manifold $M$ of dimensional $\geqslant 5$ can be triangulated as a simplicial complex if and only if $\delta\varkappa (M)=0$. Here $\varkappa(M)$ denotes  the Kirby-Siebenmann invariant of $M$.
\end{thm}

\m Manolescu~\cite{Man} proved that the above mentioned short exact sequance does not split. This allowed him to prove that, for any $n\geqslant 5$, there is a manifold $M^n$ with  $\delta\varkappa (M^n)\neq 0$. Thus, for all $n\geqslant 5$ there exists an $n$-dimensional manifold that cannot be triangulated as a simplicial complex.

Concerning explicit constructions of such manifolds. Galewski amd Stern~\cite{GaSt1} constructed a certain manifold $N^5$ with the following property: if $N$ can be triangulated as a simplicial complex then {\em every} closed manifold of dimension $\geqslant 5$ can. So, $N$ cannot be triangulated. In particular,  $\gd\varkappa(N)\neq 0$. Finally, $N\times T^k$ cannot be triangulated as a simplicial complex because $\gd\varkappa(N \times T^k)\ne 0$.

\m {\bf \centerline {Summary}}

\m Here all manifolds are assumed to be connected and having the homotopy type of a finite CW complex.

\m {\sl $1$. Every manifold $M^n$ with $n\leqslant 3$ admits a unique PL structure $($trivial assertion for $n=1$, Rado~\cite{Rad} for $n=2$, Moise~\cite{Mo} for $n=3.)$

\m $2$. There are uncountable set of mutually different PL manifolds that are homeomorphic to $\RR^4$ $($Taubes~\cite{Ta}, cf also~\cite{GS, K2}$)$. There are countably infinite set of mutually different closed 4-dimensional PL manifolds that are homeomorphic to the blow-up of $\CC\PP^2$ at the nine points of intersection of two general cubics $($Okonek--Van~de~Ven~\cite{OV}$)$.

\m $3$. For every $n\geqslant 5$ there exist closed $n$-dimensional PL manifolds that are homeomorphic but not PL homeomorphic. So, the Hauptvermutung is wrong in general. However, any topological manifold $M^n, n\geqslant 5$\, $($not necessarily closed$)$ possesses only finite number of PL structures $($Kirby--Siebenmann~\cite{KS2}$)$. 

\m $4$. For every $n\geqslant 4$ there exist closed topological $n$-dimensional manifolds that do not admit any PL structure $($Freedman~\cite{F} for $n=4$, Kirby--Siebenmann~\cite{KS2} for $n>4)$.    

\m {\rm  The item 4 can be bifurcated as follows:}
 
 \m $4a$. For every $n\geqslant 5$ there exists an $n$-manifold that does not possess any PL structure but can be triangulated as a simplicial complex $($Siebenmann~\cite{Sieb} {\rm +} Cannon\cite{Cannon}$)$. Such examples do not exist for $n\leqslant 4$.
 
 \m $4b$. For every $n\geqslant 4$ there exists an $n$-manifold that cannot be triangulated as a simplicial complex 
 $($Casson~\cite{AM, Sav} for $n=4$, Manolescu~\cite{Man} for $n\geqslant 5.)$
 }

\section[Invariance of Characteristic  Classes ]{Topological and Homotopy  Invariance of Characteristic  Classes}

Given a real vector bundle $\xi$ over a space $X$, the $k$th {\it Pontryagin class} of $\xi$ is a cohomology class  $p_k(\xi)\in H^{4k}(X)$,~\cite{MS}. In particular, for every smooth manifold $M$ we have the Pontryagin classes $p_k(M):=p_k(\tau M)$ where $\tau M$ is the tangent bundle of $M$.  Given a commutative ring $\Lambda$ with unit, we can consider $p_k(\xi)\in H^{4k}(X;\Lambda)$, the image of the Pontryagin class $p_k(\xi)\in H^{4k}(X)$ under the coefficient homomorphism $\ZZ\to \Lambda$.
In particular, we have {\it rational Pontryagin classes} $p_k(\xi)\in H^{4k}(X;\QQ)$ and {\it modulo $p$ Pontryagin classes} $p_k(\xi)\in H^{4k}(X;\ZZ/p)$. 

In this section we discuss homotopy and topological invariance of  some characteristic classes. In particular, we prove that the Novikov's Theorem~\cite{N2} on topological invariance of rational Pontryagin classes is a direct corollary of the Main Theorem. (It is worthy to note, however, that the proof of the Main Theorem uses ideas from~\cite{N2}.) Concerning other proofs of  the Novikov's theorem see~\cite{G, ST, RW}.

 \begin{df}\rm\label{d:topinvar}
 Given a class  $x\in H^*(BO;\Lambda)$, we say that $x$ is {\it  topologically invariant} if, for any two maps $f_1, f_2: B\to BO$ such that
\[
 \ga^{O}_{TOP}f_1\cong  \ga^{O}_{TOP}f_2: B \to BTOP,
 \]
   we have
 \[
f_1^*(x)= f_2^*(x) \text{ in } H^*(B; \Lambda).
 \]
 \end{df}
 
\m  Now we give some conditions for topological invariance. Similarly to the fibration \eqref{albet}, consider the fibration 
\[
\CD
TOP/O @>\gb>> BO @>\ga>> BTOP.
\endCD
\] 

 \begin{prop}\label{p:nessuf} {\rm (i) } If 
 \[
 x\in \IM\{\ga^*: H^*(BTOP; \Lambda)\to H^*(BO; \Lambda)\}
 \]
  then $x$ is topologically invariant. In particular, if $\Lambda$ is such that $\ga^*$ is epimorphic than every class $x\in H^*(BO;\Lambda)$ is topologically invariant.
 
 {\rm (ii)} If $x\in H^*(BO;\Lambda)$ is topologically invariant then $\gb^*(x)=0$ for $\gb^*: H^*(BO)\to H^*(TOP/O)$. 
 \end{prop}

 \p (i) is obvious. To prove (ii), note that $\ga\gb$ is inessential. Hence $\ga\gb\cong\ga\eps$ where $\eps: TOP/O\to BO$ is a constant map. Since $x$ is  topologically invariant, we conclude that $\gb^*(x)=\eps^*(x)=0$. 
 \qed
 
\m \propref{p:nessuf}(i) tells us a sufficient condition for topological invariance, while \ref{p:nessuf}(ii) tells us a necessary condition. We will see below that \ref{p:nessuf}(i) is not necessary and \ref{p:nessuf}(ii) is not sufficient for topological invariance. Now we give a necessary and sufficient condition for invariance. Consider the map
\[
\CD
\mu: BO\times TOP/O @>1\times \beta >> BO\times BO @>m>> BO
\endCD
\]
where $m$ is the multiplication in the $H$-space $BO$.

 \begin{thm}\label{t:shar} The class $x\in H^*(BO;\Lambda)$ is topologically invariant if and only if $\mu^*(x)=x\otimes 1\in H^*(BO;\Lambda)\otimes H^*(TOP/O;\Lambda)$
\end{thm}

\p The map $\ga\gb$ is topologically trivial, and hence $\ga\mu$ is homotopic to the map 
\[\ga\nu:  BO\times TOP/O \to BO\to BTOP
\]
 where $\nu: BO\times TOP/O \to BO$ is the projection on the first factor. Since $x$ is topologically invariant, we conclude that $\mu^*x=\nu^*(x)=x\otimes 1$.

Conversely, suppose that $\mu^*(x)=x\otimes 1$. Recall that, for all $X$,  the infinite space structure in $BO$ turns $[X, BO]$ into an abelian group. Let $f_1, f_2: B\to BO$ be two maps such that $\ga f_1\cong \ga f_2$. Recall that $[X, BO]$ is an abelian group with respect to the infinite space structure in $BO$. Then $f_2-f_1: B\to BO$ lifts to a map $B \to TOP/O$. In other words, $f_2=f_1+g$ for some $g: B  \to TOP/O$. Hence we have a homotopy commutative diagram
\[
\CD
B@>f_2>>BO\\
@V\Delta VV @AA\mu A\\
B\times  B @>f_1\times g>> BO\times TOP/O
\endCD
\]
Now
\begin{equation*}
\begin{aligned}
f_2^*(x)&=\Delta^*(f_1\times g)^*\mu^*(x)=\Delta^*(f_1\times g)^*(x\otimes 1) 
=\Delta^*(f_1^*(x)\otimes 1)\\
&=f_1^*(x).
\end{aligned}
\end{equation*}
\qed

\begin{remark}\rm
 The items \ref{d:topinvar}--\ref{t:shar} are taken from the paper of Sharma~\cite{S}.
\end{remark}

\m The following lemma plays a  crucial role for topological invariance of rational Pontryagin classes. 
 
\begin{lemma}\label{l:bo-btop} 
The forgetful map $\ga^{O}_{TOP}: BO[0]\to BTOP[0]$ is a homotopy equivalence. Thus, the forgetful map $\ga^{O}_{TOP}: BO \to 
BTOP$ induces an isomorphism 
$$ 
(\ga^O_{TOP})^*: H^*(BTOP;\QQ) \to H^*(BO;\QQ). 
$$ 
\end{lemma} 
 
\p First, note that the homotopy groups $\pi_i(PL/O)$ are finite, see \cite[IV.4.27(iv)]{Rud} for the references. Hence, the space $PL/O[0]$ is contractible. Thus, $\ga^{O}_{PL}: BO[0]\to BPL[0]$ is a homotopy equivalence.

Second, the homotopy groups $\pi_i(TOP/PL)$ are finite by the Main Theorem. Hence, the space $TOP/PL[0]$ is contractible. Thus, $\ga^{PL}_{TOP}: BPL[0]\to BTOP[0]$ is a homotopy equivalence. 

Now, since $\ga^{O}_{TOP}=\ga_{TOP}^{PL}\ga^{O}_{PL}$, we conclude that $\ga^{O}_{TOP}[0]$ is a homotopy equivalence.
\qed 
 
\m Recall that $H^*(BO;\QQ)=\QQ[p_1, \ldots, p_i, \ldots]$ where $p_k, \dim  p_k=4k$ is the universal Pontryagin class,~\cite{MS}. It follows from \lemref{l:bo-btop} that $H^*(BTOP;\QQ)=\QQ[p'_1, \ldots, p'_k, 
\ldots]$ where $p'_k$ are the cohomology classes determined by the condition 
$$
\ga^*(p'_k)=p_k\in H^*(BO;\QQ).
$$ 
Now, given an arbitrary topological $\RR^n$ bundle $\lambda$ over $B$, we define its rational 
Pontryagin classes $p'_k(\gl)\in H^{4i}(B;\QQ)$ by setting  
$$
 p'_k(\gl) = t^*p'_k
 $$
  where $t: B \to BTOP$ classifies $\gl$.
  
\begin{thm}\label{t:topiso}
 Every class in $H^*(BO;\QQ)$ is topologically invariant. In other words, if $\xi_i=\{\pi_i: E_i\to B\}, i=1,2$ be two topologically isomorphic vector bundles over a space $B$ then $p_k(\xi_1;\QQ)=p_k(\xi_2;\QQ)$.
\end{thm}

This is the famous Novikov theorem on topological invariance of rational Pontryagin classes.

\p This follows from \lemref{l:bo-btop} and \propref{p:nessuf}(ii) immediately. \qed

 \m For completeness, we state the original Novikov version of topological invariance, see~\cite{N2}.
 
\begin{thm}\label{top-pontr}
Let $f: M_1 \to M_2$ be a  homeomorphism of closed smooth manifolds. and let $f^*:H^*(M_2;\QQ) \to H^*(M_1;\QQ)$ be the induced isomorphism. Then  $f^*p_k(M_2;\QQ)=p_k(M_1;\QQ)$ for all $k$. 
\end{thm}

\p Let $t_s: M_s\to BO\to BTOP, s=1,2$ classify the stable tangent bundle of $M_s$. Then  $t_1\simeq t_2f$. Now
\begin{equation*}
\begin{aligned}
f^*p_k(M_2;\QQ)=f^*t_2^*p'_k=(t_2f)^*p'_k=t_1^*p'_k=p_k(M_1;\QQ),
\end{aligned}
\end{equation*}
and we are done. 
\qed
 
 \begin{remark}\rm \label{r:plcat}
 We can pass the previous issues to PL category. To define PL invariance, we should replace topological isomorphism by PL isomorphism of PL bundles and require $B$ to be a polyhedron in \defref{d:topinvar}. Rokhlin and \v Svarc~\cite{RS} and Thom~\cite{T} proved PL invariance of rational Pontryagin classes in 1957-58th. Of course, this result follows from  the Novikov \theoref{top-pontr} on topological invariance of rational Pontryagin classes, but the Novikov Theorem appeared almost 10 years later.
 \end{remark}

So, rational Pontryagin classes are topological invariants. What about {\it integral} Pontryagin classes? It turns out to be that they are not even PL invariant. Milnor~\cite[\S 9]{Mi3} constructed two smooth manifolds $M_1, M_2$ that are PL homeomorphic while $p_2(M_1)=0, p_2(M_2)\neq 0$ (and $7p_2(M_2)=0$). 
\m Nevertheless, there are certain topological invariance results for integral Pontryagin classes.

\begin{notation}\rm \label{not:eps} Because of \lemref{l:bo-btop}, the index of the image subgroup
\[
\IM\{(\ga^O_{TOP})^*:H^m(BTOP)\to H^m(BO)\}
\]
in $H^m(BO)$ is finite for each $m$. Let $\eps_k$ denote this index for $m=4k$. Clearly, the class $\eps_kp_k\in H^{4k}(BO)$ (the multiple of the integral Pontryagin class) is topologically invariant. 
 
Define  $e_k\in\NN$ to be the smallest number such that $e_kp_k$ is topologically invariant. 
\end{notation}
 
\begin{comment}\label{f:calc}\rm
 To evaluate $e_k$, Sharma~\cite{S} proved the following. Let $d_k$ be the smallest positive integer such that 
\[
d_kp_k\in\Ker \{\gb*: H^*(BO) \to H^*(TOP/O)\}
\]
Then $e_k={\rm LCM}(d_1, \ldots, d_k)$. In particular, $e_k|e_{k+1}$.

To compute $d_k$, let $\gga_k=(2^{2k-1}-1)\Num (B_{2k}/4k)$. Here $B_m$'s are the Bernoulli numbers in notation where $B_{2n+1}=0$ and Num denotes the numerator.  Now, if $p$ is an odd prime which divides $\gga_k$ but does not divide $\gga_i$ with $i<k$, then $\nu_p(d_k)=\nu_p(\gga_k)$. Here, as usual, $m=p^{\nu_p(m)}a$ with $(a,p)=1$.
\end{comment}

\m Sharma~\cite[Theorem 1.6]{S} used these results in order to  evaluate $e_k$ for  $k\leqslant 8$. In particular, $e_1=1$, $e_2=7$, $e_3=7\cdot 31$, $e_4=7\cdot 31\cdot 127$. It is remarkable to note that $e_4< \eps_4$ (strict inequality!),~\cite[Prop. 1.7 ff]{S}. So, there are topologically invariant classes that do not come from $BTOP$, i.e. the sufficient condition \ref{p:nessuf}(i) for topological invariamce is not necessary. To see that the necessary condition \ref{p:nessuf}(ii) is not sufficient, note that $31p_3$ is not topological invariant because $e_3$ does not divide 31, while $31p_3\in\Ker \gb$, see~\cite[Section 4, proof of Theorem 1.3]{S}.
 
\m Another kind of topological invariance appears when we consider $p_k$ mod $m$, the modulo $m$ Pontryagin classes. Here we will not give detailed proofs but give a sketch/survey only.  As a first example, note that $p_k\mod 2=w_{2k}^2$, and hence $p_k\mod 2$ is topologically (and even homotopy) invariant in view of homotopy invariance of Stiefel--Whitney classes,~\cite{MS}. So, the question about topological invariance of modulo $p$ Pontryagin classes is not vacuous. In fact, we have the following result:
  
\begin{theorem}[SS]\label{t:modp}
Given an odd prime $p$, let $n(p)$ be the smallest value of $k$ such that $p$ divides $e_k$. Then $p_k \mod p$ is a topological invariant for $k < n(p)$ and is not a topological invariant for $k \geqslant n(p)$. In particular, if $p$ does not divide $e_k$, for every $k \geqslant 1$, then $p_k \mod p$ is a topological invariant. 
\end{theorem}

\m Because of \theoref{t:modp} and Comment~\ref{f:calc}, one can prove that the classes $p_k$ mod $p$ are topologically invariant for all $k$ and $p=3,5,11,13, 17$, while $p_k mod 7$ is not a topological invariant.  (For $p=3$ it is an old theorem of Wu, see~\theoref{wu}.) 

\m Now some words about homotopy invariance. 

\begin{df}\rm\label{d:hinvar}
Given a class  $x\in H^*(BO;\Lambda)$, we say that $x$ is {\it  homotopy invariant} if,  for any two maps $f_1, f_2: B\to BO$ such that
\[
 \ga^{O}_{F}f_1\cong  \ga^{O}_{F}f_2: B \to BF,
 \]
  we have
 \[
f_1^*(x)-r_2^*(x) \text{ in } H^*(B; \Lambda).
 \]
 \end{df}
 
  \m The obvious analogs of \propref{p:nessuf} and \theoref{t:shar} remains valid if we speak about homotopy invariance instead of topological invariance and replace $TOP$ by $F$.

\begin{prop}\label{p:notinvar}
Rational Pontryagin classes are {\em not} homotopy invariant.
\end{prop}

\p Note that $\pi_i(BF)$ is isomorphic to the  stable homotopy group $\pi_{i-1}^S(S)$ and therefore is finite because of a well-known theorem of Serre, ~\cite{Se}. Hence, $\pi_i(BF)\otimes \QQ=0$, and so $BF[0]$ is contractible. Now consider the fibration $F/O \to BO \to BF$ and conclude that $\gb[0]: F/O[0]\to BO[0]$ is a homotopy equivalence, and hence $\gb^*: H^*(BO;\QQ) \to H^*(F/O;\QQ)$ is an isomorphism. Thus, because of the homotopy analog of \ref{p:nessuf}(ii),  we see that $x\in \wt H^*(BO;\QQ)$ is homotopy invariant iff $x=0$.
\qed

\m On the other hand, we have $p_i\mod 2 =w_{2i}^2$,~\cite{MS} i.e. $p_i\mod 2$ is a homotopy invariant. So, it seems reasonable to ask about homotopy invariance of $p_i\mod p$ for odd prime $p$.

\m Recall that the homotopy invariance of Stiefel-Whitney follows from the Thom-Wu formula $w_i(\xi)=\gf^{-1} Sq^iu$ where $u$ is the Thom class of $\xi$ and $\gf:  H^*(B;\ZZ/2)\to \wt H^{*+n}(T\xi;\ZZ/2)$ is the Thom isomorphism,~\cite{MS}. (Here $T\xi$ is the Thom space of the $\RR^n$-bundle $\xi$ over $B$.) 

\m We apply this idea modulo $p$. So, let $p$ be an odd prime and 
$\C P^k: H^*(;\ZZ/p)\to H^{*+2k(p-1)}(;\ZZ/p)$ be the Steenrod power. Given an oriented $\RR^n$- bundl (or an $(S^n,*)$-fibration) $\xi$ over $B$, let $T\xi$ be the Thom space of $\xi$,  let $u\in H^n(T\xi;\ZZ/p)$ be the Thom class, and let $\gf:H^*(B ;\ZZ/p)\to \wt H^{*+n}(T\xi;\ZZ/p$ be the Thom isomorphism. Then
\[
q_k(\xi):=\gf^{-1}\C P^k(u)\in H^{2k(p-1)}(X)
\]
is a characteristic class, and it is homotopy invariant by construction. For $X=BO$ we get a universal characteristic class, and it is a polynomial of universal Pontryagin classes mod $p$. Wu~\cite{Wu} proved that $q_k=p_k$ if $p=3$. So, we get the following theorem.

\begin{thm}[Wu]\label{wu}
The Pontryagin classes $p_k\mod 3, k\geqslant 1$ are homotopy invariant.
\end{thm}

Madsen~\cite{M} proved that the classes $p_k\mod 8, k\geqslant 1$ are homotopy invariant. So, we have the following result:

\begin{cory}\label{c:24}
The classes $p_k\mod 24, k\geqslant 1$ are homotopy invariant. 
\end{cory}

 \end{document}